\documentstyle{amsppt}
\magnification=1200
\hoffset=-0.5pc
\vsize=57.2truepc
\hsize=38truepc
\nologo
\spaceskip=.5em plus.25em minus.20em
\define\parti{d}
\define\ltime{\ltimes}
\define\barakont{1}
\define\canhawei{2}
\define\cheveile{3}
\define\evluwein{4}
\define\gersthtw{5}
\define\geschthr{6}
\define\getzltwo{7}
\define\golmilfo{8}
\define\grothtwo{9}
\define\poiscoho{10}
\define\duality{11}
\define\bv{12}
\define\extensta{13}
\define\crosspro{14}
\define\liribi{15}
\define\huebstas{16}
\define\kosmathr{17}
\define\kosmafou{18}
\define\kosmafiv{19}
\define\kosmagtw{20}
\define\koszulon{21}
\define\liazutwo{22}
\define\liwextwo{23}
\define\liweixu{24}
\define\luweinst{25}
\define\mackfift{26}
\define\mackxu{27}
\define\maclaboo{28}
\define\majidtwo{29}
\define\maninfiv{30}
\define\mokritwo{31}
\define\palaione{32}
\define\rinehone{33}
\define\schecone{34}
\define\tianone{35}
\define\todortwo{36}
\define\weinsfte{37}
\define\wittetwe{38}
\define\xuone{40}
\define\Bobb{\Bbb}
\define\fra{\frak}

\topmatter
\title
Twilled Lie-Rinehart algebras and\\
differential Batalin-Vilkovisky algebras 
\endtitle
\author
Johannes Huebschmann
\endauthor
\date{November 2, 1998}
\enddate
\abstract 
Twilled L(ie-)R(inehart)-algebras generalize, in the Lie-Rinehart context, 
complex structures on smooth manifolds. An almost complex manifold determines 
an \lq\lq almost twilled pre-LR algebra\rq\rq, which is a true twilled 
LR-algebra iff the almost complex structure is integrable. We characterize 
twilled LR structures in terms of certain associated differential (bi)graded 
Lie and G(erstenhaber)-algebras; in particular the G-algebra arising from an 
almost complex structure is a d(ifferential) G-algebra iff the almost complex 
structure is integrable. Such G-algebras, endowed with a generator turning them
into a B(atalin-)V(ilkovisky)-algebra, occur on the B-side of the mirror 
conjecture. We generalize a result of Koszul to those dG-algebras which arise 
from twilled LR-algebras. A special case thereof explains the relationship 
between holomorphic volume forms and exact generators for the corresponding 
dG-algebra and thus yields in particular a conceptual proof of the Tian-Todorov
lemma. We give a differential homological algebra interpretation for twilled 
LR-algebras and by means of it we elucidate the notion of generator in terms of
homological duality for differential graded LR-algebras and we indicate how 
some of our results might be globalized by means of Lie groupoids.
\endabstract
\affil 
Universit\'e des Sciences et Technologies de Lille
\\
UFR de Math\'ematiques
\\
F-59 655 VILLENEUVE D'ASCQ C\'edex/France
\\
Johannes.Huebschmann\@univ-lille1.fr
\endaffil
\keywords 
{Lie-Rinehart algebra, twilled Lie-Rinehart algebra,
Lie bialgebra, Gerstenhaber algebra, 
Batalin-Vilkovisky algebra, 
differential graded Lie algebra,
mirror conjecture,
Calabi-Yau manifold}
\endkeywords
\subjclass{primary 17B55, 17B56, 17B65, 17B66, 17B70, 
secondary 
17B81, 32C16, 32G05, 53C05, 53C15, 81T70}
\endsubjclass
\endtopmatter
\document
\rightheadtext{Twilled Lie-Rinehart and BV algebras}

\medskip\noindent{\bf Introduction}\smallskip\noindent
A version of the mirror conjecture involves  
certain differential Batalin-Vilkovisky algebras
arising from a Calabi-Yau manifold.
A crucial ingredient is what is referred to
in the literature as the
{\it Tian-Todorov\/} lemma.
Our goal is to study
such differential Batalin-Vilkovisky
algebras and generalizations thereof
in the framework of Lie-Rinehart algebras.
Now a differential Batalin-Vilkovisky algebra 
is a Gerstenhaber algebra together with 
an exact generator, and the underlying Gerstenhaber algebras
of interest for us, in turn, arise 
as (bigraded) algebras of forms on twilled Lie-Rinehart algebras
(which we introduce below).
A twilled Lie-Rinehart algebra generalizes,
in the Lie-Rinehart context,
the notion of a complex structure
on a smooth manifold.
One of our results 
(Theorem 4.4)
will say that
 an \lq\lq almost
twilled pre-Lie-Rinehart algebra\rq\rq\ 
is a true twilled Lie-Rinehart algebra
if and only if  
the corresponding
Gerstenhaber algebra is a differential Gerstenhaber algebra.
(The wording is somewhat imprecise here
and Theorem 4.4 will in fact be phrased in terms of
\lq\lq almost
twilled Lie-Rinehart algebras\rq\rq, to be introduced below.)
As a consequence, we deduce that
an almost complex structure on a smooth manifold 
is integrable if and only if
the corresponding Gerstenhaber algebra is a differential
Gerstenhaber algebra.
Now a theorem of Koszul 
\cite\koszulon\ 
establishes, on an ordinary smooth manifold,
a bijective correspondence between generators for 
the Gerstenhaber algebra of multi vector fields
and connections in the top exterior power of the tangent bundle
in such a way that
exact generators correspond to 
flat connections.
In Theorem 5.4.6 below
we will generalize this bijective correspondence
to the differential Gerstenhaber
algebras 
arising from twilled Lie-Rinehart algebras;
such Gerstenhaber algebras come into play,
for example, in the mirror conjecture.
What corresponds to a flat connection
on the line bundle in Koszul's theorem is now
a holomorphic volume form---its existence
is implied by the Calabi-Yau condition---and 
our generalization of Koszul's
theorem shows in particular how
a holomorphic volume form determines
a generator for the corresponding
differential Gerstenhaber algebra
turning it into a differential Batalin-Vilkovisky algebra.
The resulting 
differential Batalin-Vilkovisky algebra
then generalizes that which underlies what is called the B-model.
In particular,
as a consequence of our methods,
we obtain a new proof of the Tian-Todorov lemma.
We will also give a 
differential homological algebra
interpretation of twilled Lie-Rinehart algebras
and, furthermore,
of a generator for a 
differential Batalin-Vilkovisky algebra
in terms of a suitable notion of homological duality.
Finally we indicate how some of our results might be globalized
by means of Lie groupoids.
\smallskip
We now give a more detailed outline of the paper.
Let $R$ be a commutative ring.
A Lie-Rinehart algebra $(A,L)$ 
consists of a commutative $R$-algebra $A$ and an $R$-Lie algebra
$L$ together with 
an $A$-module structure
$A \otimes_R L \to L$ 
on $L$, written
$a \otimes_R\alpha \mapsto a\alpha$,
and an action
$L \to \roman{Der}(A)$ 
of $L$ on $A$ (which is a morphism of $R$-Lie algebras and)
whose adjoint 
$L \otimes_R A \to A$
is written
$\alpha\otimes_Ra \mapsto \alpha(a)$;
here
$a \in A$ and $\alpha \in L$.
These mutual actions 
are required
to satisfy certain compatibility properties modeled on
$(A,L) = (C^{\infty}(M),\roman{Vect}(M))$
where
$C^{\infty}(M)$ and $\roman{Vect}(M)$
refer to the algebra of smooth functions 
and the Lie algebra of smooth vector fields,
respectively,
on a smooth manifold $M$.
In general, the compatibility conditions read:
$$
\align
(a \alpha) b &= a \alpha (b),
\quad a, b \in A,\ \alpha \in L,
\tag0.1
\\
[\alpha, a \beta] &= \alpha (a) \beta+ a [\alpha, \beta],
\quad a \in A,\ \alpha,\beta \in L.
\tag0.2
\endalign
$$
For a Lie-Rinehart algebra $(A,L)$,
following \cite\rinehone,
we will refer to $L$ as an $(R,A)$-{\it Lie algebra\/}.
In differential geometry, $(R,A)$-Lie algebras arise as spaces of sections
of Lie algebroids.
\smallskip
Given two Lie-Rinehart algebras $(A,L')$ and $(A,L'')$, together with mutual
actions
$\cdot\colon L' \otimes_R L'' \to L''$
and
$\cdot\colon L'' \otimes_R L' \to L'$
which endow $L''$ and $L'$ with an
$(A,L')$- and
$(A,L'')$-module structure, respectively,
we will refer to
$(A,L',L'')$ as an {\it almost twilled Lie-Rinehart algebra\/};
we will call it a
{\it twilled Lie-Rinehart algebra\/}
provided the direct sum $A$-module structure on $L=L' \oplus L''$,
the sum $(L' \oplus L'')\otimes_R A \to A$
of the adjoints of the
$L'$- and $L''$-actions on $A$,
and the bracket $[\cdot,\cdot]$
on $L=L' \oplus L''$ given by
$$
[(\alpha'',\alpha'),(\beta'',\beta')]
=
[\alpha'',\beta''] +[\alpha',\beta']
+ \alpha'' \cdot \beta'
-\beta' \cdot \alpha''
+ \alpha' \cdot \beta''
-\beta'' \cdot \alpha'
\tag 0.3
$$
turn $(A,L)$
into a Lie-Rinehart algebra.
We then write
$L=L' \bowtie L''$ and refer to
$(A,L)$ as the
{\it twilled sum\/}
of
$(A,L')$ and
$(A,L'')$.
\smallskip
For illustration, consider
a smooth manifold $M$ with an almost complex structure,
let $A$ be the algebra of smooth complex functions on $M$,
$L$ the $(\Bobb C,A)$-Lie algebra of complexified smooth vector fields
on $M$, and consider the ordinary decomposition
of the complexified tangent bundle $\tau^{\Bobb C}_M$ 
as a direct sum 
$\tau'_M \oplus \tau''_M$ 
of the {\it almost holomorphic\/}
and {\it almost antiholomorphic\/}
tangent bundles
$\tau'_M$
and  $\tau''_M$, respectively;
write $L'$ and $L''$ for their spaces of smooth sections.
Then 
$(A,L',L'')$,
together with the mutual actions coming from $L$,
is a twilled Lie-Rinehart algebra
if and only if the almost complex structure is integrable,
i.~e. a true complex structure;
$\tau'_M$
and  $\tau''_M$
are then the ordinary {\it holomorphic\/}
and {\it antiholomorphic\/}
tangent bundles, respectively.
In Section 1 below we actually show that
the precise analogue of an almost complex structure
is what we will call
an {\it almost twilled pre-Lie-Rinehart algebra\/} structure.
A situation similar to that of a complex
structure on a smooth manifold
and giving rise to a twilled Lie-Rinehart algebra
arises from a smooth manifold with two
transverse foliations as well as from a Cauchy-Riemann structure;
see Section 8 below for some comments about Cauchy-Riemann structures.
Lie bialgebras
provide
another class of examples
of
twilled Lie-Rinehart algebras; 
Kosmann-Schwarzbach and F. Magri refer to these objects, or rather
to the corresponding twilled sum, as
{\it twilled extensions of Lie algebras\/} \cite\kosmagtw;
Lu and Weinstein 
call them
{\it double Lie algebras\/} \cite\luweinst;
and Majid
uses the terminology
{\it matched pairs\/}
of Lie algebras \cite\majidtwo.
Spaces of sections
of suitable pairs of Lie algebroids
with additional structure
lead to yet another class of 
examples of
twilled Lie-Rinehart algebras;
these have been studied in the literature
under the name
{\it matched pairs 
of Lie algebroids\/}
by Mackenzie \cite\mackfift\  and Mokri \cite\mokritwo.
\smallskip
An almost twilled Lie-Rinehart algebra $(A,L'',L')$
is a true 
twilled Lie-Rinehart algebra 
if and only if 
$(A,L'',L')$
satisfies three compatibility conditions,
spelled out in Proposition 1.7 below;
this proposition is merely an adaption of
earlier results in the literature
to our more general situation.
We then give another interpretation of the compatibility conditions
in terms of 
annihilation properties of the two operators
$d'$ and $d''$
which arise as formal extensions of the ordinary Lie-Rinehart
differentials
with respect to $L'$ and $L''$, respectively,
on the bigraded algebra
$\roman{Alt}_A^*(L'', \roman{Alt}_A^*(L',A))$
(but are not necessarily exact);
{\sl for the twilled Lie-Rinehart algebra
arising from the holomorphic and antiholomorphic
tangent bundles of a complex manifold,
the resulting differential
bigraded algebra\/}
$(\roman{Alt}_A^*(L'', \roman{Alt}_A^*(L',A)),\parti', \parti'')$
{\sl comes down to the ordinary Dolbeault complex\/}.
See Theorem 1.15 for details.
\smallskip
In the rest of the paper,
we show that 
other characterizations
of
twilled Lie-Rinehart algebras
explain 
the differential Batalin-Vilkovisky
algebras
mentioned before:
Let
$(A,L'',L')$ be an almost twilled Lie-Rinehart algebra
having $L'$ finitely generated and projective as an $A$-module.
Write $\Cal A'' = \roman{Alt}_A(L'',A)$
and
$\Cal L' = \roman{Alt}_A(L'',L')$.
Now
$\Cal A''$ is a graded commutative $A$-algebra
and, endowed with the Lie-Rinehart differential $d''$
(which corresponds to the
$(R,A)$-Lie algebra structure on $L''$),
$\Cal A''$ is a differential graded commutative $R$-algebra.
Moreover,
from the $(A,L'')$-module structure
on $L'$,
$\Cal L'$
inherits an obvious
differential graded
$\Cal A''$-module structure.
Furthermore,
the $(A,L')$-structure
on $L''$ induces an action
of $L'$ on
$\Cal A''$
by graded derivations.
The latter, in turn,
induces a graded $R$-Lie algebra
structure on
$\Cal L'$
and a pairing
$$
\Cal L'
\otimes
\Cal A''
@>>>
\Cal A''
$$
which turns
$(\Cal A'',\Cal L')$
into a {\it graded Lie-Rinehart algebra\/}
(in an obvious sense);
this is in fact the {\it graded crossed product\/}
Lie-Rinehart structure.
Section 2 below is devoted to differential graded Lie-Rinehart
algebras; the differential graded crossed product
Lie-Rinehart
algebra will be explained in (2.8) and (2.9) below.
Now, on
$\Cal L' = \roman{Alt}_A(L'',L')$
we have the Lie-Rinehart 
differential $d''$ which corresponds to the
$(R,A)$-Lie algebra structure on $L''$ and
the $(A,L'')$-module structure $L'$.
By symmetry, 
when $L''$ is finitely generated and projective as an $A$-module,
we have the same structure, with
$L'$ and $L''$ interchanged.
Now
Theorem 3.2 will say that {\it the statements\/}
(i), (ii) and (iii) below {\it are equivalent\/}:
(i) 
{\sl $(A,L'',L')$
is a true
twilled Lie-Rinehart algebra\/};
(ii) 
$(\Cal L',d'') = (\roman{Alt}_A(L'',L'),d'')$
{\sl is a differential graded $R$-Lie algebra\/};
(iii)
$(\Cal A'',\Cal L';d'')$ {\sl is a differential graded Lie-Rinehart
algebra\/}.
Thus, under these circumstances,
{\sl there is a bijective correspondence between
twilled Lie-Rinehart algebra and differential graded
Lie-Rinehart algebra structures.}
\smallskip
We note that, in this situation,
the Lie bracket on
$\Cal L'=\roman{Alt}_A(L'',L')$
does not just come down to the shuffle product of
forms on $L''$ and the Lie bracket on $L'$;
in fact, such a bracket would not even be well defined
since the Lie bracket of 
$L'$ 
is not $A$-linear, i.~e.
does not behave as a \lq\lq tensor\rq\rq.
When $(A,L',L'')$
is the twilled Lie-Rinehart algebra
arising from 
the holomorphic and antiholomorphic tangent bundles of
a smooth complex manifold $M$,
$(\Cal L',d'') = (\roman{Alt}_A(L'',L'),d'')$
is what is called the
{\it Kodaira-Spencer\/}
algebra in the literature;
it controls the infinitesimal deformations
of the complex structure on $M$.
The cohomology
$\roman H^*(L'',L')$
then inherits a graded Lie algebra structure
and the obstruction to deforming the complex structure
is the map
$\roman H^1(L'',L')\to \roman H^2(L'',L')$
which sends $\eta\in \roman H^1(L'',L')$ to 
$[\eta,\eta]\in \roman H^2(L'',L')$.
\smallskip
We now return  to a general almost twilled Lie-Rinehart algebra
$(A,L',L'')$ having $L'$ finitely generated and projective
as an $A$-module and
consider the graded crossed product
Lie-Rinehart algebra
$(\Cal A'',\Cal L')$.
Write $\Lambda_AL'$ for the exterior $A$-algebra
on $L'$; as in the ungraded situation,
the 
graded Lie-Rinehart bracket
on
$\Cal L' (=\roman{Alt}_A(L'',L'))$
extends to a (bigraded) bracket
on
$\roman{Alt}_A(L'',\Lambda_AL')$
which turns the latter into a
{\it bigraded Gerstenhaber algebra\/};
as a bigraded algebra,
$\roman{Alt}_A(L'',\Lambda_AL')$
could be thought as of the exterior
$\Cal A''$-algebra on
$\Cal L'$, and we will often write
$$
\Lambda_{\Cal A''} \Cal L' = \roman{Alt}_A(L'',\Lambda_AL').
$$
The Lie-Rinehart differential $d''$
which corresponds to the
Lie-Rinehart structure on $L''$
and the induced
graded $(A,L'')$-module structure
on $\Lambda_AL'$
turn
$\roman{Alt}_A(L'',\Lambda_AL')$
into a differential (bi)-graded commutative $R$-algebra.
By symmetry, 
when $L''$ is finitely generated and projective as an $A$-module,
we have the same structure, with
$L'$ and $L''$ interchanged.
Theorem 4.4 will say that
{\sl the almost twilled Lie-Rinehart algebra $(A,L'',L')$
is a true
twilled Lie-Rinehart algebra
if and only if
$(\Lambda_{\Cal A''}\Cal L',d'')$
($=(\roman{Alt}_A(L'',\Lambda_AL'),d'')$)
is a differential (bi)-graded Gerstenhaber algebra\/}.
\smallskip
When $(A,L',L'')$ 
arises from 
the holomorphic and antiholomorphic tangent bundles of
a smooth complex manifold $M$,
the resulting differential Gerstenhaber algebra
$(\roman{Alt}_A(L'',\Lambda_AL'),d'')$
is that of forms of type
$(0,*)$
with values in the holomorphic multi vector fields,
the operator $d''$
being the Cauchy-Riemann operator
(which is more usually written $\overline \partial$).
This 
differential Gerstenhaber algebra
comes into play in the mirror conjecture;
it was studied by
Barannikov-Kontsevich \cite\barakont,
Manin \cite\maninfiv,
Witten \cite\wittetwe,
and others.
\smallskip
Let now $(A,L'',L')$ be a twilled Lie-Rinehart algebra
having $L'$
finitely generated and projective as an $A$-module
of constant rank $n$ (say),
and write
$\Lambda_A^nL'$ for the top exterior power of $L'$ over $A$.
Consider the differential Gerstenhaber algebra
$(\roman{Alt}_A(L'',\Lambda_AL'),d'')$.
Our next aim is to study
generators thereof.
To this end, we observe that,
when
$\roman{Alt}_A(L',\Lambda_A^nL')$
is endowed with the obvious graded
$(A,L'')$-module structure
induced from the left $(A,L'')$-module structure
on $L'$ which is part of the structure of twilled Lie-Rinehart algebra,
the canonical isomorphism
$$
\roman{Alt}_A(L'',\Lambda_AL')
@>>>
\roman{Alt}_A(L'',\roman{Alt}_A(L',\Lambda_A^nL'))
\tag0.4
$$
of graded $A$-modules
is compatible with the differentials
which correspond to the Lie-Rinehart structure
on $L''$ and the
$(A,L'')$-module structures on the coefficients
on both sides of (0.4);
abusing notation, we denote each of these differentials by $d''$.
Theorem 5.4.6 below says the following:
{\sl The isomorphism\/} (0.4) {\sl furnishes a bijective correspondence
between generators
of the bigraded Gerstenhaber structure
on the left-hand side\/} (of (0.4)) {\sl and
$(A,L')$-connections
on $\Lambda_A^nL'$
in such a way that
exact generators
correspond to
$(A,L')$-module structures
(i.~e. flat connections).
Under this correspondence,
generators
of the differential bigraded Gerstenhaber structure
on the left-hand side 
correspond to
$(A,L')$-connections
on $\Lambda_A^nL'$
which are compatible with the
$(A,L'')$-module structure on $\Lambda_A^nL'$.
Thus, in particular,
exact generators
of the differential bigraded Gerstenhaber structure
on the left-hand side 
correspond to
$(A,L'')$-compatible $(A,L')$-module structures
on $\Lambda_A^nL'$.\/}
\smallskip
When $L''$ is trivial and
$L'$ the Lie algebra of smooth vector fields
on a smooth manifold,
the statement of this theorem comes down to the result
of Koszul \cite\koszulon\ mentioned earlier.
Our result not only provides many examples 
of differential Batalin-Vilkovisky algebras
but also explains how
every differential Batalin-Vilkovisky algebra 
having an underlying bigraded $A$-algebra of the kind
$\roman{Alt}_A(L'',\Lambda_AL')$
arises.
\smallskip
When $(A,L',L'')$ is the twilled Lie-Rinehart algebra
which comes from the holomorphic and antiholomorphic tangent bundles
of a smooth complex manifold $M$ as explained above,
the theorem gives a bijective correspondence
between 
generators
of the 
differential bigraded Gerstenhaber algebra
$(\roman{Alt}_A(L'',\Lambda_AL'),d'')$
of forms of type
$(0,*)$
with values in the holomorphic multi vector fields,
the differential $d''$ being the
Cauchy-Riemann operator $\overline \partial$,
and holomorphic connections
on the highest exterior power of the holomorphic tangent bundle
in such a way that
exact generators correspond to flat holomorphic connections.
In particular, 
suppose
that $M$ is a {\it Calabi-Yau\/}
manifold,
that is,
admits
a holomorphic volume form $\Omega$ (say).
This holomorphic volume form
identifies
the highest exterior power of the holomorphic tangent bundle
with the algebra of 
smooth complex functions on $M$
as a module over $L=L''\oplus L'$,
hence induces
a flat holomorphic connection
thereupon 
and thence
an exact generator
$\partial_\Omega$
for
$(\roman{Alt}_A(L'',\Lambda_AL'),d'')$,
turning the latter into a differential (bi)graded
Batalin-Vilkovisky algebra.
This is precisely the
differential (bi)graded
Batalin-Vilkovisky algebra
coming into play 
on the B-side of the mirror conjecture
and studied in the cited sources.
The fact that the holomorphic volume form induces a generator
for the differential Gerstenhaber structure
is referred to in the literature as
the {\it Tian-Todorov lemma\/}.
In our approach, 
this lemma drops out
as a special case of our generalization
of Koszul's theorem to the bigraded setting,
and this generalization indeed provides a conceptual proof of the lemma.
This lemma implies that, for a K\"ahlerian Calabi-Yau manifold $M$,
the deformations of the complex structure are unobstructed, that is to say,
there is an open subset of
$\roman H^1(M,\tau_M)$
parametrizing the deformations of the complex structure;
here  $\roman H^1(M,\tau_M)$
is the first cohomology group
of $M$ with values in the holomorphic tangent bundle
$\tau_M$.
Under these circumstances,
after a choice of holomorphic volume form
$\Omega$
has been made,
the canonical isomorphism
(0.4), combined
with the isomorphism
$$
\Omega^{\flat}\colon
\roman{Alt}_A(L'',\roman{Alt}_A(L',\Lambda_A^nL'))
@>>>
\roman{Alt}_A(L'',\roman{Alt}_A(L',A))
$$
induced by $\Omega$
identifies
$(\roman{Alt}_A(L'',\Lambda_AL'),d'',\partial_\Omega)$
with the Dolbeault complex
of $M$ and hence
the cohomology
$\roman H^*(\roman{Alt}_A(L'',\Lambda_AL'),d'',\partial_\Omega)$
with the ordinary complex valued cohomology of $M$.
This is nowadays well understood;
see also 5.4.8 below.
The cohomology
$\roman H^*(\roman{Alt}_A(L'',\Lambda_AL'),d'',\partial_{\Omega})$
is referred to 
in the literature as the {\it extended moduli space of complex
structures\/} \cite\wittetwe;
it underlies what is called the B-{\it model\/}
in the theory of mirror symmetry.
\smallskip
In Section 6 we will give differential homological algebra
interpretations of some of our
earlier results.
In particular, we will show that,
{\sl for a twilled Lie-Rinehart algebra\/}
$(A,L',L'')$ {\sl having\/} $L'$ 
{\sl finitely generated and projective as an\/} $A$-{\sl module,
the differential
bigraded algebra\/}
$(\roman{Alt}_A^*(L'', \roman{Alt}_A^*(L',A)),\parti', \parti'')$
{\sl computes
the differential graded Lie-Rinehart cohomology\/}
$\roman H^*(\Cal L',\Cal A'')$,
where
$(\Cal A'',\Cal L';d'')$ 
is the differential graded
crossed product Lie-Rinehart algebra
$(\roman{Alt}_A(L'',A),\roman{Alt}_A(L'',L');d'')$
mentioned before.
When $L''$ is trivial,
so that 
$\roman H^*(\Cal L',\Cal A'')$
is an ordinary (ungraded) Lie-Rinehart algebra
$(A,L)$,
the differential graded Lie-Rinehart cohomology
boils down to 
the ordinary  Lie-Rinehart cohomology
$\roman H^*(L,A)$.
Moreover, 
for the special case 
when
$A$ and $L$ are the algebra of smooth functions and smooth vector
fields on a smooth manifold,
the Lie-Rinehart cohomology
$\roman H^*(L,A)$
amounts to the de Rham cohomology;
this fact has been established by Rinehart \cite\rinehone.
In our more general situation,
when
the twilled Lie-Rinehart algebra
$(A,L',L'')$ 
arises from the holomorphic and antiholomorphic tangent bundles
of a smooth complex manifold,
the complex
calculating
the differential graded Lie-Rinehart cohomology
$\roman H^*(\Cal L',\Cal A'')$
of the differential graded
crossed product Lie-Rinehart algebra
$(\Cal A'',\Cal L';d'')= (\roman{Alt}_A(L'',A),\roman{Alt}_A(L'',L');d'')$
is the Dolbeault complex,
and 
the differential graded Lie-Rinehart cohomology
amounts to the Dolbeault cohomology.
Thus our approach provides, in particular, 
an interpretation of the Dolbeault complex in the framework
of differential homological algebra.
In Section 7, 
generalizing results in our earlier paper \cite\bv, 
we will elucidate
the concept of generator of a differential bigraded
Batalin-Vilkovisky algebra
in the framework of homological duality for differential graded
Lie-Rinehart algebras.
In particular,
we will show that
{\sl an exact generator amounts to the differential
in a standard complex computing differential graded Lie-Rinehart
homology\/} (!) 
{\sl with appropriate coefficients\/};
see Proposition 7.13 below.
Further, we will see that,
when the appropriate additional structure 
(in terms of Lie-Rinehart differentials and dBV-generators)
is taken into account,
the above isomorphism (0.4) is essentially just a duality
isomorphism 
in the (co)homology of the differential graded crossed
product Lie-Rinehart algebra
$(\Cal A'',\Cal L')$;
see Proposition 7.14 for details.
In particular, the Tian-Todorov Lemma comes down to a statement
about differential graded
(co)homological duality.
\smallskip
Twilled Lie-Rinehart algebras thus generalize Lie bialgebras,
and the twilled sum
is an analogue,
even a generalization, of the Manin double of a Lie bialgebra.
The 
Lie bialgebroids
introduced by Mackenzie and Xu 
\cite\mackxu\ generalize Lie bialgebras as well,
and there is a corresponding notion of Lie-Rinehart bialgebra,
which we explain at the end of Section 4 below.
However,
twilled Lie-Rinehart algebras and Lie-Rinehart bialgebras
are different, in fact non-equivalent notions
which both generalize
Lie bialgebras.
In a sense,
Lie-Rinehart bialgebras
generalize Poisson 
and in particular symplectic
structures
while twilled Lie-Rinehart algebras
generalize complex structures.
In Theorem 4.8 below we characterize twilled Lie-Rinehart algebras
in terms of Lie-Rinehart bialgebras.
For the special case where the twilled Lie-Rinehart algebra
under 
consideration arises from a matched pair of Lie algebroids,
this characterization may be deduced from what is said in
\cite\mackfift.
In Section 8 below, we  use our characterization
of twilled Lie-Rinehart algebras in terms of Lie-Rinehart bialgebras
to indicate how some of the results
of this paper might be globalized
in terms of Lie groupoids.
\smallskip
As in Mac Lane's book
\cite\maclaboo,
(bi)graded objects
will always be understood as being externally (bi)graded.
\smallskip
I am indebted to Y. Kosmann-Schwarzbach and K. Mackenzie for
discussions, and to
J. Stasheff and A. Weinstein for some
e-mail correspondence
about various topics
related with the paper.
Most of the results to be given below
have been presented at the \lq\lq Poissonfest\rq\rq\ 
(Warsaw, August 1998), and at that occasion, Y. Kosmann-Schwarzbach
introduced me to the recent manuscript
\cite\schecone\ 
which treats topics somewhat related to
the present paper.
There is little overlap, though.
\medskip\noindent{\bf 1. Twilled Lie-Rinehart algebras}
\smallskip\noindent
The complexified tangent bundle $\tau^{\Bobb C}_M$ of a smooth
complex manifold $M$ decomposes
as a direct sum 
$\tau'_M \oplus \tau''_M$ 
of the {\it holomorphic\/}
and the {\it antiholomorphic\/}
tangent bundles
$\tau'_M$
and  $\tau''_M$, respectively;
both $\tau'_M$
and  $\tau''_M$ yield
smooth complex Lie algebroids over $M$,
and 
the integrability condition
amounts to 
these two Lie algebroid structures being compatible
in a very precise sense.
Similar situations arise from a smooth manifold with two
transverse foliations,
from Lie bialgebras, and
from a matched pair of Lie algebroids,
cf. what has been said in the introduction.
We now 
develop a theory 
incorporating,  generalizing,
and unifying these special cases.
\smallskip
As before, let $R$ be a commutative ring,
let $A$ be a commutative algebra, and let
$L'$ and $L''$ be two $A$-modules.
We will study and answer the following two related questions:
\smallskip
\noindent
{\smc Question 1.1.}
Given 
a Lie-Rinehart algebra $(A,L)$ and
a direct sum decomposition
$L= L' \oplus L''$ of $A$-modules
inducing
$(R,A)$-Lie algebra structures
on $L'$ and $L''$,
what kind of additional 
structure
relates $L'$ and $L''$?
The decomposition
$L= L' \oplus L''$ 
will then be referred to as an {\it integrable decomposition\/}
of $L$.
\newline\noindent
{\smc Question 1.2.}
Given $(R,A)$-Lie algebra structures
on $L'$ and $L''$,
which kind of additional structure
turns the direct sum
$L= L' \oplus L''$
of $A$-modules
into an $(R,A)$-Lie algebra
in such a way that
the action of $L$ on $A$ amounts to the sum of the
$L'$ and $L''$-actions
and that the bracket on $L$ restricts to the
given brackets on $L'$ and $L''$?
The new structure will then be referred to as
a {\it twilled \/} $(R,A)$-{\it Lie algebra\/} and 
the resulting 
$(R,A)$-Lie algebra
will be called
the {\it twilled sum\/}
of $L'$ and $L''$.
\smallskip\noindent
{\smc Example 1.3.}
Let $R$ be the ring $\Bobb C$ of complex numbers, $A$ the algebra of
smooth complex functions on a smooth almost complex manifold $M$,
$L$ be the $(\Bobb C, A)$-Lie algebra of smooth complexified
vector fields, and let
$L= L'' \oplus L'$
be the customary eigenspace decomposition 
(of the spaces of sections 
of the complexified tangent bundle)
arising from the almost complex structure.
The $A$-modules
$L'$ and $L''$ inherit
$(\Bobb C,A)$-Lie algebra structures
in such a way that
$L$ is their twilled sum 
(in a sense to be made precise)
if and only if
the almost complex structure is integrable.
\smallskip
We now proceed towards a description of a twilled
Lie-Rinehart algebra, the basic concept of the present paper:
Let $A$ be a commutative $R$-algebra.
Consider two $A$-modules $L'$ and $L''$,
together with skew-symmetric $R$-bilinear brackets
of the kind (1.4.1)---not necessarily Lie brackets---and
$R$-bilinear pairings
of the kind
(1.4.2$'$), (1.4.2$''$), (1.4.3), and (1.4.4) below:
$$
\alignat1
[\cdot,\cdot]'
&\colon L' \otimes_R L'
@>>>
L'
\tag1.4.1$'$
\\
[\cdot,\cdot]''
&\colon L'' \otimes_R L''
@>>>
L''
\tag1.4.1$''$
\\
L' \otimes_R A &@>>> A,
\quad x\otimes_R a \mapsto x(a),
\quad x\in L',\  a \in A
\tag1.4.2$'$
\\
L'' \otimes_R A &@>>> A,
\quad \xi\otimes_R a \mapsto x(a),
\quad \xi \in L'', a \in A
\tag1.4.2$''$
\\
\cdot\colon &L' \otimes_R L'' @>>> L''
\tag1.4.3
\\
\cdot\colon &L'' \otimes_R L' @>>> L'
\tag1.4.4
\endalignat
$$
We will refer to a $(A,L',L'')$
as an {\it almost twilled pre-Lie-Rinehart algebra\/},
provided 
$(A,L',L'')$
satisfies (i), (ii), and (iii) below.
\newline\noindent
(i) The values of
the adjoints
$L' @>>> \roman{End}_R(A)$
and
$L'' @>>> \roman{End}_R(A)$
of (1.4.2$'$) and (1.4.2$''$)
respectively
lie in 
$\roman{Der}_R(A)$;
\newline\noindent
(ii)
(1.4.1$'$), (1.4.2$'$)
and the $A$-module structure on $L'$
and, likewise,
(1.4.1$''$), (1.4.2$''$) 
and
the $A$-module structure on $L''$ 
are related by conditions of the kind (0.1)
and (0.2);
\newline\noindent
(iii) (1.4.3) and (1.4.4)
behave formally like connections.
\newline\noindent
Requirement (ii) is made precise by
(1.4.5$'$), 
(1.4.6$'$), 
(1.4.5$''$), 
(1.4.6$''$)  below,
and
(iii) is made precise by
(1.4.7) and (1.4.8).
$$
\alignat1
(a x) (b) &= a (x (b)),
\quad a, b \in A,\  x\in L',
\tag1.4.5$'$
\\
[x, a y]' &= x (a) y+ a [x, y]',
\quad a \in A,\ x,y \in L'
\tag1.4.6$'$
\\
(a \xi) (b) &= a (\xi (b)),
\quad a, b \in A,\  \xi\in L'',
\tag1.4.5$''$
\\
[\xi, a \eta]'' &= \xi (a) \eta+ a [\xi, \eta]'',
\quad a \in A,\ \xi,\eta \in L'',
\tag1.4.6$''$
\\
x \cdot (a\xi) &= (x(a)) \xi + a (x\cdot \xi),
\quad
(ax) \cdot \xi =  a(x\cdot \xi),
\quad x \in L',\ \xi \in L'',
\tag1.4.7
\\
\xi \cdot (a x) &= (\xi(a)) x + a (\xi\cdot x),
\quad
(a\xi) \cdot x =  a(\xi\cdot x),
\quad x \in L',\ \xi \in L''.
\tag1.4.8
\endalignat
$$
When $(A,L',L'')$ is
an almost twilled pre Lie-Rinehart algebra,
the pair $(L',L'')$, together with the other structure,
will be called an
{\it almost twilled pre\/}-$(R,A)$-{\it Lie algebra\/}.
\smallskip
Given an almost twilled pre-Lie-Rinehart algebra
$(A,L',L'')$,
let $L=L'\oplus L''$ be the direct sum of $A$-modules, and
extend the brackets on $L'$ and $L''$ to an
$R$-bilinear alternating
bracket
$$
[\cdot,\cdot]
\colon L \otimes_R L
@>>>
L
\tag1.5.1
$$
by means of 
the formula
$$
[(\alpha'',\alpha'),(\beta'',\beta')]
=
[\alpha'',\beta'']'' +[\alpha',\beta']'
+ \alpha'' \cdot \beta'
-\beta' \cdot \alpha''
+ \alpha' \cdot \beta''
-\beta'' \cdot \alpha',
\tag1.5.2
$$
and the two 
pairings
(1.4.2$'$) and (1.4.2$''$)
to a pairing
$$
L \otimes_R A @>>> A
\tag1.5.3
$$
in the obvious way,
that is,
by means of the assignment
$$
(\xi,x) \otimes_Ra \mapsto \xi(a) +x(a),
\quad x \in L',\, \xi \in L'',\, a \in A.
$$
By construction, the 
values of the
adjoint of (1.5.3) 
then lie in
$\roman{Der}_R(A)$, that is this adjoint
is then
of the form
$$
L = L'' \oplus L' @>>> \roman{Der}_R(A) .
\tag1.5.4
$$
An almost twilled pre-Lie-Rinehart algebra
$(A,L',L'')$
will be said to be 
an
{\it almost twilled Lie-Rinehart algebra\/}
provided
$(A,L')$, endowed
with the structure
(1.4.1$'$) and (1.4.2$'$),
and $(A,L'')$, endowed with (1.4.1$''$) and (1.4.2$''$),
are true Lie-Rinehart algebras and,
furthermore,
(1.4.3) is 
a 
left
$(A,L')$-module structure
on $L''$
and
(1.4.4)
a 
left
$(A,L'')$-module structure
on $L'$.
The pair $(L',L'')$, together with the two 
module structures
(1.4.3) and (1.4.4),
will then be called an
{\it almost twilled\/} $(R,A)$-{\it Lie algebra\/}.
An almost twilled Lie-Rinehart algebra
$(A,L',L'')$
will be said to be 
a {\it twilled Lie-Rinehart algebra\/}
provided
$(A,L)$, together
with the bracket (1.5.1) and 
the assignment
(1.5.4),
is a Lie-Rinehart algebra;
this Lie-Rinehart algebra will then be called
the {\it twilled sum\/}
of $(A,L')$ and
$(A,L'')$;
likewise,
$(L',L'')$
will then be called
a
{\it twilled\/} $(R,A)$-{\it Lie algebra\/}
and
$L$, written $L' \bowtie L''$, the {\it twilled sum\/}
of $L'$ and
$L''$.
\smallskip
A direct sum decomposition $L= L' \oplus L''$ 
of an $(R,A)$-Lie algebra $L$
yields in an obvious fashion
an almost twilled pre-Lie-Rinehart algebra structure
on $(A,L',L'')$:
The brackets
(1.4.1$'$) and
(1.4.1$''$)
result from restriction
and projection;
the pairings
(1.4.2$'$) and
(1.4.2$''$)
are as well obtained 
by restriction;
further,
the requisite pairings (1.4.3) and (1.4.4)
are given by the
composites
$$
\cdot\colon
L'\otimes_R L''
@>{[\cdot,\cdot]|_{L'\otimes_R L''}}>>
L'\oplus L''
@>>\roman{pr_{L''}}>
L''
\tag1.6.1
$$
and
$$
\cdot\colon
L''\otimes_R L'
@>{[\cdot,\cdot]|_{L''\otimes_R L'}}>>
L''\oplus L'
@>>\roman{pr_{L'}}>
L'
\tag1.6.2
$$
where,
for
$M= L'\otimes_R L''$ and
$M=L''\otimes_R L'$,
$[\cdot,\cdot]|_{M}$
denotes the restriction of the Lie bracket to $M$.
The formula (1.5.1) is then merely a decomposition of
the  initially given bracket on $L$
into components corresponding to the direct sum decomposition 
of $L$ into $L'$ and $L''$, and
(1.5.4) is accordingly a decomposition
of the $L$-action on $A$.
\smallskip
The following result is a mere adaption 
to our situation
of earlier results in the literature;
it is therefore labelled as a proposition.

\proclaim{Proposition 1.7}
An almost twilled Lie-Rinehart algebra
$(A,L',L'')$
is a true twilled Lie-Rinehart algebra if and only if
the Lie brackets
$[\cdot,\cdot]'$
and
$[\cdot,\cdot]''$
on $L'$ and
$L''$, respectively,
and the actions
{\rm (1.4.3)}
and
{\rm (1.4.4)}
are
related by
$$
\alignat1
\xi (x(a))- x(\xi(a)) &= (\xi \cdot x)(a) -(x \cdot \xi)(a) 
\tag1.7.1
\\
x\cdot[\xi,\eta]''
&=
[x\cdot\xi,\eta]''
+
[\xi,x\cdot \eta]''
-(\xi \cdot x)\cdot \eta
+(\eta \cdot x)\cdot \xi
\tag1.7.2
\\
\xi\cdot[x,y]'
&=
[\xi\cdot x,y]'
+
[x,\xi\cdot y]'
-(x \cdot \xi)\cdot y
+(y \cdot \xi)\cdot x ,
\tag1.7.3
\endalignat
$$
where $a\in A,\ x,y \in L',\ \xi,\eta \in L''$.
\endproclaim

An argument for the special case of this proposition
where
$L'$ and $L''$
are ordinary Lie algebras
may be found in
\cite\kosmagtw.
In fact, (1.7.2) and (1.7.3) then come down to
(1.3.1) and (1.3.2)
in \cite\kosmagtw.
More generally,
the case 
where
$L'$ and $L''$
arise from two Lie algebroids
has been established in Theorem 4.2 of \cite\mokritwo.

\demo{Proof}
The bracket (1.5.1) is plainly skew-symmetric.
Hence the proof comes down
to relating the Jacobi identity in $L$
and the Lie-Rinehart compatibility properties
with (1.7.1) -- (1.7.3).
\smallskip
Thus, suppose that the bracket
$[\cdot,\cdot]$ on $L=L'\oplus L''$ 
given by (1.5.1) satisfies the Jacobi identity.
Then, with a slight abuse of the notation $[\cdot,\cdot]$,
$$
\align
x\cdot [\xi,\eta]-[\xi,\eta] \cdot x 
&= [x,[\xi,\eta]]
\\
&= [[x,\xi],\eta] + [\xi,[x,\eta]]
\\
&= [x\cdot\xi -\xi \cdot x,\eta] + [\xi,x\cdot \eta - \eta \cdot x]
\\
&=  [x\cdot\xi,\eta] 
   -[\xi \cdot x,\eta]
   +[\xi,x\cdot \eta]
   -[\xi,\eta \cdot x]
\\
&= [x\cdot\xi,\eta] + [\xi,x\cdot \eta]
  -(\xi\cdot x)\cdot\eta -\eta\cdot (\xi \cdot x)
  +(\eta\cdot x)\cdot \xi - \xi\cdot(\eta \cdot x)
\endalign
$$
whence, comparing components in $L'$ and $L''$, we conclude
$$
\align
x\cdot [\xi,\eta]
&=
 [x\cdot\xi,\eta] + [\xi,x \cdot \eta]
  -(\xi\cdot x)\cdot\eta +(\eta\cdot x)\cdot \xi
\\
[\xi,\eta] \cdot x 
&=
 \xi\cdot(\eta \cdot x) -\eta\cdot (\xi \cdot x)
\endalign
$$
that is,
(1.7.2) 
holds and, furthermore,
$\cdot$ is a left $(A,L'')$-module structure
on $L'$ (but this is true already by assumption).
By symmetry,
(1.7.3) 
holds
as well.
Conversely,
suppose that
the two actions
are related by (1.7.2) and (1.7.3). 
We
can then read the above calculation backwards and conclude
that
the bracket $[\cdot,\cdot]$ on $L$ satisfies the Jacobi identity. 
\smallskip
We leave the rest of the proof to the reader. 
The arguments given in \cite\mokritwo\ are actually formal
and carry over. \qed
\enddemo

Theorem 1.7 thus gives a complete answer
 to Question 1.2,
as well as
to Question 1.1, as 
the following shows:

\proclaim{Corollary 1.8}
For an integrable decomposition
$L=L'\oplus L''$
of an $(R,A)$-Lie algebra $L$,
the resulting
almost
twilled pre-Lie-Rinehart algebra 
$(A,L',L'')$
is a true
twilled Lie-Rinehart algebra.
\endproclaim

They reader might ask: Why bother at all?
The answer is this:
We will show in
Section 5 below that the additional 
structure
relating the summands $L'$ and $L''$
of an integrable decomposition
explains in particular 
certain Batalin-Vilkovisky algebras
related with the mirror conjecture.
\smallskip \noindent
{\smc Remark 1.9.}
Let $\fra g$ be an ordinary Lie algebra,
finitely generated and projective over the ground ring $R$,
with Lie bracket $[\cdot,\cdot]$,
let $\Delta \colon \fra g \to \fra g \otimes_R \fra g$
be a Lie coalgebra structure on $\fra g$,
write
$[\cdot,\cdot]_*$
for the corresponding Lie bracket on
the dual $\fra g^*$,
and consider the pair $(\fra g, \fra g^*)$
together with the ordinary actions 
$\cdot\colon\fra g\otimes_R\fra g^* \to \fra g^*$
of
$\fra g$ on $\fra g^*$
and 
$\cdot \colon \fra g^*\otimes_R\fra g \to \fra g$
of
$\fra g^*$ on $\fra g$
induced by the Lie brackets on
$\fra g$ and $\fra g^*$, respectively,
the 
$\fra g$-  and $\fra g^*$-actions on the ground ring $R$
being taken trivial.
Then
(1.7.2) is equivalent to the customary requirement that
$\Delta$ be a 1-cocycle for $\fra g$ with values in
$\fra g \otimes_R \fra g$,
that is, to
$$
\Delta [x,y] = x \cdot \Delta y - y\cdot \Delta x
\tag1.9.1
$$
or, equivalently, to
$$
d_*[x,y] = [d_* x, y] + [x, d_* y]
\tag1.9.2
$$
where, on the right-hand side, $[\cdot,\cdot]$
refers to the corresponding Gerstenhaber bracket
on $\Lambda_R \fra g$;
here $d_*$ denotes the Chevalley-Eilenberg differential on
$\roman{Alt}_R(g^*,R)\cong \Lambda_R \fra g$.
Likewise,
(1.7.3) is equivalent to the requirement that the dual
$\Delta_* \colon \fra g^* \to \fra g^* \otimes_R \fra g^*$ 
of the Lie bracket $[\cdot,\cdot]$ on $\fra g$
be a 1-cocycle for $\fra g^*$ with values in
$\fra g^* \otimes_R \fra g^*$,
that is, to
$$
\Delta_* [\xi,\eta]_* = \xi \cdot \Delta_* \eta - \eta\cdot \Delta_* \xi
\tag1.9.3
$$
or, equivalently, to
$$
d[\xi,\eta]_* = [d \xi, \eta]_* + [\xi, d \eta]_*
\tag1.9.4
$$
where, on the right-hand side, $[\cdot,\cdot]_*$
refers to the corresponding Gerstenhaber bracket
on $\Lambda_R \fra g^*$, and
where
$d$ denotes the Chevalley-Eilenberg differential on
$\roman{Alt}_R(g,R)\cong \Lambda_R \fra g^*$.
Moreover,
(1.9.1) and (1.9.3)
are equivalent as well.
All these fact are nowadays well known.

\smallskip
\noindent
{\smc Example 1.10.}
An ordinary Lie bialgebra $(\fra g,\fra g^*)$
is as well a twilled $(R,R)$-Lie algebra,
as the corresponding Manin triple shows.
\smallskip
However, given a Lie-Rinehart algebra $(A,L)$
together with an
$(R,A)$-Lie algebra structure on 
$D = \roman{Hom}_A(L,A)$,
when the action of $L$ on $A$ 
(or that of $D$ on $A$, or that of both $L$ and $D$ on $A$) is non-trivial,
(1.9.1) and (1.9.3)
will not even make sense,
and a twilled Lie-Rinehart algebra structure 
on $(A,L,D)$ 
will {\it not\/} satisfy the obvious generalizations
of (1.9.2)
or (1.9.4).
In fact,
the obvious generalizations
of (1.9.2)
or (1.9.4)
lead to a different concept,
that of what we will call
a {\it Lie-Rinehart bialgebra\/};
see the end of Section 4 below and \cite\liribi.
Lie-Rinehart bialgebras generalize {\it Lie bialgebroids\/},
introduced in \cite\mackxu.

\smallskip
There is yet another way to understand
the integrability
of a decomposition of
a Lie-Rinehart algebra.
To explain it, we reproduce briefly the
Rinehart complex,
having as module variable a graded object:
Let $(A,L)$ be an (ungraded)
Lie-Rinehart algebra.
A graded $A$-module $M$,
together with a graded left $L$-module structure
$L \otimes_RM \to M$
is said to be a {\it graded\/}
({\it left\/})
$(A,L)$-{\it module\/}, provided
the actions are compatible, 
that is, for 
$\alpha
 \in  L, \,a \in A,\,m \in M$, we have
$$
\align
(a\,\alpha)(m) &= a(\alpha(m)),
\tag1.11.1
\\
\alpha(a\,m) &=   a\,\alpha(m) + \alpha(a)\,m.
\tag1.11.2
\endalign
$$
When $M$ is concentrated in degree zero,
we simply talk about a ({\it left\/}) $(A,L)$-{\it module\/}.
In particular, 
with the obvious structure,
the 
algebra $A$ itself is
a
(left) $(A,L)$-module.
Given  a graded $(A,L)$-module
$M$, the (bi)-graded
$R$-multilinear
alternating functions from $L$ into $M$
with the ordinary {\smc Cartan-Chevalley-Eilenberg\/}~\cite\cheveile\ 
differential $d$ given by 
$$
\aligned
(df)(\alpha_1,\dots,\alpha_n)
&=
\quad 
(-1)^n
\sum_{i=1}^n (-1)^{(i-1)}
\alpha_i(f (\alpha_1, \dots\widehat{\alpha_i}\dots, \alpha_n))
\\
&\phantom{=}+\quad
(-1)^n
\sum_{j<k} (-1)^{(j+k)}f(\lbrack \alpha_j,\alpha_k \rbrack,
\alpha_1, \dots\widehat{\alpha_j}\dots\widehat{\alpha_k}\dots,\alpha_n)
\endaligned
\tag1.11.3
$$
constitute a (graded) chain
complex
$\roman{Alt}_R(L,M)$
where
as usual \lq $\ \widehat {}\ $\rq \ 
indicates omission of the corresponding term.
As observed first by
{\smc Palais}~\cite\palaione\
(for the ungraded setting),
the defining properties
(0.1) and (0.2)
of a Lie-Rinehart algebra
entail that the differential $d$ on
$\roman{Alt}_R(L,M)$
passes to an
$R$-linear
differential on the
(bi)graded $A$-submodule
$\roman{Alt}_A(L,M)$
of $A$-multilinear functions,
written 
$$
d
\colon
\roman{Alt}_A(L,M)
@>>>
\roman{Alt}_A(L,M),
\tag1.11.4.1
$$
too,
and referred to henceforth as
{\it Lie-Rinehart\/} differential;
this differential will not be
$A$-linear unless $L$ acts trivially on $A$, though.
We will call
the resulting (co)chain complex
$$
(\roman{Alt}_A(L,M), d)
\tag1.11.4.2
$$
the {\it Rinehart complex of\/} $M$-{\it valued forms on\/} $L$;
often we write 
this complex
more simply in the form
$\roman{Alt}_A(L,M)$.
For $M=A$, with its obvious
$(A,L)$-module structure,
the differential $d$ turns
$\roman{Alt}_A(L,A)$
into a differential graded commutative $R$-algebra, and 
a general graded
$(A,L)$-module pairing $M_1 \otimes_AM_2 \to M$
induces a (bi)graded pairing
$$
\roman{Alt}_A(L,M_1) \otimes_R
\roman{Alt}_A(L,M_2)
@>>>
\roman{Alt}_A(L,M)
\tag1.11.5
$$
of $R$-chain complexes,
in fact of differential graded
$\roman{Alt}_A(L,A)$-modules.
The sign $(-1)^n$ in (1.11.3) has been introduced 
according to 
the customary Eilenberg-Koszul
convention in differential homological algebra,
since the Rinehart complex (1.11.4.2)
involves graded objects.
See also our paper \cite\extensta.
In the classical approach such a sign does not occur.
More generally, given a graded $A$-module $M$,
a graded pairing
$L\otimes_R M \to M$,
not necessarily a graded left $L$-module structure
but still satisfying (1.11.1) and (1.11.2),
is referred to as an 
$(A,L)$-{\it connection\/}, cf. \cite\poiscoho,
or, somewhat more precisely, as a
{\it graded left\/} $(A,L)$-{\it connection\/};
in this language,
an $(A,L)$-module structure
(or a graded one)
is a 
{\it flat\/} $(A,L)$-{\it connection\/}
(or a graded one).
Given a graded $A$-module $M$, together with an
$(A,L)$-connection,
we extend the definition of the
Lie-Rinehart operator to an operator
$$
d\colon
\roman{Alt}_A(L,M)
@>>>
\roman{Alt}_A(L,M)
\tag1.11.6
$$ 
by means of the formula (1.11.3),
with the $(A,L)$-connection
instead of the
$(A,L)$-action on $M$.
The resulting operator $d$ is well defined;
it is a differential if and only if the $(A,L)$-connection
on $M$ is flat, i.~e. a true $(A,L)$-module structure.
\smallskip

Let $(A,L',L'')$ be an almost twilled pre-Lie-Rinehart algebra.
Consider
the
bigraded
$A$-module
$$
\roman{Alt}_A^{*,*}(L'' \oplus L',A)
\cong
\roman{Alt}_A^*(L'', \roman{Alt}_A^*(L',A)).
\tag1.12.1
$$
Henceforth we spell out a particular homogeneous constituent as
$$
\roman{Alt}_A^q(L'', \roman{Alt}_A^p(L',A)),
\tag1.12.2
$$
keeping in mind
that, under the circumstances
of (1.3),
when the almost complex structure is a true complex structure,
the notations $p$ and $q$ have become standard for
the \lq\lq holomorphic\rq\rq\  and 
\lq\lq antiholomorphic\rq\rq\ degrees, respectively; 
for intelligibility,
we follow this
convention, see below.
The pairings (1.4.3) and (1.4.4) induce
graded pairings
$$
\align
L' \otimes_R\roman{Alt}_A^*(L'', A)
&@>>>
\roman{Alt}_A^*(L'',A)
\tag1.13.1
\\
L''\otimes_R \roman{Alt}_A^*(L',A)
&@>>>
\roman{Alt}_A^*(L',A)
\tag1.13.2
\endalign
$$
on
$\roman{Alt}_A^*(L'',A)$
and $\roman{Alt}_A^*(L',A)$,
respectively,
when (1.4.3) and (1.4.4) 
are formally treated like connections.
Via (1.11.6), applied formally,
that is, by a formal evaluation
of the expression given
on the right-hand side of (1.11.3),
with 
(1.4.1$'$)
and
(1.4.1$''$)
instead of the Lie brackets,
and
(1.4.2$'$)
and
(1.4.2$''$)
instead of the requisite
module structures,
these pairings, in turn,
induce two operators
$$
\align
\parti'
\colon
\roman{Alt}_A^q(L'', \roman{Alt}_A^p(L',A))
&@>>>
\roman{Alt}_A^q(L'', \roman{Alt}_A^{p+1}(L',A))
\tag1.14.1
\\
\parti''
\colon
\roman{Alt}_A^q(L'', \roman{Alt}_A^p(L',A))
&@>>>
\roman{Alt}_A^{q+1}(L'', \roman{Alt}_A^p(L',A)).
\tag1.14.2
\endalign
$$
A little thought reveals that,
in view of
(1.4.5$'$),
(1.4.5$''$),
(1.4.6$'$),
(1.4.6$''$),
(1.4.7),
(1.4.8),
these operators, which are at first
defined only on the $R$-multilinear alternating functions,
in fact pass to operators
on $A$-multilinear alternating functions.
Then the requirement that $d=d'+d''$ be a differential,
i.~e. that $dd = 0$, amounts to
$$
\align
\parti'\parti' &= 0
\tag1.15.1
\\
\parti''\parti'' &= 0
\tag1.15.2
\\
[\parti',\parti''] &= 0,
\tag1.15.3
\endalign
$$
where as usual
$[\parti',\parti''] =\parti'\parti''+\parti''\parti'$;
in other words,
$d$ being a differential
is equivalent to
$$
(\roman{Alt}_A^*(L'', \roman{Alt}_A^*(L',A)),
\parti', \parti'')
\tag1.15.4
$$
being a bicomplex.
\smallskip
An $A$-module $M$ will be said to have {\it property\/} P 
provided for  $x\in M$,
$\phi(x) = 0$ for every $\phi \colon M \to A$
implies that $x$ is zero.
For example, a projective $A$-module has property P,
or a reflexive $A$-module has this property
or, more generally, an $A$-module $M$ such that the canonical map
from $M$ into its double $A$-dual is injective.
On the other hand, for example, 
for a smooth manifold $X$,
the $C^{\infty}(X)$-module $D$ of formal (= K\"ahler) differentials
does {\it not\/} have property P:
On the real line, with coordinate $x$,
consider the functions $f(x) = \sin x$ and $g(x) = \cos x$.
The formal differential
$d f -g dx$ is non-zero in $D$;
however, the $C^{\infty}(X)$-linear maps from $D$ to
$C^{\infty}(X)$
are the smooth vector fields, whence every 
such $C^{\infty}(X)$-linear map
annihilates
the formal differential
$d f -g dx$.

\proclaim{Theorem 1.15}
If $(A,L',L'')$ is a twilled Lie-Rinehart algebra,
{\rm (1.15.4)} is a bicomplex which
then necessarily computes
the cohomology
$\roman H^*(\roman{Alt}_A(L,A))$
of the twilled sum $L$ of $L'$ and $L''$.
Conversely,
$(A,L',L'')$ being an almost twilled pre-Lie-Rinehart algebra,
if {\rm (1.15.4)} is a bicomplex,
and if $L'$ and $L''$ have property {\rm P},
$(A,L',L'')$ is a true twilled Lie-Rinehart algebra.
\endproclaim

\demo{Proof}
If $(A,L',L'')$ is twilled Lie-Rinehart algebra,
{\rm (1.15.4)} is plainly a bicomplex which
then necessarily computes the indicated cohomology.
We now prove the converse.
Thus suppose that
{\rm (1.15.4)} is a bicomplex.
Consider the operator
$$
d''d''
\colon
\roman{Alt}_A^j(L'', A)
@>>>
\roman{Alt}_A^{j+2}(L'', A)
$$
for $j=0$ and $j=1$.
Notice that
$\roman{Alt}_A^j(L'', A)$
equals
$\roman{Alt}_A^j(L'', \roman{Alt}_A^0(L',A))$
and that
$\roman{Alt}_A^{j+2}(L'', A)$
equals $\roman{Alt}_A^{j+2}(L'', \roman{Alt}_A^0(L',A))$.
For $j=1$,
given
$\xi,\eta,\vartheta \in L''$ and  
$\phi \in \roman{Hom}_A(L'',A)=\roman{Alt}_A^j(L'', A)$,
we find
$$
(d''d'' \phi)(\xi,\eta,\vartheta)
=
\phi(
[[\xi,\eta]'',\vartheta]'' +
[[\eta,\vartheta]'',\xi]'' +
[[\vartheta,\xi]'',\eta]''). 
$$
Since $L''$ has property P,
we conclude that the bracket on
$L''$
satisfies the Jacobi identity,
that is,
$L''$ is an $R$-Lie algebra.
Likewise,
for $j=0$,
given
$\xi,\eta\in L''$ and $a \in A$, we find
$$
(d''d'' a)(\xi,\eta)
=
\xi (\eta (a))-\eta(\xi (a)) - [\xi,\eta](a).
$$
Consequently
the adjoint
$L'' \to \roman{Der}_R(A)$
of (1.4.2$'')$
is a morphism of $R$-Lie algebras.
In view of (1.4.5$''$) and
(1.4.6$''$),
we conclude that $(A,L'')$
is a Lie-Rinehart algebra.
The same kind of reasoning shows that
$(A,L')$
is a Lie-Rinehart algebra.
\smallskip
Next,
consider the operator
$$
d''d''
\colon
\roman{Alt}_A^0(L'', \roman{Alt}_A^1(L',A))
@>>>
\roman{Alt}_A^2(L'', \roman{Alt}_A^1(L',A)).
$$
We note that
$\roman{Alt}_A^0(L'', \roman{Alt}_A^1(L',A))=
\roman{Alt}_A^1(L',A)=
\roman{Hom}_A(L',A)$.
Let $x\in L'$, $\xi,\eta \in L''$, and
$\phi \in \roman{Hom}_A(L',A)$.
A straightforward calculation gives
$$
(d''d'' \phi)(\xi,\eta)
=
\phi(\eta \cdot (\xi \cdot x) - \xi \cdot (\eta \cdot x)
+[\xi,\eta] \cdot x).
$$
Since $L'$ is assumed to have property P,
we conclude that,
for every $x\in L',\, \xi,\eta \in L''$,
$$
[\xi,\eta] \cdot x =  \xi \cdot (\eta \cdot x)-\eta \cdot (\xi \cdot x),
$$
that is,
(1.4.4)
is a left $(A,L'')$-module structure
on $L'$.
The same kind of reasoning shows that
(1.4.3)
is a left $(A,L')$-module structure
on $L''$.
\smallskip
Pursuing the
same kind of reasoning, consider the operator
$$
d'd''+d''d'
\colon
A=\roman{Alt}_A^0(L'', \roman{Alt}_A^0(L',A))
@>>>
\roman{Alt}_A^1(L'', \roman{Alt}_A^1(L',A)).
$$
Let $a \in A,\,x \in L',\,\xi \in L''$.
Again a calculation shows that
$$
((d'd''+d''d')a)(\xi,x)
=
\xi (x(a))- x(\xi(a)) - ((\xi \cdot x)(a) -(x \cdot \xi)(a))
$$
whence
the vanishing of
$d'd''+d''d'$
in bidegree (0,0)
entails the compatibility property (1.7.1).
Likewise
consider the operator
$$
d'd''+d''d'
\colon
\roman{Hom}_A(L'',A) =\roman{Alt}_A^1(L'', \roman{Alt}_A^0(L',A))
@>>>
\roman{Alt}_A^2(L'', \roman{Alt}_A^1(L',A)).
$$
Again a calculation shows that, for
 $x \in L',\,\xi,\eta \in L'',
\phi\in \roman{Hom}_A(L'',A)$,
$$
((d'd''+d''d')\phi)(\xi,x)
=
\phi\left(x\cdot[\xi,\eta]
-
([x\cdot\xi,\eta]
+
[\xi,x\cdot \eta]
-(\xi \cdot x)\cdot \eta
+(\eta \cdot x)\cdot \xi)\right)
$$
whence
the vanishing of
$d'd''+d''d'$
in bidegree (1,0)
entails the compatibility property (1.7.2).
Likewise,
the vanishing of
$d'd''+d''d'$
in bidegree (0,1)
entails the compatibility property (1.7.3). \qed
\enddemo

\smallskip\noindent
{\smc Example 1.3} (continued).
An almost complex structure determines
an almost twilled pre-Lie-Rinehart algebra
$(A,L',L'')$, and the almost complex structure is integrable
if and only if
$L= L''\oplus L'$ is an integrable decomposition
(whence the terminology).
In the integrable case,
the operator $\parti''$---the corresponding
Cauchy-Riemann operator---defines 
a holomorphic structure on the manifold $M$,
and the bigraded object
$$
\roman{Alt}_A^{*,*}(L''\oplus L',A) \cong
\roman{Alt}_A^*(L'',\roman{Alt}_A^*(L',A))
$$
amounts to the 
smooth complex valued forms
of type $(0,*)$
with values in the exterior powers of the holomorphic
cotangent bundle.
The resulting bicomplex
$\roman{Alt}_A^*(L'',(\roman{Alt}_A^*(L',A),\parti'),\parti'')$
is the customary
{\it Dolbeault\/} complex 
computing the  sheaf hypercohomology
of $M$ with values in the complex 
of sheaves of germs of holomorphic 
differential forms on $M$
and,
by virtue of the Poincar\'e lemma, 
the total complex of the 
Dolbeault complex 
yields a resolution of the constant sheaf
of complex numbers
whence 
the cohomology 
of the Dolbeault complex 
coincides with the ordinary
smooth complex valued (de Rham) cohomology of $M$, viewed
as a real manifold.
All this is classical, cf. e.~g. \cite\grothtwo.
Our description in terms of Lie-Rinehart structures 
seems to be new, though.
\smallskip
Thus, twilled Lie-Rinehart algebras
generalize complex manifolds
in the same sense as 
Lie bialgebroids
or more generally Lie-Rinehart bialgebras
(see Section 4 below or \cite\liribi)
generalize Poisson 
and in particular symplectic structures.
Almost twilled pre-Lie-Rinehart algebras
have been spelled out above
as the exact analogue 
of almost complex structures.
In the rest of the paper, 
almost twilled pre-Lie-Rinehart algebras will no longer
come into play explicitly and
only almost twilled 
and twilled Lie-Rinehart algebras
will be considered.

\medskip\noindent{\bf 2. Differential graded Lie-Rinehart algebras}
\smallskip\noindent
There are various concepts of
differential graded Lie algebras in the literature.
To introduce notation, we reproduce a description tailored to our purposes.
To simplify the exposition somewhat, we will assume that
the primes 2 and 3 are invertible in the ground ring $R$.
If $x$ is an element in a graded module then
$|x|$ denotes its degree.
\smallskip
Let 
$L$
be a chain complex over $R$,
and let
$$
\lbrack \cdot,\cdot\rbrack
\colon
L \otimes_R L
\longrightarrow
L
$$
be a 
pairing of 
chain complexes
of degree zero.
We will say that
$(L,\lbrack \cdot,\cdot\rbrack)$
is  a {\it  differential graded Lie algebra\/}
provided it is skew-symmetric in the graded sense
and satisfies the {\it graded Jacobi identity\/},
that is,
$$
\alignat 1
\lbrack x,y \rbrack  
&=
- (-1)^{|x||y|}\lbrack y,x \rbrack, \quad
\text{for all $x$ and $y$ in $L$},
\tag2.1.1
\\
\lbrack x,\lbrack y,z \rbrack\rbrack
&= 
\lbrack \lbrack x,y \rbrack,z\rbrack
  + (-1)^{|x||y|} 
\lbrack y,\lbrack x,z \rbrack \rbrack,
\quad \text{for all\ } x,\,y,\,z \in L.
\tag2.1.2
\endalignat
$$
Here are two immediate consequences 
of the definition:
$$
\alignat 1
\lbrack x,x \rbrack  
 &= 0,\quad
\text{for all homogeneous $x$ in $L$ of even degree},
\tag2.1.1.$\tfrac 12$
\\
\lbrack x,\lbrack x,x \rbrack\rbrack &=  0,
\quad \text{for all homogeneous $x$ in $L$ of odd degree.}
\tag2.1.2.$\tfrac 13$
\endalignat
$$
The pairing
$\lbrack \cdot,\cdot\rbrack$
is what is called a ({\it graded\/}) {\it Lie bracket\/}.
Given two  differential graded Lie algebras $L$ and $L'$,
a {\it morphism\/} $\phi \colon L \longrightarrow L'$
{\it of  differential graded Lie algebras\/} over $R$
is the obvious
thing, i.~e. it is a morphism of chain complexes
 which is compatible with the
graded Lie brackets.
We note that, when $2$ is not invertible in the ground ring,
there are two
non-equivalent notions of
graded Lie algebra depending on whether
or not
(2.1.1.$\tfrac 12$)
is required to hold and, likewise,
when $3$ is not invertible in the ground ring,
(2.1.2.$\tfrac 13$)
is an additional requirement.
\smallskip
For 
a 
 differential graded
algebra $U$
over $R$,
the {\it associated  differential graded Lie algebra\/}
over $R$,
written $LU$ or, with an abuse of notation, just $U$,
has the same underlying chain complex as $U$,
while its bracket $\lbrack \cdot,\cdot \rbrack$
is given by
$$
\lbrack u,v\rbrack  =
uv - (-1)^{|u||v|}vu,
\quad                                 
\text{for \ }  u, v \in U.
\tag2.2
$$
Whenever we say that a
 differential graded
 algebra is viewed as a 
 differential graded Lie algebra,
this structure will be the intended one.
In particular, 
for a  chain complex $M$, the object
$\roman{End}_R(M)$ is a  differential graded algebra, and hence 
$L(\roman{End}_R(M))$ is a  differential graded Lie algebra.
Furthermore, if $L$ is 
a  differential graded Lie algebra
and $M$ a chain complex,
a {\it  differential graded $L$-module structure on \/}
$M$ is a morphism
$L \longrightarrow
L(\roman{End}_R(M))$ of differential graded Lie algebras.
\smallskip
Let $U$ be a  differential graded algebra over $R$.
Recall that a {\it (homogeneous) derivation\/}  of $U$
is a (homogeneous)  morphism
$
\delta
\colon
U
\longrightarrow
U
$
of chain complexes
so that for $u,\,v \in U$,
$$
\delta(uv) = (\delta(u))v + (-1)^{|u||\delta|}u\delta(v).
\tag2.3
$$
The 
graded submodule
$\roman{Der}(U)$ 
of 
derivations of $U$
is a  graded submodule of
$\roman{End}_R(U)$;
moreover, it inherits a differential from the latter,
and it is well known that
the  bracket
(2.2)
induces a bracket
$$
\lbrack\cdot,\cdot\rbrack  
\colon
\roman{Der}(U)\otimes_R
\roman{Der}(U)
\longrightarrow
\roman{Der}(U)
\tag2.4
$$
for
$\roman{Der}(U)$ which
turns $\roman{Der}(U)$ 
into a  differential graded Lie algebra over $R$.
Further,
if $L$ is a  differential graded Lie algebra over $R$
and if $U$ is a differential graded $R$-algebra,
as usual, a morphism
$  
L 
\longrightarrow
\roman{Der}(U)
$
of  differential graded
Lie algebras over $R$
is called
an {\it action\/} of
$L$ on $U$ (by derivations);
on elements, we will always write
the adjoint
$L \otimes_RU \to U$
of an $L$-action on $U$
in the form
$\alpha \otimes_R x\mapsto \alpha(x), \,\alpha \in L,\, x \in U$.
\smallskip
Given 
two  differential graded algebras $U$ and $U'$
over $R$,
 differential graded 
Lie algebras $L$ and  $L'$ 
over $R$, and actions of
$L$ and $L'$ on $U$ and $U'$ respectively,
a {\it morphism \/}
$$
(\phi,\psi)
\colon
(U,L)
\longrightarrow
(U',L')
$$
({\it of actions\/})
is the
obvious thing, i.e. it consists of a morphism
$\phi \colon
U \longrightarrow U'
$
of  differential graded
$R$-algebras
and a morphism
$\psi \colon
L \longrightarrow L'
$
of 
 differential graded Lie algebras over $R$,
 so that the diagram
$$
\CD
L \otimes_R U 
@>>>
U
\\
@V{\psi \otimes_R \phi}VV
@V{\phi}VV
\\
L' \otimes_R U'
@>>>
U'
\endCD
$$
is commutative;
here the unlabelled horizontal arrows refer to the corresponding
structure maps.
\smallskip
Given  a differential graded Lie algebra $L$
and a chain complex $M$ over $R$,
as usual, a morphism
$  
L 
\longrightarrow
L\roman{End}(M)
$
of  differential graded
Lie algebras over $R$
is called an {\it action\/} of
$L$ on $M$,
and $M$ is said to be a
{\it  differential graded (left)\/} $L$-{\it module\/};
we will always write
the adjoint
$L \otimes_RM \to M$
in the form
$\alpha \otimes_R x\mapsto \alpha(x), \,\alpha \in L,\, x \in M$.
The precise definition of the concept of
a morphism of 
differential graded $L$-modules
is obvious and left to the reader.
\smallskip
We now generalize
the notion of Lie-Rinehart algebra
to that of differential graded Lie-Rinehart algebra.
For intelligibility,
ordinary (ungraded) Lie-Rinehart algebras will be denoted by
$(A,L)$ etc. and
differential graded Lie-Rinehart algebras by $(\Cal A,\Cal L)$, etc.
\smallskip
Let $\Cal A$ be a  differential graded commutative
$R$-algebra,
let ${\Cal L}$ be a  differential graded Lie algebra over $R$,
let 
$
\Cal A \otimes_R {\Cal L} \to {\Cal L}
$
be a   differential graded left $\Cal A$-module structure on ${\Cal L}$,
written 
${a\otimes_R \alpha \mapsto a\,\alpha}$, 
and 
let 
$
{\Cal L} 
\to
\roman{Der}(\Cal A)
$
be an
action of
${\Cal L}$ on $\Cal A$
whose adjoint 
${\Cal L} \otimes_R\Cal A \to \Cal A$
is written
$\alpha\otimes_R a \mapsto \alpha(a),\,\alpha \in {\Cal L},\, a \in \Cal A$. 
We will refer to ${\Cal L}$ as a {\it  differential graded \/}
$(R,\Cal A)$-{\it Lie algebra\/},
provided
$$
\alignat1
(a\,\alpha)(b) &= a\,(\alpha(b)),
\quad \alpha \in {\Cal L}, \,a,b \in \Cal A,
\tag2.5.a
\\
\lbrack 
\alpha,
a\,\beta
\rbrack
&=
(-1)^{|a||\alpha|}
a\,
\lbrack 
\alpha,
\beta
\rbrack
+
\alpha(a)\,\beta,\quad \alpha,\beta \in {\Cal L}, \,a \in \Cal A.
\tag2.5.b
\endalignat
$$
Extending terminology introduced in our paper
~\cite\poiscoho\ 
(for the ungraded case), we will refer to a pair
$(\Cal A,{\Cal L})$,
where $\Cal A$ is a differential graded commutative algebra and
${\Cal L}$ a differential graded
$(R,\Cal A)$-Lie algebra, 
as a {\it differential graded Lie-Rinehart algebra\/}.
An example of a
differential graded Lie-Rinehart algebra
is the pair $(\Cal A,\roman{Der}(\Cal A))$,
where $\Cal A$ is a differential graded commutative algebra
and 
$\roman{Der}(\Cal A)$ the differential graded $\Cal A$-module
of graded derivations of $\Cal A$, with the obvious structures.
\smallskip
Given 
two
differential graded Lie-Rinehart algebras
$(\Cal A,{\Cal L})$
and
$(\Cal A',{\Cal L}')$,
a {\it morphism \/}
$
(\phi,\psi)
\colon
(\Cal A,{\Cal L})
\longrightarrow
(\Cal A',{\Cal L}')
$
{\it of differential graded Lie-Rinehart algebras\/}
is the
obvious thing, that is, it 
is a morphism
of actions in the above sense so that,
in addition,
$\psi \colon
{\Cal L} \longrightarrow {\Cal L}'
$
is a morphism of  differential graded
left $\Cal A$-modules where
$\Cal A$ acts on ${\Cal L}'$ via $\phi$.
\smallskip
Let $(\Cal A,{\Cal L})$ be a differential graded
Lie-Rinehart algebra and
let $M$ be a
chain complex over $R$ having 
a differential graded left $\Cal A$-module
structure
and, furthermore, a differential graded
left ${\Cal L}$-module structure.
Then $M$ is said to be a
{\it differential graded\/} ({\it left\/})
 $(\Cal A,{\Cal L})$-{\it module\/}, provided
the actions are compatible, 
that is, for 
$\alpha
 \in {\Cal L}, \,a \in \Cal A,\,m \in M$, we have
$$
\align
(a\,\alpha)(m) &= a(\alpha(m)),
\tag2.6.a
\\
\alpha(a\,m) &=   (-1)^{|\alpha||a|}a\,\alpha(m) + \alpha(a)\,m.
\tag2.6.b
\endalign
$$
In particular, 
with the obvious structure,
the 
differential graded
algebra $\Cal A$ itself is
a
 differential graded
(left) $(\Cal A,{\Cal L})$-module.
Furthermore, there is an obvious notion of
 morphism of modules
over differential graded Lie-Rinehart
algebras;
we leave the details to the reader.

\smallskip
For a differential graded Lie algebra
$L$ over $R$,
given differential graded (left) $L$-modules
$M'$ and $M''$,
the customary formula 
$$
\alpha(x \otimes_R y) = \alpha(x) \otimes_R y +
(-1)^{|\alpha||x|} x \otimes_R \alpha(y),
\quad
\alpha \in L, 
\,x \in M',\, y \in M'',
\tag2.7
$$
endows
the 
differential graded
tensor product $M' \otimes_R M''$ with a 
differential graded
(left) $L$-module structure;
this is just the ordinary (differential graded) tensor
product $L$-module structure.
If
$M$
is another differential graded $L$-module, 
a pairing 
$
\mu
\colon
M'
\otimes_R
M''
\longrightarrow
M
$
of $R$-modules
which is a morphism of
differential graded
 $L$-modules
(with respect to (2.7)) will be said to be a
{\it a pairing of 
differential graded $L$-modules\/}.
For an ungraded Lie-Rinehart algebra $(A,L)$,
viewed as
a differential graded Lie-Rinehart algebra concentrated
in degree zero with zero differential,
given differential graded $(A,L)$-modules 
$M'$ and $M''$,
a little thought reveals that
the formula
(2.7) turns the (graded) tensor product
$M' \otimes_A M''$
into a differential graded 
$(A,L)$-module;
we refer to
$M' \otimes_A M''$
with this structure as the
{\it tensor product of\/}
$M'$ and $M''$
{\it in the category of
differential graded $(A,L)$-modules\/}.
Given 
differential graded
$(A,L)$-modules 
$M$,\ $M'$,\ and $M''$,
a pairing
${
\mu_A
\colon
M'
\otimes_A
M''
\longrightarrow
M
}$
of $A$-modules
 which is compatible with the 
differential graded
$L$-structures
will be said to be a 
{\it pairing of 
differential graded $(A,L)$-modules\/}.
See our paper \cite\extensta\  for more details.
\smallskip
\noindent
{\smc 2.8. The graded crossed product extension.}
For later reference, we reproduce
briefly a description of the graded crossed product Lie-Rinehart algebra
extension
tailored to our purposes;
see \cite\poiscoho\ for the ungraded case.
Let $(A,L)$ be a Lie-Rinehart algebra, and let
$\Cal A$ be a graded commutative 
$A$-algebra 
which is endowed with a graded $(A,L)$-module structure
in such a way that (i) $L$ acts
on $\Cal A$ by derivations---this is equivalent
to requiring the structure map
from $\Cal A \otimes_A \Cal A$ to $\Cal A$
to be a morphism of graded $(A,L)$-modules---and that (ii)
the canonical map from $A$ to $\Cal A$ is
a morphism of left $(A,L)$-modules.
Let 
$\Cal L = \Cal A \otimes_A L$, and define a bigraded bracket
$$
[\cdot,\cdot]
\colon
\Cal L
\otimes_R
\Cal L
@>>>
\Cal L
\tag2.8.1
$$
of bidegree $(0,-1)$
by means of the formula
$$
[\alpha \otimes_A x, \beta \otimes_A y]
=
(\alpha \beta) \otimes_A[x,y]
+\alpha (x\cdot \beta) \otimes_A y
- (-1)^{|\alpha||\beta|} \beta (y \cdot \alpha)\otimes_A x
\tag2.8.2
$$
where
$\alpha,\beta \in \Cal A$ and $x,y \in L$.
A calculation shows that,
for every $\beta \in \Cal A$ and every $x,y,z \in L$,
$$
[[x,y],\beta \otimes_Az] -
([x,[y,\beta \otimes_Az]] -[y,[x,\beta \otimes_A z]])
=\big([x,y] (\beta) - x(y(\beta))-y(x(\beta))\big )\otimes_A z,
$$
whence 
(2.8.1) being a graded Lie bracket is actually equivalent
to the structure map $L\otimes_R\Cal A \to \Cal A$ being a Lie algebra
action.
Moreover,
let
$$
\Cal A \otimes_R \Cal L @>>> \Cal L
\tag2.8.3
$$
be the obvious graded left
$\Cal A$-module structure
arising from extension of scalars,
that is from extending $L$ to a (graded)
$\Cal A$-module,
and
define a pairing
$$
\Cal L \otimes_R \Cal A @>>> \Cal A
\tag2.8.4
$$
by
$$
(\alpha \otimes_A x) \otimes_R \beta \mapsto
(\alpha \otimes_A x) (\beta)
=\alpha (x(\beta)).
\tag2.8.5
$$
Then $(\Cal A,\Cal L)$,
together with (2.8.1),
(2.8.3) and (2.8.4),
constitutes a graded Lie-Rinehart algebra.
We refer to 
$(\Cal A,\Cal L)$
as the {\it crossed product\/}
of  $\Cal A$ and $(A,L)$
and to the corresponding
$(R,\Cal A)$-Lie algebra $\Cal L$
as the {\it crossed product\/} of $\Cal A$ and $L$.
\smallskip
\noindent
{\smc Remark 2.8.6.} We must be a little circumspect here:
The three terms on the right-hand side of (2.8.2)
are {\it not\/} well defined individually; only their sum is well
defined. For example,
if we take $ax$ instead of $x$, where $a \in A$,
on the left-hand side,
$\alpha \otimes_A(ax)$ equals
$(\alpha a)\otimes_A x$
but
$(\alpha \beta) \otimes_A[ax,y]$
differs from
$(\alpha a\beta) \otimes_A[x,y]$.
\smallskip
\noindent
(2.9) A special case of the differential graded crossed product
arises as follows:
Let $(A,L)$ be a Lie-Rinehart algebra, let
$M$ be a left $(A,L)$-module,
and consider the graded $A$-algebra $\Cal A = \roman{Alt}_A(M,A)$,
endowed with the induced
left $(A,L)$-module structure;
this is in fact an $L$-action on
$\Cal A$ by derivations.
We then have the crossed product 
$(R,\Cal A)$-Lie algebra $\Cal L$
which, as a graded
$\Cal A$-module, has the form
$\Cal A \otimes_A L$.
When $L$ is finitely generated and projective,
the canonical morphism
$$
\Cal A \otimes_A L
@>>>
\roman{Alt}_A(M,L)
\tag2.9.1
$$
is an isomorphism, and (2.8.3) comes down to the ordinary shuffle
pairing
$$
\roman{Alt}_A(M,A)
\otimes_R
\roman{Alt}_A(M,L)
@>>>
\roman{Alt}_A(M,L).
\tag2.9.2
$$

\medskip\noindent{\bf 3. The integrability condition reexamined}
\smallskip\noindent
Let $(A,L',L'')$ be an almost twilled Lie-Rinehart algebra.
Suppose that, as an $A$-module, $L'$ is 
finitely generated and projective.
Let
$\Cal A'' =\roman{Alt}_A^*(L'',A)$ and,
with reference to the
left
$(A,L')$-module structure
(1.4.3) on $L''$,
consider the graded crossed product Lie-Rinehart algebra
$$
(\Cal A'',\Cal L')
=
(\roman{Alt}_A^*(L'',A),\roman{Alt}_A^*(L'',A) \otimes_A L')
\cong
(\roman{Alt}_A^*(L'',A),\roman{Alt}_A^*(L'',L'))
\tag3.1
$$
explained in (2.9),
$L''$ playing the role of $M$ in (2.9).
The left $(A,L'')$-module structure 
(1.4.4) on $L'$
induces the corresponding
Lie-Rinehart
operator
$\parti''$
(cf. 1.11.4.1)
on
$\Cal L'=\roman{Alt}_A^*(L'',L')$
and, with reference to the graded left
$\Cal A''$-module structure
(2.8.3)
on $\Cal L'$,
$\Cal L'$ is a differential graded $\Cal A''$-module,
$\Cal A''$ being endowed with the ordinary Lie-Rinehart differential
(1.11.4.1).

\proclaim{Theorem 3.2}
Under these circumstances, the following are equivalent.
\newline\noindent
{\rm (i)}
$(A,L',L'')$ is a twilled Lie-Rinehart algebra.
\newline\noindent
{\rm (ii)}
$(\Cal L',[\cdot,\cdot]',\parti'')$
is a differential graded
$R$-Lie algebra.
\newline\noindent
{\rm (iii)}
$(\Cal A'',\Cal L';\parti'')$
is a differential graded
Lie-Rinehart algebra.
\endproclaim

We note that,
when $L''$ is finitely generated and projective as an $A$-module, 
with the roles of $L'$ and $L''$ interchanged,
the same statements as those spelled out in Theorem 3.2 are true.
Under the circumstances of Theorem 3.2,
we will refer to
$(\Cal A'',\Cal L';\parti'')$
as a {\it differential graded
crossed product Lie-Rinehart algebra\/}.
Thus the theorem says that,
provided $L'$ is finitely generated and projective as an $A$-module,
twilled Lie-Rinehart algebras
and differential graded crossed product Lie-Rinehart algebras 
are equivalent notions.
\smallskip
The proof requires some preparation.

\proclaim{Lemma 3.3}
Given $x$ and $y$ in $L'$
$$
\parti''[x,y]'
=[\parti''x,y]' +[x,\parti''y]'
\tag3.3.1$'$
$$
if and only if
$$
\xi \cdot [x,y]'
=
[\xi \cdot x,y]' 
+[x,\xi \cdot y]' 
- (x \cdot \xi) \cdot y
+(y \cdot \xi) \cdot x
$$
for {\it every\/} $\xi \in L''$.
Consequently the truth of {\rm (3.3.1$'$)}
for every 
$x$ and $y$ in $L'$
is equivalent to
the compatibility condition
{\rm (1.7.3)}.
\endproclaim

\demo{Proof}
Let $x,y \in L'$ and write
$$
\align
\parti''x &= \sum \alpha_i \otimes_A x_i 
\in \roman{Alt}_A(L'',A) \otimes_A L'
\\
\parti''y &= \sum \alpha_j \otimes_A x_j 
\in \roman{Alt}_A(L'',A) \otimes_A L'.
\endalign
$$
Then
$$
\align
[\parti''x,y]' 
&= \sum [\alpha_i \otimes_A x_i ,y]'
=
\sum \alpha_i \otimes_A [x_i,y]' - (y \cdot \alpha_i) \otimes_A x_i
\\
[x,\parti''y]' 
&= \sum [x,\beta_j \otimes_Ay_j]'
=
\sum \beta_j \otimes_A [x,y_j]'
+ (x \cdot \beta_j) \otimes_A y_j
\endalign
$$
Thus, given
$\xi \in L''$, we have
$$
\align
[\parti''x,y]'(\xi)
&=
\sum \alpha_i (\xi)\otimes_A[x_i,y]'-((y \cdot \alpha_i)(\xi))\otimes_A x_i
\\
&=\sum \alpha_i (\xi)\otimes_A[x_i,y]'
-y(\alpha_i(\xi)) \otimes_A x_i
+ \alpha_i(y \cdot \xi)) \otimes_A x_i 
\\
&=
[(\parti''x)(\xi),y]' 
+ \sum (\alpha_i(y \cdot \xi)) \otimes_A x_i 
\\
&=
[(\parti''x)(\xi),y]' 
+ (\parti''x)(y \cdot \xi)
\\
&=
[\xi \cdot x,y]' 
+(y \cdot \xi) \cdot x
\\
[x,\parti''y]'(\xi) 
&=
\sum \beta_j(\xi) \otimes_A [x,y_j]'
+ ((x \cdot \beta_j)(\xi)) \otimes_A y_j
\\
&=\sum \beta_j(\xi) \otimes_A [x,y_j]'
+x(\beta_j(\xi)) \otimes_A y_j
- (\beta_j(x \cdot \xi)) \otimes_A y_i 
\\
&=
[x,(\parti''y)(\xi)]' 
- \sum (\beta_j(x \cdot \xi)) \otimes_A y_i 
\\
&=
[x,(\parti''y)(\xi)]' 
- (\parti''y)(x \cdot \xi)
\\
&=
[x,\xi \cdot y]' 
- (x \cdot \xi) \cdot y
\endalign
$$
Consequently
$$
\align
([\parti''x,y]' +[x,\parti''y]')(\xi)
&=
[\xi \cdot x,y]' 
+(y \cdot \xi) \cdot x
+[x,\xi \cdot y]' 
- (x \cdot \xi) \cdot y
\\
&=
[\xi \cdot x,y]' 
+[x,\xi \cdot y]' 
- (x \cdot \xi) \cdot y
+(y \cdot \xi) \cdot x,
\endalign
$$
On the other hand
$$
(\parti''[x,y]')(\xi) =\xi \cdot [x,y]'.
$$
Hence
$$
\align
(\parti''[x,y]'
&-[\parti''x,y]' -[x,\parti''y]')(\xi)
\\
&=
\xi \cdot [x,y]'
-
([\xi \cdot x,y]' 
+[x,\xi \cdot y]' 
- (x \cdot \xi) \cdot y
+(y \cdot \xi) \cdot x)
\endalign
$$
This completes the proof of the Lemma. \qed
\enddemo

For later reference, from the proof of Lemma 3.3,
we record the following formulas
$$
\align
[\parti''x,y]'(\xi)
&=  
[\xi \cdot x,y]' 
+(y \cdot \xi) \cdot x
\tag3.3.2$'$
\\
[x,\parti''y]'(\xi) &=
[x,\xi \cdot y]' 
- (x \cdot \xi) \cdot y,
\tag3.3.3$'$
\endalign
$$
where
$x,y \in L'$ and $\xi \in L''$.

\proclaim{Lemma 3.4}
Given $x \in L'$,
the following are equivalent.
\newline\noindent
{\rm (i)} For every
homogeneous $\alpha\in \Cal A''=\roman{Alt}_A(L'',A)$,
$$
d''(x \cdot \alpha) =
(d''x) \cdot \alpha + x \cdot (d''\alpha).
\tag3.4.1$'$
$$
\newline\noindent
{\rm (ii)} For every $a \in A$ and
every $\xi \in L''$,
$$
\xi(x(a))
-
x(\xi(a))
=(\xi\cdot x - x \cdot \xi)(a)
$$
and,
for every $\xi,\eta \in L''$, and 
every $\beta \in \roman{Alt}^1_A(L'',A)=\roman{Hom}_A(L'',A)$,
$$
\beta(x\cdot[\xi,\eta]
-
([x\cdot\xi,\eta]
+
[\xi,x\cdot \eta]
-(\xi \cdot x)\cdot \eta
+(\eta \cdot x)\cdot \xi))
=0
$$
Consequently  {\rm (3.4.1$'$)}
holds for every 
$x \in L'$
and every
$\alpha \in \roman{Alt}_A(L'',A)$
if and only if the
compatibility conditions
{\rm (1.7.1)}
and
{\rm (1.7.2)}
are satisfied.
\endproclaim

\demo{Proof}
Let $\xi \in L''$, and $a \in A$. Then
$$
\align
(d''(x(a)))(\xi) &=
\xi(x(a))
\\
((d''x)(a))(\xi) &=
((d''x)(\xi))(a) = (\xi \cdot x)(a)
\\
(x\cdot (d''a))(\xi) &=
x((d''a)(\xi))-(d''a)(x \cdot \xi)
=
x(\xi(a))-(x \cdot \xi)(a)
\endalign
$$
whence
$$
\align
(d''(x (a))
&-
(d''x) (a)
-
x \cdot(d''a))(\xi)
\\
&=
\xi(x(a))
-
x(\xi(a))
-(\xi\cdot x - x \cdot \xi)(a)
\endalign
$$
Consequently, given $a \in A$,
$$
d''(x (a))
=
(d''x) (a)
+
x \cdot(d''a)
$$
if and only if
$$
\xi(x(a))
-
x(\xi(a))
=(\xi\cdot x - x \cdot \xi)(a)
$$
for every $\xi \in L''$.
Thus the statement of the Lemma is true
when $\alpha$ has degree 0.
\smallskip
Let  $\xi,\eta \in L''$, and 
$\beta \in \roman{Alt}^1_A(L'',A)=\roman{Hom}_A(L'',A)$.
Then
$$
\align
(x \cdot \beta)(\xi) &= x (\beta (\xi))- \beta (x \cdot \xi)
\\
(d''(x \cdot \beta))(\xi,\eta)&=
\xi ((x \cdot \beta)(\eta))
-
\eta ((x \cdot \beta)(\xi))
-
(x \cdot \beta)[\xi,\eta]
\\
&=
\xi (x (\beta(\eta))-\beta(x \cdot \eta))
-
\eta (x (\beta(\xi))-\beta(x \cdot \xi))
\\
&\quad
-
x(\beta[\xi,\eta])+ \beta(x\cdot [\xi,\eta])
\\
((d''x) \cdot \beta)(\xi,\eta)
&=
((d''x)(\xi)) \cdot \beta)(\eta)
-
((d''x)(\eta)) \cdot \beta)(\xi)
\\
&=
((\xi \cdot x) \cdot \beta)(\eta)
-
((\eta \cdot x) \cdot \beta)(\xi)
\\
&=
(\xi \cdot x) (\beta(\eta))
-
\beta((\xi \cdot x) \cdot \eta)
-
(\eta \cdot x) (\beta(\xi))
+
\beta((\eta \cdot x) \cdot \xi)
\\
(x \cdot (d''\beta))(\xi,\eta)
&=
x((d''\beta)(\xi,\eta))
-
(d''\beta)(x \cdot \xi,\eta)
-
(d''\beta)(\xi,x \cdot\eta)
\\
&=
x(\xi(\beta(\eta))-\eta(\beta(\xi))-\beta[\xi,\eta])
\\
&\quad
-
(x \cdot \xi)(\beta(\eta))
+
\eta(\beta(x\cdot \xi))
+
\beta[x \cdot \xi,\eta]
\\
&\quad
-
\xi(\beta(x\cdot \eta))
+
(x \cdot \eta)(\beta(\xi))
+
\beta[ \xi,x \cdot\eta]
\endalign
$$
A straightforward comparison of terms gives
$$
\align
(d''(x \cdot \beta)
&-
(d''x) \cdot \beta
-
x \cdot (d''\beta))(\xi,\eta)
\\
&=
\xi(x(\beta(\eta)))
-
x(\xi(\beta(\eta)))
-
(\xi \cdot x - x \cdot \xi)(\beta(\eta))
\\
&\quad
-
(\eta(x(\beta(\xi)))
-
x(\eta(\beta(\xi)))
-
(\eta \cdot x - x \cdot \eta)(\beta(\xi)))
\\
&\quad
+\beta
(
x\cdot[\xi,\eta]
-
([x\cdot\xi,\eta]
+
[\xi,x\cdot \eta]
-(\xi \cdot x)\cdot \eta
+(\eta \cdot x)\cdot \xi))
\endalign
$$
This shows that the statement of the Lemma is true
when $\alpha$ has degree 1.
Since for two homogeneous elements $\alpha_1, \alpha_2$
of $\Cal A''$
$$
x(\alpha_1 \alpha_2) =
(x(\alpha_1))\alpha_2 + \alpha_1(x\alpha_2)
$$
the statement of the Lemma is true for
homogeneous $\alpha$ of arbitrary degree. \qed
\enddemo

\demo{Proof of Theorem {\rm 3.2}}
Suppose first that
$(A,L',L'')$ is a twilled Lie-Rinehart algebra,
that is to say, the compatibility conditions
(1.7.1) -- (1.7.3) are satisfied.
Then, by Lemma 3.3,
the identity (3.3.1$'$)
holds for every $x, y \in L$ and,
by Lemma 3.4,
the identity (3.4.1$'$)
holds for every $x\in L$ and
every homogeneous $\alpha\in \Cal A''$.
Since as a graded $\Cal A'' =\roman{Alt}_A(L'',A)$-module,
in fact as a
differential graded $\Cal A''$-module,
$\Cal L'= \roman{Alt}_A(L'',L')$
is generated by $L'$,
this implies that
$\Cal L'$ is a differential graded $R$-Lie algebra.
In fact, a straightforward calculation 
involving
(3.3.1$'$)
and
(3.4.1$'$)
shows that,
given $\alpha,\beta \in \Cal A''$ and
$x,y \in L'$, in view of (2.8.2),
$$
d''[\alpha \otimes_A x, \beta \otimes_A y]'
=
[d''(\alpha \otimes_A x), \beta \otimes_A y]'
+(-1)^{|\alpha|}[\alpha \otimes_A x, d''(\beta \otimes_A y)]'.
\tag3.2.1
$$
Moreover,
the truth of
the identity (3.4.1$'$)
for every $x\in L$ and
every homogeneous $\alpha\in \Cal A''$
implies that
$(\Cal A'',\Cal L')$
is a differential graded Lie-Rinehart algebra.
Conversely,
suppose that $\Cal L'$
is a differential graded
$R$-Lie algebra.
Then (3.2.1)
is manifestly true
for every $\alpha,\beta \in \Cal A''$ and
$x,y \in L'$.
Thus the identity (3.3.1$'$) then
holds a fortiori for every $x, y \in L'$, 
and from Lemma 3.3 we deduce at once 
that the compatibility condition (1.7.1) holds.  Furthermore, a 
straightforward comparison of terms involving only the differential graded 
$\Cal A''$-module structure and the identity
(3.3.1$'$)  gives, 
for $x,y\in L'$ and homogeneous 
$\alpha \in \Cal A''$, 
$$ 
\align 
d''[1 \otimes_A x, \alpha \otimes_A y]' 
 & -[d''(1\otimes_A x), \alpha \otimes_A y]' 
 -[1 \otimes_A x, d''(\alpha \otimes_A y)]'
\\ 
&= \left(d''(x\cdot \alpha)- (d''x)\cdot \alpha)- x\cdot 
(d''\beta)\right) y .  
\endalign 
$$ 
Since 
$d''[1 \otimes_A x, \alpha \otimes_A y]' 
- [d''(1\otimes_A x), \alpha \otimes_A y]' -[1 \otimes_A x, d''(\alpha 
\otimes_A y)]'$ 
is actually zero by assumption, 
$L'$ being finitely generated and projective as an $A$-module,
we conclude 
$$ 
d''(x\cdot \alpha)= 
(d''x)\cdot \alpha +x\cdot (d''\alpha) 
$$ 
for every $x\in L'$ and every 
homogeneous $\alpha \in \Cal A''$.  
Thus the identity (3.4.1$'$) holds for 
every $x\in L$ and every homogeneous 
$\alpha \in \roman{Alt}_A(L'',A)$.  From 
Lemmata 3.3 and 3.4 we deduce at once that the compatibility properties
(1.7.1) -- (1.7.3) hold.  \qed \enddemo 

As for the identity
(3.2.1),
we must again be a bit circumspect: a comment of formally the
same kind as that spelled out in Remark 2.8.6 is to be made here.
\smallskip
We remind the reader that the property P has been introduced
before Theorem 1.15.

\proclaim{Corollary 3.5}
Under the circumstances of {\rm (3.2)},
if in addition $L''$ has property {\rm P},
each of the three equivalent statements
in {\rm (3.2)} is equivalent to the operator
$[\parti',\parti'']$ 
in {\rm (1.15.4)} being zero.
\endproclaim

\demo{Proof} This is an immediate consequence of
Theorems 1.15 and 3.2. \qed
\enddemo

\smallskip\noindent
{\smc Example 3.6.} 
We return to the circumstances of
Example 1.3
and consider the
twilled Lie-Rinehart algebra
$(A,L',L'')$ arising from
a complex structure
on a smooth manifold $M$.
Thus,
to adjust the notation, the Cauchy-Riemann operator $\overline \partial$
now being identified with the operator $\parti''$,
$L'$ and $L''$ are the spaces of sections
of the holomorphic and antiholomorphic tangent bundles
of $M$, respectively.
The resulting differential graded Lie algebra
$$
(\Cal L', [\cdot,\cdot]', \parti'')
=(\roman{Alt}_A(L'',L'),[\cdot,\cdot]', \parti'')
$$
spelled out in Theorem 3.2
is
what is called
the {\it Kodaira-Spencer\/} algebra, cf. e.~g.  p. 337 of \cite\golmilfo;
it controls the infinitesimal deformations
of the complex structure of $M$.
In particular, $\Cal L'$
is the space of $\overline \partial$-forms
with values in the holomorphic tangent bundle of $M$,
and the differential 
$\parti''$ in $\Cal L'$
is the ordinary Cauchy-Riemann operator,
more customarily written $\overline \partial$.
Further, the Lie bracket is given by the formula (2.8.2)
above;
this bracket is {\it not\/}
given just by the shuffle product of
$\overline \partial$-forms
and the Lie-bracket of sections of the holomorphic
tangent bundle!
Cf. what is said at various places in the literature.
In fact, such a bracket would not
even be well defined
since the Lie bracket of vector fields is not a tensor.
It is also worthwhile pointing that,
in view of Theorem 3.2,
{\it the compatibility properties defining part of the
structure of the Kodaira-Spencer algebra
are equivalent to the integrability condition of the
initially only almost complex structure\/}.
\smallskip
\noindent
{\smc Remark 3.7.}
By means of a suitable generalization
of the graded crossed product Lie-Rinehart algebra
extension to the case where, in the notation of (2.9),
the left $(A,L)$-module $M$ is no longer assumed
to be finitely generated and projective
as an $A$-module, 
in a follow-up paper \cite\crosspro\ we will prove that
the statement of Theorem 3.2 still holds {\it without\/}
the hypothesis that $L'$ be finitely generated
and projective as an $A$-module.
Thus the bijective correspondence between twilled Lie-Rinehart 
and differential graded crossed product Lie-Rinehart structures
is valid in general.

\medskip\noindent{\bf 4. Differential Gerstenhaber algebras}\smallskip\noindent
The ground ring $R$ being fixed,
recall that a {\it Gerstenhaber algebra\/}
is a graded commutative $R$-algebra
$\Cal A$
together with a graded Lie bracket
from $\Cal A \otimes_R \Cal A$ to $\Cal A$
of degree $-1$ (in the sense that,
if $\Cal A$ is regraded down by one,
$[\cdot,\cdot]$ is an ordinary graded Lie bracket)
such that, for each homogeneous element $a$ of $\Cal A$,
$[a,\cdot]$ is a derivation of $\Cal A$ of degree $|a|-1$
where $|a|$ refers to the degree of $a$;
see \cite\geschthr\ 
where these objects are called G-algebras,
or \cite{\kosmathr,\liazutwo,\xuone};
for a Gerstenhaber algebra
$\Cal A$,
the bracket
from $\Cal A \otimes_R \Cal A$ to $\Cal A$
will henceforth be referred to as its {\it Gerstenhaber\/}
bracket.
Below we will interchangeably
talk about Gerstenhaber algebras and G-algebras.
We recall from
Theorem 5 of  \cite\geschthr\ 
that (i) the assignment
to
a Gerstenhaber algebra $\Cal A$ of the pair
$(A_0, A_1)$
consisting of the homogeneous degree zero and degree one
components
$A_0$ and $A_1$, respectively,
yields a functor
from Gerstenhaber algebras to Lie-Rinehart algebras,
and that this functor has a left adjoint which assigns
the exterior algebra
$\Lambda_AL$ 
in the category of $A$-modules
to the Lie-Rinehart algebra $(A,L)$,
together with the obvious bracket operation
on
$\Lambda_AL$ induced by the Gerstenhaber bracket structure.
Here $L$ is viewed to be concentrated in degree one.
\smallskip
\noindent
{\smc Definition 4.1.} Let $\Cal A$
be a bigraded commutative $R$-algebra.
We will say that
a bigraded bracket
$[\cdot,\cdot]\colon\Cal A \otimes_R \Cal A \to \Cal A$
of bidegree $(0,-1)$
is a {\it bigraded\/}
{\it Gerstenhaber bracket\/}
provided 
$[\cdot,\cdot]$ is an ordinary bigraded Lie bracket
when the second degree of
$\Cal A$ is regraded down by one, the first one being kept,
such that, for each homogeneous element $a$ of $\Cal A$
of bidegree $(p,q)$,
$[a,\cdot]$ is a derivation of $\Cal A$ of bidegree $(p,q-1)$;
a bigraded
$R$-algebra
with a
bigraded Gerstenhaber
bracket will be referred to
as
a {\it bigraded Gerstenhaber algebra\/}.

\medskip\noindent{\bf 4.2. The bigraded crossed product Gerstenhaber algebra}
\smallskip\noindent
Let $(A,L)$ be an (ungraded) Lie-Rinehart algebra, and let
$\Cal A$ be a graded commutative 
$A$-algebra together with an $L$-action
$L \otimes_R\Cal A \to \Cal A$ by derivations
such that the canonical map from $A$ to $\Cal A$ is
a morphism of left $(A,L)$-modules.
Let $\Cal L$ be the corresponding
crossed product
of $\Cal A$ and $L$ given in (2.8);
it is a graded $(R,\Cal A)$-Lie algebra.
Consider the bigraded $A$-algebra
$$
\Lambda_{\Cal A}\Cal L
=
\Cal A \otimes _A \Lambda_AL.
$$
It is the graded exterior $\Cal A$-algebra
on $\Cal L$
in an obvious sense whence the notation.
The graded Lie bracket
$[\cdot,\cdot]$
and the $L$-action on $\Cal A$
induce
a {\it bigraded\/} Gerstenhaber bracket
$$
[\cdot,\cdot]
\colon
\Lambda_{\Cal A}\Cal L
\otimes_R
\Lambda_{\Cal A}\Cal L
@>>>
\Lambda_{\Cal A}\Cal L
\tag4.2.1
$$
by means of the formulas
$$
\aligned
[\alpha \beta,\gamma] 
&=
\alpha [\beta,\gamma] +{(-1)}^{|\alpha||\beta|} \beta[\alpha,\gamma],
\quad \alpha,\beta,\gamma \in \Lambda_{\Cal A}\Cal L,
\\
[x,a] &= x(a), \quad x \in L,\, a \in\Cal A,
\\ 
[\alpha,\beta] 
&= -{(-1)}^{(|\alpha|-1)(|\beta|-1)}[\beta,\alpha],\quad
\alpha,\beta \in \Lambda_{\Cal A}\Cal L,
\endaligned
\tag4.2.2
$$
where $|\cdot|$ refers to the total degree, i.e. the sum of the
two bidegree components.
We refer to
the bracket (4.2.1)
as the {\it (bigraded)
crossed product bracket extension\/}
and to
$\Lambda_{\Cal A}\Cal L$
as the {\it (bigraded) crossed product Gerstenhaber algebra\/}
of $\Cal A$ with
the ordinary Gerstenhaber algebra
$\Lambda_A L$.
In terms of the latter
and the graded left $(A,L)$-module structure
on $\Cal A$,
the bigraded Gerstenhaber bracket
(4.2.1)
may be described in the following way
which, among others, gives an explicit formula:
Let $a,b \in \Cal A$ and 
$u = \alpha_1 \wedge \ldots \wedge\alpha_\ell \in \Lambda^{\ell}_AL$ and
$v = \alpha_{\ell+ 1}\wedge \ldots \wedge \alpha_n \in \Lambda^{n-\ell}_AL$,
where $\alpha_1, \ldots, \alpha_n \in L$;
then the ordinary Gerstenhaber bracket $[u,v]$ in $\Lambda_AL$
is given by the expression
$$
[u,v]
=
(-1)^{\ell} 
\sum_{j\leq \ell <k} (-1)^{(j+k)}
\lbrack \alpha_j,\alpha_k \rbrack \wedge 
\alpha_1\wedge \ldots \widehat{\alpha_j} \ldots \widehat{\alpha_k}
\ldots \wedge \alpha_n,
\tag4.2.3
$$
where $\ell = |u|$ is the degree of $u$,
cf. \cite\bv\ (1.1).
Writing
$au=a \otimes_Au$ and
$bv=b \otimes_Av$,
from
(4.2.2),
we obtain
$$
\aligned
[u,b] &=
[\alpha_1 \wedge \ldots \wedge\alpha_\ell, b]
\\
&=
\sum_{j=1}^{\ell} (-1)^{\ell -j}
\alpha_1 \wedge \ldots \widehat{\alpha_j}\ldots\wedge\alpha_\ell
[\alpha_j, b]
\\
&=
\sum _{j=1}^{\ell}(-1)^{\ell -j + (\ell-1)|b|}
\alpha_j (b)
\alpha_1 \wedge \ldots \widehat{\alpha_j}\ldots\wedge\alpha_\ell
\\
[a,v]&=[a,\alpha_{\ell+ 1}\wedge \ldots \wedge \alpha_n]
\\
&=
\sum_{j=1}^{n-\ell} {(-1)}^{j-1}
[a,\alpha_{\ell+ j}]\alpha_{\ell+ 1}\wedge \ldots 
\widehat{\alpha_{\ell+ j}}\ldots
\wedge \alpha_n
\\
&=
\sum_{j=1}^{n-\ell} {(-1)}^{j} 
\alpha_{\ell+ j}(a)\alpha_{\ell+ 1}\wedge \ldots 
\widehat{\alpha_{\ell+ j}}\ldots \wedge \alpha_n
\\
&=
\sum_{k=\ell+1}^n {(-1)}^{k-\ell}
\alpha_k(a)\alpha_{\ell+1}\wedge\ldots\widehat{\alpha_k}\ldots\wedge\alpha_n
\endaligned
$$
whence
$$
\aligned
[au, bv] &= 
{(-1)}^{(\ell-1)|b|} ab[u,v]
+a [u,b] v 
+{(-1)}^{(\ell-1+|a|)|b| +(\ell-1)(n-\ell-1)}b [a,v] u
\\
&=
{(-1)}^{(\ell-1)|b|} ab[u,v]
\\
&\quad
+
(-1)^{(\ell-1)|b|}
\sum _{j=1}^{\ell} (-1)^{\ell -j} a \alpha_j (b)
\alpha_1 \wedge \ldots \widehat{\alpha_j}\ldots\wedge\alpha_\ell \wedge v 
\\
&\quad
+{(-1)}^{(\ell-1)(|b| +n-\ell-1)} \sum _{k=\ell+1}^n (-1)^{k-\ell} 
\alpha_k(a)b\alpha_{\ell+1}\wedge\ldots\widehat{\alpha_k}\ldots\wedge\alpha_n
\wedge u .
\endaligned
$$
Here the sign in the expression
${(-1)}^{(\ell-1+|a|)|b| +(\ell-1)(n-\ell-1)}b [a,v] u$
arises from
first interchanging  $b$ with $au$
and thereafter interchanging
$u$ and $v$,
which necessitate the
signs 
${(-1)}^{(\ell-1+|a|)|b|}$
and
${(-1)}^{(\ell-1)(n-\ell-1)}$, respectively.
Consequently
$$
\aligned
[au, bv] &=
{(-1)}^{(\ell-1)|b|} ab[u,v]
\\
&\quad
+(-1)^{(\ell-1)|b|}\sum _{j\leq \ell} (-1)^{\ell -j} a \alpha_j (b)
\alpha_1 \wedge \ldots \widehat{\alpha_j}\ldots \wedge \alpha_n  
\\
&\quad
+ (-1)^{(\ell-1)|b|}\sum _{j>\ell} 
{(-1)}^{j-\ell}
\alpha_j(a) b\alpha_1 \wedge \ldots\widehat{\alpha_j}\ldots\wedge\alpha_n .
\endaligned
\tag4.2.4
$$

\smallskip
\noindent
{\smc Remark 4.2.5.}
We note that
$[au, bv]$ is the sum of
$\pm ab[u,v]$---which  involves only the product
$ab$ in $\Cal A$ and the Gerstenhaber bracket $[u,v]$ in
$\Lambda_AL$---and {\it two other terms\/},
which involve the action of $L$ on $\Cal A$
and the product in $\Cal A$.
Thus the 
(crossed product)
Gerstenhaber bracket
on 
$\Lambda_{\Cal A}\Cal L$
{\it cannot\/}
be written
just in terms of the product
on $\Cal A$ and the Gerstenhaber bracket on $\Lambda_AL$.
\smallskip
In particular,
under the circumstances of (2.9),
with
$\Cal A = \roman{Alt}_A(M,A)$,
the crossed product Gerstenhaber structure is available
for the crossed product
$(R,\Cal A)$-Lie algebra
$\Cal L= \Cal A \otimes _A L$
and yields a bigraded Gerstenhaber bracket
on the bigraded $A$-algebra
$$
\Lambda_{\Cal A}\Cal L
=
\roman{Alt}_A(M,A) \otimes _A \Lambda_AL.
\tag4.2.6
$$
When $L$ is finitely generated
and projective
as an $A$-module,
the bigraded
$A$-algebra
$\Lambda_{\Cal A}\Cal L$
may in fact
be written in the form
$$
\roman{Alt}_A(M,\Lambda_AL).
\tag4.2.7
$$
\smallskip
Recall that a {\it differential\/} Gerstenhaber algebra
$(\Cal A,[\cdot,\cdot],d)$ 
consists of a Gerstenhaber algebra
$(\Cal A,[\cdot,\cdot])$
together with a 
differential $d$
(of degree $+1$)
which endows 
$\Cal A$ with a differential graded
$R$-algebra structure
\cite\kosmathr, \cite\xuone;
in \cite\getzltwo, these objects are studied  under the name
{\it braid algebras\/}.
We will say that
a differential Gerstenhaber algebra
$(\Cal A,
[\cdot,\cdot], d)$ 
is {\it strict\/}
provided
$d$ behaves as a derivation for
the Gerstenhaber bracket $[\cdot,\cdot]$, that is,
$$
d[x,y] = [dx,y] -(-1)^{|x|} [x,dy], \quad x,y \in \Cal A.
$$
\smallskip
\noindent
{\smc Definition 4.3.} 
A bigraded Gerstenhaber algebra
$(\Cal A,[\cdot,\cdot])$
together with a 
differential $d$
of bidegree $(1,0)$
which endows $\Cal A$ with
a differential graded
$R$-algebra structure
will be said to be
a {\it differential bigraded\/} Gerstenhaber algebra
(or differential bigraded G-algebra),
written $(\Cal A,[\cdot,\cdot], d)$, 
provided
$d$ behaves as a derivation for
the bigraded Gerstenhaber bracket $[\cdot,\cdot]$, that is,
$$
d[x,y] = [dx,y] -(-1)^{|x|} [x,dy], \quad x,y \in \Cal A,
$$
where the total degree $|x|$ is the sum of the bidegrees.
\smallskip
We now return to the circumstances of (3.2):
Thus $(A,L',L'')$ 
is an almost twilled Lie-Rinehart algebra,
$L'$ being finitely generated
and projective as an $A$-module,
and
$$
(\Cal A'',\Cal L')=
(\roman{Alt}_A(L'',A),\roman{Alt}_A(L'',L'))
$$
is the corresponding graded crossed product Lie-Rinehart algebra,
cf. (3.1).
Furthermore,
$\Cal A'' = \roman{Alt}_A(L'',A)$
is endowed
with the ordinary Lie-Rinehart differential
$\parti''$
(cf. 1.11.4.1)
and 
$\Cal L' =\roman{Alt}_A(L'',L')$
with the Lie-Rinehart differential
$\parti''$
which comes from the given
left $(A,L'')$-module structure 
(1.4.4)
on
$L'$
and,
moreover, with the graded Lie bracket $[\cdot,\cdot]'$
given in {\rm (2.8.2)}
which comes from the given
left $(A,L')$-module structure 
(1.4.3) on
$L''$.
Consider the
resulting bigraded crossed product
Gerstenhaber algebra
(4.2.6),
with $M=L''$ and $L=L'$.
As a graded $A$-module,
this
Gerstenhaber algebra may be written
$$
\Lambda_{\Cal A''} \Cal L' 
=\roman{Alt}_A(L'',\Lambda_A L'),
$$
and the operator $\parti''$
induced 
by the given
left $(A,L'')$-module structure on
$L'$
is of bidegree $(1,0)$
and turns
$\roman{Alt}_A(L'',\Lambda_A L')$
into a differential graded $R$-algebra. 
We denote the bigraded Gerstenhaber bracket
on the latter by $[\cdot,\cdot]'$.

\proclaim{Theorem 4.4}
For an almost twilled Lie-Rinehart algebra
$(A,L',L'')$ having $L'$ finitely generated and projective
as an $A$-module,
$(\roman{Alt}_A(L'',\Lambda_A L'),[\cdot,\cdot]',\parti'')$
is a differential bigraded
Gerstenhaber algebra
if and only if
$(A,L',L'')$ is a twilled Lie-Rinehart algebra.
\endproclaim

It is clear that,
by symmetry,
when
$L''$ is finitely generated (and projective) as an $A$-module,
with the roles of $L'$ and $L''$ interchanged,
exactly the same statement as that spelled out in Theorem 4.4 is true.

\demo{Proof}
This is an immediate consequence of (3.2) above. \qed
\enddemo

\smallskip\noindent
{\smc Example 4.5.} 
We return to the circumstances of
Example 3.6 
and consider the
twilled Lie-Rinehart algebra
$(A,L',L'')$ arising from
a complex structure
on a smooth manifold $M$.
The corresponding
differential bigraded Gerstenhaber algebra
spelled out
in (4.4)
has as underlying bigraded $A$-module,
where $A$ is the algebra of smooth complex functions on $M$,
that of smooth complex valued 
$\overline \partial$-forms
with values in the exterior powers of the holomorphic
tangent bundle;
such  differential Gerstenhaber algebras
were studied in \cite\barakont\ and elsewhere.
For reasons explained in
(4.2.5) above,
the Gerstenhaber bracket
does {\it not\/} just involve
the shuffle product of 
$\overline \partial$-forms
and the Schouten-Nijenhuis
bracket of sections of the
exterior powers of the holomorphic
tangent bundle, though,
and 
in the corresponding Gerstenhaber bracket (4.2.1)
two others terms come into play, cf. the description
(4.2.4).
By symmetry,
interchanging the roles
of the holomorphic and antiholomorphic
tangent bundles,
we obtain as well a differential
(bigraded) Gerstenhaber algebra
which consists
of
smooth complex valued
$\partial$-forms
with values in the exterior powers of the antiholomorphic
tangent bundle.
It is also worthwhile pointing that,
in view of Theorem 4.4,
{\it the compatibility properties defining part of the
structure of the differential bigraded
Gerstenhaber algebra
are equivalent to the integrability condition of the
initially only almost complex structure\/}.

\smallskip
\noindent
{\smc Remark 4.6.}
In a follow-up paper \cite\crosspro\ 
we will generalize
the bigraded crossed product Gerstenhaber algebra
extension to the case where, in the notation of (2.9),
the left $(A,L)$-module $M$ is no longer assumed
to be finitely generated and projective
as an $A$-module, 
so that,
for arbitrary $M$,
a bigraded crossed product Gerstenhaber algebra
structure on a bigraded $A$-algebra of the kind
(4.2.7) 
results.
By means of this generalization,
we will then prove that
the statement of Theorem 4.4 still holds {\it without\/}
the hypothesis that $L'$ be finitely generated
and projective as an $A$-module.

\smallskip
We conclude this section with a description
of twilled Lie-Rinehart algebras
in terms of
Lie-Rinehart
bialgebras:
Let $L$ and $D$ be $(R,A)$-Lie algebras
which,
as $A$-modules, are  finitely generated and projective,
in such a way that,
as an $A$-module, $D$ is isomorphic to $L^* = \roman{Hom}_A(L,A)$.
We say
that $L$ and $D$ {\it are in duality\/}.
We write $d$ for the differential on
$\roman{Alt}_A(L,A)\cong \Lambda_AD$
coming from the Lie-Rinehart structure on $L$
and
$d_*$ for the differential on
$\roman{Alt}_A(D,A)\cong \Lambda_AL$
coming from the Lie-Rinehart structure on $D$.
Likewise we denote 
the Gerstenhaber bracket
on
$\Lambda_A L$
coming from 
the Lie-Rinehart structure on 
$L$ by
$[\cdot,\cdot]$
and
that
on
$\Lambda_A D$
coming from 
the Lie-Rinehart structure on 
$D$ by
$[\cdot,\cdot]_*$.

\proclaim{Proposition 4.7}
Given $L$ and $D$ in duality, the following are equivalent.
\newline
\noindent
{\rm (4.7.1)} The differential $d$ on
$\roman{Alt}_A(L,A) \cong \Lambda_A D$
and the Gerstenhaber bracket
$[\cdot,\cdot]$
on
$\Lambda_A D$
are related by
$$
d[x,y]_* = [dx,y]_* + [x,dy]_*, \quad x,y \in D.
$$
\newline
\noindent
{\rm (4.7.2)} The differential $d_*$ on
$\roman{Alt}_A(D,A) \cong \Lambda_A L$
and the Gerstenhaber bracket
$[\cdot,\cdot]$
on
$\Lambda_A L$
are related by
$$
d_*[x,y] = [d_*x,y] + [x,d_*y], \quad x,y \in L.
$$
\newline
\noindent
{\rm (4.7.3)} The differential $d$ on
$\roman{Alt}_A(L,A) \cong \Lambda_A D$
behaves as a derivation for
the Gerstenhaber bracket
$[\cdot,\cdot]_*$
in all degrees, that is to say
$$
d[x,y]_* = [dx,y]_* -(-1)^{|x|} [x,dy]_*, \quad x,y \in \Lambda_AD.
$$
\newline
\noindent
{\rm (4.7.4)} The differential $d_*$ on
$\roman{Alt}_A(D,A) \cong \Lambda_A L$
behaves as a derivation for
the Gerstenhaber bracket
$[\cdot,\cdot]$
in all degrees, that is to say
$$
d_*[x,y] = [d_*x,y] -(-1)^{|x|} [x,d_*y], \quad x,y \in \Lambda_AL.
$$
\endproclaim

\demo{Proof}
Since $d$ is a derivation for the algebra structure in
$\roman{Alt}_A(L,A)$, 
(4.7.1) and (4.7.3) are manifestly equivalent;
likewise,
since $d_*$ is a derivation for the algebra structure in
$\roman{Alt}_A(D,A)$,
(4.7.2) and (4.7.4) are equivalent.
The equivalence of
(4.7.3) and (4.7.4)
(the self-duality)
is established in
\cite\mackxu, \cite\kosmathr\ (3.3),
and \cite\liwextwo,
for the special case where $L$ and $D$ come from
real Lie algebroids. The argument given in these sources
is formal and carries
over. \qed
\enddemo

We will say that $(A,L,D)$
constitutes a
 {\it Lie-Rinehart\/}
bialgebra if
one (and hence any) of the (equivalent) conditions 
(4.7.1)--(4.7.4) is satisfied.
Thus, for a Lie-Rinehart bialgebra $(A,L,D)$,
$$
(\Lambda_A L, [\cdot,\cdot], d_*)
=
(\roman{Alt}_A(D,A), [\cdot,\cdot], d_*)
$$
is a strict {\it differential Gerstenhaber algebra\/},
and the same is true of
$$
(\Lambda_A D, [\cdot,\cdot]_*, d)
=
(\roman{Alt}_A(L,A), [\cdot,\cdot]_*, d);
$$
see \cite\kosmathr\ (3.5) for details.
In fact, a straightforward extension of an observation of
Y. Kosmann-Schwarzbach \cite\kosmathr\ 
shows that
Lie-Rinehart bialgebra structures
on $(A,L,D)$
and
strict differential Gerstenhaber algebra structures
on
$(\Lambda_A L, [\cdot,\cdot], d_*)$
or, what amounts to the same,
on
$(\Lambda_A D, [\cdot,\cdot]_*, d)$,
are equivalent notions.
This parallels the well known fact that Lie-Rinehart structures on $(A,L)$
are in bijective correspondence with differential graded $R$-algebra
structures on $\roman{Alt}_A(L,A)$.
\smallskip
Let $(A,L',L'')$ be an almost twilled Lie-Rinehart algebra,
having $L'$ and $L''$ finitely generated and projective
as $A$-modules.
The $(A,L')$-module structure (1.4.3) on $L''$
induces an
$(A,L')$-module on the dual ${L''}^*$ which, in turn, 
${L''}^*$ being viewed as an abelian Lie algebra and hence
abelian $(R,A)$-Lie algebra, gives rise to
the semi direct product $(R,A)$-Lie algebra
$L'\ltime {L''}^*$.
Likewise the $(A,L'')$-module structure (1.4.4) on $L'$
determines the corresponding
semi direct product $(R,A)$-Lie algebra
$L''\ltime {L'}^*$.
Plainly $L=L'\ltime {L''}^*$ and
$D=L''\ltime {L'}^*$ are in duality.
Consider the obvious adjointness isomorphisms
$$
\roman{Alt}_A(L'',\Lambda_A L')
@>>>
\roman{Alt}_A(L''\ltime {L'}^*, A )
=
\roman{Alt}_A(D, A )
\tag4.8.1
$$
and
$$
\Lambda_A L =
\Lambda_A (L'\ltime {L''}^*) 
@>>>
\roman{Alt}_A(L'',\Lambda_A L')
\tag4.8.2
$$
of bigraded $A$-algebras;
these isomorphisms are independent of the 
Lie-Rinehart
semi direct product
constructions and 
instead of $L'\ltime {L''}^*$ and
$L''\ltime {L'}^*$,
we could as well have written
$L'\oplus {L''}^*$
and
$L''\oplus {L'}^*$, respectively.
However, incorporating these semi direct product structures,
we see that, under (4.8.1), 
the Lie-Rinehart differential $d''$
on 
$\roman{Alt}_A(L'',\Lambda_A L')$
passes to the Lie-Rinehart differential
$d_*$ on
$\roman{Alt}_A(D,A)$
and that
under
(4.8.2) 
the (bigraded) Gerstenhaber bracket
$[\cdot,\cdot]$
on
$\Lambda_A L$
passes to the
bigraded Gerstenhaber bracket
$[\cdot,\cdot]'$
$\roman{Alt}_A(L'',\Lambda_A L')$.
Moreover, by construction, the differentials
on both sides of
(4.8.1) are derivations with respect to the multiplicative
structures.

\proclaim{Theorem 4.8}
For an almost twilled Lie-Rinehart algebra
$(A,L',L'')$ having $L'$ and $L''$ finitely generated and projective
as $A$-modules,
$(\roman{Alt}_A(L'',\Lambda_A L'),[\cdot,\cdot]',\parti'')$
is a differential bigraded
Gerstenhaber algebra
if and only if
$(A,L,D)$
is a Lie-Rinehart bialgebra.
\endproclaim

\demo{Proof}
The  property (4.7.2)
characterizing  $(A,L,D)$ to be a Lie-Rinehart bialgebra 
is plainly equivalent to 
$(\roman{Alt}_A(L'',\Lambda_A L'),[\cdot,\cdot]',\parti'')$
being a  
differential bigraded
Gerstenhaber algebra, cf. (4.3). \qed
\enddemo

The following is now immediate.

\proclaim{Corollary 4.9}
An almost twilled Lie-Rinehart algebra
$(A,L',L'')$ having $L'$ and $L''$ finitely generated and projective
as $A$-modules
is a true
twilled Lie-Rinehart algebra
if and only if
$(A,L,D) = (A,L'\ltime {L''}^*,L''\ltime {L'}^*)$
is a Lie-Rinehart bialgebra. \qed
\endproclaim

This result may be proved directly, i.~e. without 
the intermediate
differential bigraded
Gerstenhaber algebra
in (4.8). The reasoning is formally the same, though.
For the special case
where
$L'$ and $L''$ arise from Lie algebroids,
the statement of Corollary 4.9 may be deduced from what is said in
\cite\mackfift.

\medskip\noindent{\bf 5. Differential Batalin-Vilkovisky algebras}
\smallskip\noindent
In Section 1 of \cite\bv\  we
obtained an interpretation of the notion of a generator 
of a Gerstenhaber algebra 
$\Lambda_A L$  arising from a Lie-Rinehart algebra
$(A,L)$
which, for the special case where
$A$ is the ring of smooth functions and $L$ the Lie algebra
of smooth vector fields on a smooth manifold,
comes down to a result of Koszul \cite\koszulon.
In this section,
we will first
generalize this interpretation
to bigraded Gerstenhaber algebras.
We will then show that,
in the holomorphic context,
this extension is crucial
for an understanding of the Tian-Todorov Lemma,
of the Calabi-Yau condition,
and of the Batalin-Vilkovisky algebras arising from the mirror
conjecture.
\smallskip
For a bigraded Gerstenhaber algebra $\Cal A$
over $R$,
with bracket operation written $[\cdot,\cdot]$,
an $R$-linear operator $\Delta$
on $\Cal A$  
of bidegree $(0,-1)$
will be said to {\it generate\/}
the Gerstenhaber bracket
provided, for every homogeneous $a, b \in \Cal A$,
$$
[a,b] = (-1)^{|a|}\left(
\Delta(ab) -(\Delta a) b - (-1)^{|a|} a (\Delta b)\right);
\tag5.1
$$
the operator $\Delta$ is then called a {\it generator\/}.
A generator $\Delta$  is said to be {\it exact\/}
provided $\Delta\Delta$ is zero, that is, $\Delta$ is a differential;
an exact generator will henceforth be written $\partial$.
A bigraded Gerstenhaber algebra $\Cal A$
together with a generator $\Delta$
will be 
called a {\it weak bigraded Batalin-Vilkovisky\/} algebra
(or weak bigraded BV-algebra);
when
the generator is exact,
written $\partial$,
we will  refer to 
$(\Cal A,\partial)$
(more simply)
as a
{\it bigraded Batalin-Vilkovisky\/} algebra
(or bigraded BV-algebra).
\smallskip
It is clear that a generator
determines the bigraded Gerstenhaber bracket.
An observation due to Koszul \cite\koszulon\ 
(p. 261)
carries over to the bigraded case:
for any  bigraded Batalin-Vilkovisky algebra
$(\Cal A, [\cdot,\cdot], \partial)$,
the operator $\partial$ 
(which is exact by assumption)
behaves as a derivation for
the bigraded Gerstenhaber bracket $[\cdot,\cdot]$,
that is,
$$
\partial [x,y] = [\partial x,y] -(-1)^{|x|} [x,\partial y], 
\quad x,y \in \Cal A.
\tag5.2
$$
An exact generator $\partial$
does in general {\it not\/} behave
as a derivation for the multiplication
of $\Cal A$, though.

\medskip\noindent{\bf 5.3. The crossed product (weak) Batalin-Vilkovisky 
algebra}\smallskip\noindent
Let $(A,L)$ be a Lie-Rinehart algebra, and let
$\Cal A$ be a graded commutative 
$A$-algebra together with an $L$-action
$L \otimes_R\Cal A \to \Cal A$ by derivations
such that the canonical map from $A$ to $\Cal A$ is
a morphism of left $(A,L)$-modules.
Consider the 
crossed product $(R,\Cal A)$-Lie algebra
$\Cal L = \Cal A\otimes_A L$
(given in (2.8))
and the corresponding
bigraded crossed product Gerstenhaber algebra
$\Lambda_{\Cal A} \Cal L = \Cal A\otimes_A\Lambda_A L$
introduced in (4.2),
with
bigraded Gerstenhaber bracket
{\rm (4.2.1)}.

\proclaim{Proposition 5.3.1}
A generator $D$ for the  Gerstenhaber bracket of the
ordinary Gerstenhaber algebra
$\Lambda_AL$ 
admits a unique extension to a generator
$D_{\Cal A}$
of the bigraded crossed product Gerstenhaber algebra
bracket
{\rm (4.2.1)}
on
$\Lambda_{\Cal A} \Cal L = \Cal A\otimes_A\Lambda_A L$.
This extension
may be described
by means of the formula
$$
D_{\Cal A}(a \alpha) = a D(\alpha)
+
\sum_{i=1}^n (-1)^{i}(\alpha_i(a))
\alpha_1 \wedge \ldots \widehat{\alpha_i} \ldots \wedge \alpha_n 
\tag5.3.2
$$
where
$a \in \Cal A$
and 
$\alpha = \alpha_1 \wedge \ldots \wedge \alpha_n \in \Lambda_AL$.
Further, every
generator of the bigraded crossed product Gerstenhaber algebra
bracket on
$\Lambda_{\Cal A}\Cal L$ arises in this way.
\endproclaim

\demo{Proof} This is left to the reader. \qed
\enddemo
Additional insight
into
the generator $D_{\Cal A}$
will be offered in Section 7 below; see in particular
Theorem 7.6.
\smallskip
We refer to the generator $D_{\Cal A}$
of the 
bracket on
$\Lambda_{\Cal A}\Cal L$
given by (5.3.2) as the
{\it bigraded crossed product extension\/}
of the generator $D$
for the bracket on $\Lambda_AL$. 
The resulting
weak bigraded Batalin-Vilkovisky algebra
$(\Lambda_{\Cal A}\Cal L,D_{\Cal A})$
will be referred to
as the {\it bigraded crossed product\/}
of
$\Cal A$ and
$(\Lambda_L,D)$.
A
bigraded crossed product
of
$\Cal A$ and a true Batalin-Vilkovisky algebra $(\Lambda_L,\partial)$
is manifestly a
true bigraded Batalin-Vilkovisky algebra
$(\Lambda_{\Cal A}\Cal L,\partial_{\Cal A})$.
\smallskip
Suppose that,
as an $A$-module,
$L$ is finitely generated and projective
of finite constant rank, $n$.
The canonical pairing
$$
(\cdot,\cdot)
\colon
\Lambda_A^* L
\otimes_A
\Lambda_A^{n-*} L
@>>>
\Lambda_A^nL
\tag5.3.3
$$
of graded $A$-modules
is perfect and its adjoint
$$
\Lambda_A^* L
@>>>
\roman{Hom}_A(\Lambda_A^{n-*} L,\Lambda_A^nL)
=
\roman{Alt}_A^{n-*}(L,\Lambda_A^nL)
\tag5.3.4
$$
is an isomorphism
of graded $A$-modules. Given $x \in \Lambda_A^* L$,
write $\phi_x \in
\roman{Alt}_A^{n-*}(L,\Lambda_A^nL)$
for the image of $x$
under this isomorphism.
For an $(A,L)$-connection
\linebreak
$\nabla \colon M \to \roman{Hom}_A(L,M)$
on a left $A$-module $M$ we denote its
operator of covariant derivative
by
$$
d^{\nabla}
\colon
\roman{Alt}_A(L,M)
@>>>
\roman{Alt}_A(L,M).
$$

\proclaim{Proposition 5.3.5}
The relationship
$$
\phi_{\Delta(x)} =
d^{\nabla} \phi_x
\tag5.3.6
$$
establishes a bijective correspondence
between
generators $\Delta$ for the (ordinary) Gerstenhaber bracket
on
$\Lambda_AL$ and
$(A,L)$-connections $\nabla$ on
$\Lambda_A^nL$
in such a way that
exact generators $\Delta$
correspond to
left $(A,L)$-module structures
$\nabla$,
i.~e. flat $(A,L)$-connections, on
$\Lambda_A^nL)$.
\endproclaim

\demo{Proof}
See the Corollary in Section 2 of \cite\bv. \qed
\enddemo

\smallskip
Combining
Propositions 5.3.1 and 5.3.5
we arrive at the following.

\proclaim{Theorem 5.3.7}
The relationship {\rm (5.3.6)},
combined with
{\rm (5.3.2)},
establishes a bijective correspondence
between
generators $\Delta_{\Cal A}$ for the 
bigraded crossed product Gerstenhaber algebra
$\Lambda_{\Cal A}\Cal L$ 
and
$(A,L)$-connections $\nabla$ on
$\Lambda_A^nL$
in such a way that
exact generators $\Delta_{\Cal A}$
correspond to
left $(A,L)$-module structures
$\nabla$,
i.~e. flat $(A,L)$-connections, on
$\Lambda_A^n L$. \qed
\endproclaim

Another proof will be given at the end of Section 7 below;
see what is said after Corollary 7.10.

\proclaim {Corollary 5.3.8}
When $L$ is finitely generated and projective of constant rank
$n$ as an $A$-module,
generators 
for the bigraded Gerstenhaber bracket
on the bigraded crossed product
$\Lambda_{\Cal A}\Cal L$ 
always exist.
\endproclaim

\demo{Proof}
In fact, when $L$ is (finitely generated and) projective,
so is $\Lambda^n_AL$, whence
$(A,L)$-connections on $\Lambda^n_AL$ then always exist. \qed
\enddemo
\medskip\noindent{\bf 5.4. Incorporation of differentials} \smallskip\noindent
Let $(\Cal A,\Delta)$ be a weak bigraded Batalin-Vilkovisky algebra, 
write $[\cdot,\cdot]$
for the bigraded Gerstenhaber bracket generated by $\Delta$,
and let $d$ be a differential of bidegree $(+1,0)$
which endows $(\Cal A,[\cdot,\cdot])$ with a differential bigraded 
Gerstenhaber algebra structure. Consider the graded commutator
$$
[d,\Delta] = d\Delta + \Delta d
\tag5.4.1
$$ 
on $\Cal A$;
it is an 
operator of bidegree $(1,-1)$ and hence of total degree zero.
We will say that
$(\Cal A,\Delta,d)$
is a {\it weak differential\/}
bigraded Batalin-Vilkovisky algebra provided the commutator
$[d,\Delta]$
is zero.
In particular,
a weak differential
bigraded Batalin-Vilkovisky algebra 
$(\Cal A,\partial,d)$
which has $\partial$ exact
will be called a 
{\it differential bigraded Batalin-Vilkovisky algebra\/}.
On 
the underlying bigraded object $\Cal A$ of a differential
bigraded Batalin-Vilkovisky algebra 
$(\Cal A,\partial,d)$,
the graded commutator
$[d,\partial]$ 
manifestly behaves as a derivation for the 
bigraded
Gerstenhaber bracket
since
$d$ and $\partial$
both
behave as derivations for this bracket.

\smallskip
Various notions of differential Batalin-Vilkovisky algebras 
may be found in the literature,
cf. \cite\maninfiv\ (6.1.1)
(dGBV-algebras),
\cite\kosmathr\ 
(cf. e.~g. 
the differential exact Gerstenhaber
algebras 
on p. 154),
\cite\xuone\ (Question 4 in Section 5
where a concept of strong differential BV-algebra
occurs).
Our notion of differential bigraded Batalin-Vilkovisky algebra
does {\it not\/}
coincide with any of these.

\proclaim {Lemma 5.4.2}
Given
a weak bigraded
Batalin-Vilkovisky algebra
$(\Cal A,\Delta)$
and  an operator
$\delta$ of bidegree $(1,0)$
which behaves
as a derivation
of degree $1$
for the
bigraded $R$-algebra $\Cal A$, the following are equivalent:
\newline\noindent
{\rm (i)} The operator
$\delta$
behaves as a derivation for
the bigraded Gerstenhaber
bracket $[\cdot,\cdot]$
on $\Cal A$ generated by $\Delta$, that is to say,
$$
\delta[x,y] = [\delta x,y] -(-1)^{|x|} [x,\delta y], \quad x,y \in \Cal A.
$$
\newline\noindent
{\rm (ii)} The graded commutator
$[\delta,\Delta]$
behaves as a derivation
of degree $0$
for the
bigraded $R$-algebra $\Cal A$, that is,
$$
[\delta,\Delta](ab)
= 
([\delta,\Delta]a)b
+a([\delta,\Delta]b),
\quad
a,b \in \Cal A.
$$
\endproclaim

\demo{Proof}
Let $a,b$ be homogeneous elements of $\Cal A$. Then
$$
\align
\delta [a,b]&=
(-1)^{|a|}\left(\delta\Delta(ab) -\delta((\Delta a) b) \right)
- \delta(a (\Delta b))
\\
&=
(-1)^{|a|}\left(
-\Delta\delta(ab)
+[\delta,\Delta](ab)
-(\delta\Delta a)b+(-1)^{|a|}(\Delta a)(\delta b)\right)
\\
&\qquad
-(\delta a)(\Delta b) -(-1)^{|a|} a (\delta\Delta b)
\\
&=
(-1)^{|a|-1}\left(\Delta ((\delta a)b +(-1)^{|a|} a (\delta b))
-[\delta,\Delta](ab)
-(\Delta \delta  a)b 
+([\delta,\Delta]a)b 
\right)
\\
&\qquad
+
(\Delta a) (\delta b) 
- (\delta a)(\Delta b) +(-1)^{|a|} a (\Delta \delta  b)
-(-1)^{|a|} a([\delta,\Delta]b) 
\\
&=
(-1)^{|a|-1}(\Delta ((\delta a)b)
-\Delta (a \delta b) 
+(-1)^{|a|}
(\Delta \delta a)b
\\
&\qquad
+
(\Delta a) (\delta b) 
- (\delta a)(\Delta b) 
+(-1)^{|a|} a (\Delta \delta b)
\\
&\qquad
+(-1)^{|a|} \left(
 [\delta,\Delta](ab)-([\delta,\Delta]a)b-a([\delta,\Delta]b) 
\right)
\\
&=
(-1)^{|\delta a|}(\Delta ((\delta a)b)
-(\Delta \delta a)b)
- (\delta a)(\Delta b) 
\\
&\qquad 
-\Delta (a \delta b) 
+(\Delta a) (\delta b) 
+(-1)^{|a|} a (\Delta \delta b)
\\
&\qquad
+(-1)^{|a|} \left(
 [\delta,\Delta](ab)-([\delta,\Delta]a)b-a([\delta,\Delta]b) 
\right)
\\
&=
[\delta a,b] - (-1)^{|a|}[a,\delta b]
+(-1)^{|a|} \left(
 [\delta,\Delta](ab)-([\delta,\Delta]a)b-a([\delta,\Delta]b) 
\right)
\endalign
$$
This establishes the claim. \qed
\enddemo

\proclaim {Corollary 5.4.3}
For any 
weak differential bigraded Batalin-Vilkovisky algebra
$(\Cal A,\Delta,d)$,
the differential
$d$
behaves as a derivation for
the bigraded Gerstenhaber
bracket $[\cdot,\cdot]$
on $\Cal A$ generated by $\Delta$, that is to say,
$$
d[x,y] = [dx,y] -(-1)^{|x|} [x,dy], \quad x,y \in \Cal A.
$$
In other words,
$(\Cal A ,[\cdot,\cdot], d)$
is a differential bigraded Gerstenhaber algebra.
\endproclaim

Notice that,
under the circumstances of
(5.4.3),
$\Delta$ need not behave as a derivation for the
bigraded Gerstenhaber bracket
unless $\Delta$ is exact.
\smallskip

\proclaim{Theorem 5.4.4}
Suppose that
$L'$ is finitely generated and projective
as an $A$-module, 
let $\Delta'$ be a generator for the bigraded 
Gerstenhaber
bracket $[\cdot,\cdot]'$
of
the
bigraded
crossed product Gerstenhaber algebra
$(\roman{Alt}_A(L'',\Lambda_A L'),[\cdot,\cdot]')$,
and write
$\parti''$ for the Lie-Rinehart differential
{\rm (1.11.4.1)}
induced by the $(A,L'')$-action on  $\Lambda_A L'$.
Then 
$[\parti'',\Delta'] (=\parti ''\Delta'+\Delta'\parti'')$
is 
a derivation
(of bidegree $(1,-1)$)
for the bigraded $R$-algebra
$\roman{Alt}_A(L'',\Lambda_A L')$
if and only if
$(A,L',L'')$ is a twilled Lie-Rinehart algebra.
In particular,
when
$(\roman{Alt}_A(L'',\Lambda_A L'),\Delta',\parti'')$
is a weak differential bigraded
Batalin-Vilkovisky algebra
(i.e.~ when $[\parti'',\Delta']$ is zero),
$(A,L',L'')$ is necessarily a twilled Lie-Rinehart algebra.
\endproclaim

\demo{Proof}
This is a consequence of Lemma 5.4.2, combined with Theorem 4.4.\qed
\enddemo

It is clear that, when
$L''$
is finitely generated and projective as an $A$-module,
the same statements as those given in Theorem 5.4.4 can be made,
with the roles of $L'$ and $L''$ interchanged.
Exploiting the generalization of the bigraded crossed product
Gerstenhaber algebra mentioned already in (4.6),
in \cite\crosspro,
we will prove that the statement of Theorem 5.4.4 holds {\it without\/}
the hypothesis
that $L'$ be finitely generated and projective as an $A$-module;
this then
yields a result which is symmetric in $L'$ and $L''$.

\proclaim{Proposition 5.4.5}
Under the circumstances of {\rm (5.4.4)}, the adjoint
$$
\Lambda_A^* L'
@>>>
\roman{Hom}_A(\Lambda_A^{n-*} L',\Lambda_A^nL')
=
\roman{Alt}_A^{n-*}(L',\Lambda_A^nL')
$$
of the corresponding pairing
{\rm (5.3.3)}
is an isomorphism
of graded $(A,L'')$-modules,
whence
$(\roman{Alt}_A(L'',\Lambda_A L'),\Delta',\parti'')$
is a weak differential bigraded
Batalin-Vilkovisky algebra
(i.~e. $[\Delta',\parti'']$ is zero
but $\Delta'$ is not necessarily exact)
if and only if the $(A,L')$-connection
$$
\nabla'\colon
\Lambda^n_A L' 
@>>>
\roman{Hom}_A(L',\Lambda^n_A L')
$$
on $\Lambda^n_A L'$
corresponding to
$\Delta'$
(spelled out in {\rm (5.3.7)})
is $L''$-invariant,
with reference to the 
induced 
$L''$-actions
on
$\Lambda^n_A L'$
and
$\roman{Hom}_A(L',\Lambda^n_A L')$.
\endproclaim

\demo{Proof}
This is straightforward and left to the reader. \qed \enddemo

\proclaim{Theorem 5.4.6}
Under the circumstances of {\rm (5.4.4)} if, in addition, 
as an $A$-module, $L'$ 
(being finitely generated 
and projective) is of finite
constant rank $n$ (say),
the adjoint {\rm (5.3.4)\/}
induces an isomorphism
$$
(\roman{Alt}_A(L'',\Lambda^*_A L'),\parti'')
@>>>
(\roman{Alt}_A(L'',\roman{Alt}_A^{n-*}(L',\Lambda_A^nL')),\parti'')
\tag5.4.6.1
$$
of chain complexes 
which establishes a bijective correspondence between
generators
$\Delta'$
for the bigraded Gerstenhaber bracket
on the left-hand side 
and  
operators
$\partial_{\nabla'}$
of covariant derivative
on the right-hand side, 
for 
a uniquely determined $(A,L')$-connection 
$\nabla'$ on
$\Lambda_A^nL'$.
Under this correspondence,
generators $\Delta'$ turning
$(\roman{Alt}_A(L'',\Lambda_A L'),\Delta',\parti'')$
into a weak differential bigraded
Batalin-Vilkovisky algebra
correspond to $L''$-invariant connections
$\nabla'$
(on $\Lambda_A^nL'$)
and
generators $\Delta'$ turning
$(\roman{Alt}_A(L'',\Lambda_A L'),\Delta',\parti'')$
into a true differential bigraded
Batalin-Vilkovisky algebra
correspond to flat 
$L''$-invariant connections $\nabla'$.
\endproclaim

\demo{Proof}
This follows readily from Proposition 5.4.5. \qed
\enddemo

An interpretation of (5.4.6.1)
within the framework of
(co)homological duality for differential graded Lie-Rinehart
algebras will be given in Section 7 below.
\smallskip
A special case is worthwhile spelling out.
For this purpose
we observe that
any closed $L''$-invariant $A$-valued
1-form $\alpha \colon L' \to A$
on $L'$
determines 
an $L''$-invariant
$(A,L')$-module structure on $A$,
i.~e.
a flat 
$L''$-invariant $(A,L')$-connection $\nabla_\alpha$
on $A$,
whose
operator 
$\parti'_\alpha$ 
of covariant derivative
is determined by 
$$
\parti'_\alpha
\colon
A @>>>
\roman{Hom}_A(L', A),
\quad
(\parti'_\alpha (1))(x) = \alpha (x),
\quad x \in L'.
\tag5.4.7.1
$$
This operator plainly extends
to the corresponding Lie-Rinehart operator
on
\linebreak
$\roman{Alt}_A^*(L'',\roman{Alt}_A^{*}(L',A))$
(determined by the
$L''$-invariant
$(A,L')$-module structure on $A$),
and we continue to denote this operator by
$\parti'_\alpha$;
thus
$(\roman{Alt}_A^*(L'',\roman{Alt}_A^{*}(L',A)),\parti'_\alpha, d'')$
is a bicomplex.
When $\alpha$ is zero, this is just the
ordinary bicomplex of the kind (1.15.4).

\proclaim{Theorem 5.4.7}
Under the circumstances of {\rm (5.4.6)}, suppose
in addition that
there is
an $A$-valued $n$-form $\Lambda_A^nL' \to A$
on $L'$
yielding 
an isomorphism of $A$-modules which is
invariant under $L''$
(i.~e. which is an isomorphism of $(A,L'')$-modules).
Then a choice of such an $n$-form 
$\Omega\in \roman{Alt}_A^n(L',A)$
induces an isomorphism
$$
\Omega_\flat\colon
(\roman{Alt}_A^*(L'',\Lambda^*_AL'),\parti'')
@>>>
(\roman{Alt}_A^*(L'',\roman{Alt}_A^{n-*}(L',A)),\parti'')
\tag5.4.7.2
$$
of chain complexes over the ground ring $R$,
in fact, of differential graded
\linebreak
$(\roman{Alt}_A^*(L'',A),\parti'')$-modules.
Under this isomorphism,
the operator $\parti'$ on the 
right-hand side
$\roman{Alt}_A^*(L'',\roman{Alt}_A^{n-*}(L',A))$
of {\rm (5.4.7.2)}
corresponds to a uniquely determined exact
generator $\partial_\Omega$ 
for the  Gerstenhaber bracket on
the differential bigraded Gerstenhaber algebra
$(\roman{Alt}_A^*(L'',\Lambda^*_AL'),[\cdot,\cdot]',\parti'')$
on the left-hand side of 
{\rm (5.4.7.2)}.
Furthermore,
if $\alpha \colon L' \to A$
is any closed $L''$-invariant $A$-valued
{\rm 1-}form 
on $L'$,
the corresponding Lie-Rinehart operator 
$\parti'_\alpha$ 
on 
the 
right-hand side
$(\roman{Alt}_A^*(L'',\roman{Alt}_A^{n-*}(L',A)),\parti'')$
of {\rm (5.4.7.2)}
induces
as well a uniquely determined exact
generator $\partial_{(\Omega,\alpha)}$ 
(say)
for the  Gerstenhaber bracket 
$[\cdot,\cdot]'$
on the left-hand side
$(\roman{Alt}_A^*(L'',\Lambda^*_AL'),[\cdot,\cdot]',\parti'')$
of {\rm (5.4.7.2)},
and every exact generator
for this  Gerstenhaber bracket 
arises in this way.
\endproclaim

Thus the choice of $\Omega$ enables us to rewrite
the differential bigraded algebra
$(\roman{Alt}_A^*(L'',\roman{Alt}_A^{n-*}(L',A)),\parti',\parti'')$
as a differential bigraded Batalin-Vilkovisky algebra.
\smallskip

Theorem 5.4.7
is a special case of Theorem 5.4.6,
with $\Delta'$ corresponding to a flat
$(A,L')$-connection on $\Lambda^n L$
which is invariant under
the $L''$-action.
A direct proof of Theorem 5.4.7 will be given after (5.4.11).
\smallskip
Let $M$ be a smooth complex $n$-manifold, and write
$\tau_M$ and $\overline \tau_M$
for the holomorphic and antiholomorphic tangent bundles of 
$M$. For consistency with notation used in the literature,
we momentarily write  $\partial$ and $\overline \partial$
for the operators which correwspond, under our more general
circumstances, to our operators
$\parti'$ and $\parti''$, respectively.
Conflict with the notation $\partial$ for an exact generator of a 
Batalin-Vilkovkisky algebra will be avoided since
such a generator will
be written
$\partial_{\cdot}$
with an appropriate subscript.

\proclaim{Corollary 5.4.8} {\rm (Tian-Todorov Lemma)}
If $M$ admits
a holomorphic volume form, 
a choice $\Omega$
of holomorphic volume form
induces an isomorphism
$$
\Omega_\flat\colon
(\Gamma (\Lambda^* {\overline\tau^*_M} \otimes \Lambda^*\tau_M),
\overline \partial)
@>>>
(\Gamma (\Lambda^* {\overline\tau^*_M} \otimes \Lambda^{n-*}\tau^*_M),
\overline \partial)
\tag5.4.8.1
$$
of chain complexes,
in fact, of modules
over the differential graded algebra
$(\Gamma (\Lambda^* {\overline\tau^*_M}),\overline \partial)$
of $\overline \partial$-forms
defined only on the
antiholomorphic tangent bundle $\overline \tau_M$.
Under this isomorphism,
the operator $\partial$
on the right hand side 
of {\rm (5.4.8.1)}
corresponds, on the left-hand side,
to an exact  operator $\partial_\Omega$
which turns
$$
(\Gamma (\Lambda^* {\overline\tau^*_M} \otimes \Lambda^*\tau_M),
\partial_\Omega,\overline \partial)
$$
into a differential bigraded Batalin-Vilkovisky algebra.
\endproclaim

\demo{Proof} This is just a special case of the first statement
of (5.4.7). \qed
\enddemo

The resulting
isomorphism
$$
\Omega_\flat\colon
(\Gamma (\Lambda^* {\overline\tau^*_M} \otimes \Lambda^*\tau_M);
\overline \partial, \partial_{\Omega})
@>>>
(\Gamma (\Lambda^* {\overline\tau^*_M} \otimes \Lambda^{n-*}\tau^*_M);
\overline \partial, \partial)
\tag5.4.8.2
$$
identifies the 
differential bigraded Batalin-Vilkovisky algebra
on the left-hand side with the Dolbeault complex
(spelled out on the right-hand side), as pointed out
in the introduction.

\proclaim{Addendum 5.4.9} 
Under the circumstances
of {\rm (5.4.8)},
if $\alpha$ is any holomorphic $1$-form on $M$,
the operator 
$\parti'_\alpha$ 
of covariant derivative
(given by {\rm (5.4.7.1)})
on 
the right-hand side of
{\rm (5.4.8.1)},
for the corresponding  flat 
holomorphic
connection $\nabla_\alpha$
(say)
on $\Lambda^n\tau^*_M$,
corresponds as well to a uniquely determined exact
generator $\partial_{(\Omega,\alpha)}$ 
for the  Gerstenhaber bracket 
on 
the left-hand side of
{\rm (5.4.8.1)},
and every exact generator
for this  Gerstenhaber bracket 
arises in this way.
\endproclaim

\demo{Proof} This is indeed a special case of the 
\lq\lq Furthermore\rq\rq\ 
statement
of (5.4.7). \qed
\enddemo

The corresponding bigraded Gerstenhaber algebra
is of course just the corresponding
crossed product Gerstenhaber algebra.
The existence of a holomorphic volume form is
a strong kind of orientability condition; 
it is implied by  the {\it Calabi-Yau\/} condition
$c_1=0$.
The statement of Corollary 5.4.8 includes what is referred to in
the literature as the
{\it Tian-Todorov\/} lemma
\cite{\barakont,\,\golmilfo,\,\tianone,\,\todortwo}.
This lemma arises here as a natural consequence
of our theory of 
differential bigraded Batalin-Vilkovisky algebras
having as underlying bigraded Gerstenhaber
algebra a crossed product Gerstenhaber algebra.
Notice that the 
description (4.2.4) of
the Gerstenhaber bracket
(4.2.1)
generated by
$\partial_\Omega$ 
shows that
this bracket
does {\it not\/} just involve
the shuffle product of 
$\overline \partial$-forms
and the Schouten-Nijenhuis
bracket of sections of the
exterior powers of the holomorphic
tangent bundle, 
and two other terms 
(spelled out in greater generality in (4.2.4))
come into play,
cf. (4.5) above.

\smallskip\noindent
{\smc Remark.}
When $M$ is compact, by Serre duality,
the statement of
Corollary 5.4.8
holds as well
with the holomorphic and antiholomorphic tangent
bundles interchanged.
\smallskip

Whether or not 
$M$  has a holomorphic volume form, 
under the present circumstances,
{\rm (5.4.6.1)\/}
has the form
$$
(\Gamma (\Lambda^* {\overline\tau^*_M} \otimes \Lambda^*\tau_M),
[\cdot,\cdot]',\overline \partial)
@>>>
(\Gamma (\Lambda^* {\overline\tau^*_M} 
\otimes \Lambda^{n-*}\overline\tau_M
\otimes \Lambda^n \tau_M),\overline\partial)
\tag5.4.10
$$
and is in fact an isomorphism of chain complexes 
from
the $\overline\partial$-forms
with values in the exterior powers of the holomorphic
tangent bundle
onto 
the $\overline\partial$-forms
with values in the exterior powers of the holomorphic
cotangent bundle, tensored with
the highest exterior power
of the holomorphic tangent bundle.

\proclaim{Corollary 5.4.11} 
Suppose that
the highest exterior power
of the holomorphic tangent bundle
has merely a {\it holomorphic\/} connection
$\nabla'$
which is not necessarily flat.
Via {\rm (5.4.10)}, 
the corresponding operator
$d^{\nabla'}$ of covariant derivative
on the right-hand side
of {\rm (5.4.10)} 
induces on the left-hand side
$(\Gamma \Lambda^* {\overline\tau^*_M} \otimes \Lambda^*\tau^* M,
\overline\partial)$
(of {\rm (5.4.10)}), 
that is, on the  differential bigraded Gerstenhaber
algebra of $\overline\partial$-forms
with values in the exterior powers of the holomorphic
tangent bundle,
a generator $\Delta'$ so that
$(\Gamma (\Lambda^* {\overline\tau^*_M} \otimes \Lambda^*\tau_M),
\Delta',\overline \partial)$
is a weak differential bigraded Batalin-Vilkovisky
algebra.
\endproclaim

We note that the compatibility condition
$[\Delta',\overline\partial]=0$,
which defines part of the structure of
the weak differential bigraded Batalin-Vilkovisky
algebra occuring in the statement of Corollary 5.4.11,
corresponds precisely to the
holomorphicity of 
the connection $\nabla'$.

\demo{Proof} This follows at once from Theorem 5.4.6. \qed
\enddemo

\demo{Direct proof of Theorem {\rm 5.4.7}}
The $n$-form
$\Omega$
induces an isomorphism
$$
(\roman{Alt}_A^*(L'',\roman{Alt}_A^{n-*}(L',\Lambda_A^nL')),\parti'')
@>>>
(\roman{Alt}_A^*(L'',\roman{Alt}_A^{n-*}(L',A)),\parti'')
\tag5.4.7.3
$$
of chain complexes.
This relies on the fact that $\parti'' \Omega = 0$
(i.~e. the holomorphicity of $\Omega$ under the circumstances of (5.4.8)).
The composite of (5.4.7.3)
with
(5.4.6.1)
yields the asserted
isomorphism
(5.4.7.2)
of chain complexes from
$(\roman{Alt}_A^*(L'',\Lambda^*_AL'),\parti'')$
onto
$(\roman{Alt}_A^*(L'',\roman{Alt}_A^{n-*}(L',A)),\parti'')$.
Under this isomorphism,
to the operator $\parti'$
on
$(\roman{Alt}_A^*(L'',\roman{Alt}_A^{n-*}(L',A)),\parti'')$
corresponds
to a generator $\partial_\Omega$ 
for the bigraded Geratenhaber algebra structure
on 
$(\roman{Alt}_A^*(L'',\Lambda^*_AL'),\parti'')$.
Since 
$(\roman{Alt}_A^*(L'',\roman{Alt}_A^{n-*}(L',A)),\parti'',\parti')$
is a bicomplex,
$[\parti',\parti''] = 0$.
Consequently
$[\partial_\Omega,\parti''] = 0$
on
$\roman{Alt}_A(L'',\Lambda_AL')$,
that is,
$(\roman{Alt}_A(L'',\Lambda_AL'),\partial_\Omega, \parti'')$
is a differential bigraded
Batalin-Vilkovisky algebra. \qed
\enddemo

\smallskip\noindent
{\smc Addendum to the proof.} The $n$-form $\Omega$ 
plainly endows
$\Lambda_A^nL'$
with a left $(A,L')$-module structure.
In view of Theorem 5.3.7,
$\Omega$ 
thus induces a generator
for the 
bigraded Gerstenhaber bracket
$[\cdot,\cdot]'$
on $\roman{Alt}_A(L'',\Lambda_AL')$;
this generator is just $\partial_\Omega$.

\medskip\noindent {\bf 6. Twilled Lie-Rinehart algebras and 
differential homological algebra}
\smallskip\noindent
In this section we give an interpretation of twilled Lie-Rinehart
algebras in the framework of differential homological algebra.
In the next section,
we will use this interpretation to deduce
a 
differential homological algebra
interpretation of the generator
of a differential bigraded Batalin-Vilkovisky algebra
arising from a twilled Lie-Rinehart algebra.
\smallskip
Let 
$(\Cal A,{\Cal L})$
be a differential graded
Lie-Rinehart algebra.
The
{\it universal object\/}
${(U(\Cal A,{\Cal L}),\iota_{\Cal L},\iota_{\Cal A})}$
for $(\Cal A,{\Cal L})$
is a
 differential graded
 $R$-algebra $U(\Cal A,{\Cal L})$ together with a morphism
${
\iota_{\Cal A}
\colon 
\Cal A
\longrightarrow
U(\Cal A,{\Cal L})
}$
of 
 differential graded
$R$-algebras
and
a morphism
${
\iota_{\Cal L}
\colon 
{\Cal L}
\longrightarrow
U(\Cal A,{\Cal L})
}$
of  differential graded Lie algebras over $R$
having the properties
$$
\align
\iota_{\Cal A}(a)\iota_{\Cal L}(\alpha) &= \iota_{\Cal L}(a\,\alpha),
\\
\iota_{\Cal L}(\alpha)\iota_{\Cal A}(a) - 
(-1)^{|\alpha||a|}
\iota_{\Cal A}(a)\iota_{\Cal L}(\alpha) &= \iota_{\Cal A}(\alpha(a)),
\endalign
$$
and
${(U(\Cal A,{\Cal L}),\iota_{\Cal L},\iota_{A})}$
is {\it universal\/} among triples
${(B,\phi_{\Cal L},\phi_{\Cal A})}$
having these properties.
More precisely:
\proclaim{6.1.1}
Given 
{\rm (i)}
another
 differential graded
 $R$-algebra $B$, viewed at the same time as a 
 differential graded
Lie algebra over
$R$,
\newline
\noindent
{\rm (ii)}
a morphism
$
\phi_{\Cal L} 
\colon
{\Cal L}
\longrightarrow
B
$ 
of  differential graded Lie algebras over $R$, and
\newline
\noindent
{\rm (iii)} a morphism
$
\phi_{\Cal A} 
\colon
\Cal A
\longrightarrow
B
$
of 
 differential graded
$R$-algebras,
\smallskip
\noindent
so that, for ${\alpha \in {\Cal L},\  a \in \Cal A}$,
$$
\align
\phi_{\Cal A}(a)\phi_{\Cal L}(\alpha) &= \phi_{\Cal L}(a\,\alpha),
\tag6.1.2
\\
\phi_{\Cal L}(\alpha)\phi_{\Cal A}(a) - 
(-1)^{|\alpha||a|}
\phi_{\Cal A}(a)\phi_{\Cal L}(\alpha) &= \phi_{\Cal A}(\alpha(a)),
\tag6.1.3
\endalign
$$
there is a unique morphism
${
\Phi 
\colon
U(\Cal A,{\Cal L})
\longrightarrow
B
}$
of 
differential graded
$R$-algebras
so that
$$
\Phi\,\iota_{\Cal A} = \phi_B,\quad
\Phi\,\iota_{\Cal L} = \phi_{\Cal L}.
$$
\endproclaim

\smallskip
The universal differential graded algebra 
$U(\Cal A,\Cal L)$
may be obtained
in the
customary way
as the quotient of the differential graded tensor $R$-algebra
$\roman T(\Cal A \oplus \Cal L)$
of the direct sum
$\Cal A \oplus \Cal L$,
viewed merely as an $R$-module,
by the differential graded ideal
generated in
$\roman T(\Cal A \oplus \Cal L)$
by all elements
$$
\alpha\otimes_R \beta - (-1)^{|\alpha||\beta|}
\beta\otimes_R \alpha -[\alpha,\beta],
\quad
\alpha \otimes_R a - (-1)^{|\alpha||a|} a \otimes_R \alpha   - \alpha(a),
$$
for $a \in \Cal A$ and
$\alpha,\beta \in \Cal L$.
The morphisms
$\iota_{\Cal A}$
and
$\iota_{\Cal L}$
are then the obvious ones.
Thus, as a graded $R$-algebra,
$U(\Cal A,\Cal L)$
is generated by the $a \in \Cal A$ and
the $\alpha \in \Cal L$ subject to the relations
$$
\align
\alpha \beta - (-1)^{|\alpha||\beta|}
\beta \alpha &= [\alpha,\beta]
\tag6.1.4.1
\\
\alpha a - (-1)^{|\alpha||a|} a \alpha &= \alpha(a),
\tag6.1.4.2
\endalign
$$
for $a \in \Cal A$ and
$\alpha,\beta \in \Cal L$.
Furthermore,
since 
$U(\Cal A,\Cal L)$
is generated by
$\Cal A$ and
$\Cal L$,
the differential $d$ (say) 
on
$\Cal A$ and
$\Cal L$
extends to a
unique differential on
$U(\Cal A,\Cal L)$
provided it extends at all.
However, the differential is compatible with
(6.1.4.1) by assumption, and
a little thought reveals that it is compatible
with the relations
(6.1.4.2),
whence
the differential extends to a
unique differential $d$ on
$U(\Cal A,\Cal L)$.
In particular, as a graded $\Cal A$-module,
$U(\Cal A,\Cal L)$
is generated by monomials of the kind
$\alpha_1 \dots \alpha_m$ of arbitrary length $m$,
where $\alpha_j \in \Cal L$,
subject to certain relations
involving commutators of various kinds;
such a monomial is the class 
of 
$\alpha_1 \otimes_R \dots \otimes_R \alpha_m$
in $U(\Cal A,\Cal L)$.
The interpretation of the term \lq\lq monomial\rq\rq\ 
requires some care, though,
since for example when
$\alpha \in \Cal L$ has odd degree,
$\alpha^2$ is zero in
$U(\Cal A,\Cal L)$.
A more explicit description 
of the universal graded algebra
$U(\Cal A,\Cal L)$
for the special case
where $(\Cal A,\Cal L)$ is a crossed product Lie-Rinehart
algebra will be given below.
\smallskip
If $\Cal A= R$ with trivial ${\Cal L}$-action,
so that ${\Cal L}$ is just an ordinary
differential graded Lie algebra over $R$,
the object ${(U(R,\Cal L),\iota_{\Cal L},\iota_R)}$ is the 
{\it ordinary universal differential graded algebra\/} for ${\Cal L}$ 
(over $R$).
If $(\Cal A,\Cal L)$ is concentrated in degree zero,
that is, an ordinary (ungraded) Lie-Rinehart algebra,
the universal algebra
$U(\Cal A,\Cal L)$ comes down to the 
corresponding ordinary ungraded universal algebra;
an explicit description thereof may be found e.~g. in \cite\poiscoho,
\cite\rinehone.
When $A$ is the algebra of smooth functions and $L$ the Lie algebra
of smooth vector fields on a smooth manifold $M$,
$U(A,L)$ is the algebra of globally defined differential operators
on $M$.
\smallskip
It is obvious that,
for an arbitrary
differential graded Lie-Rinehart algebra
$(\Cal A,\Cal L)$,
there is a one--one
correspondence between
differential graded
(left) ${(\Cal A,\Cal L)}$-modules and 
differential  graded (left) $U(\Cal A,\Cal L)$-modules;
this correspondence is 
an equivalence of categories.
In particular, 
the obvious differential graded 
(left)
$(\Cal A,{\Cal L})$-module structure on $\Cal A$ mentioned above
turns $\Cal A$ into a differential graded left $U(\Cal A,{\Cal L})$-module;
the corresponding structure map is given by
$$
U(\Cal A,{\Cal L})
\otimes_R
\Cal A
\longrightarrow
\Cal A,
\quad
\alpha \otimes_R a \mapsto \alpha(a),
\tag6.1.5
$$
where ${\alpha \in {\Cal L},\, a \in \Cal A}$.
Next, let
$
\varepsilon
\colon
U(\Cal A,{\Cal L}) \longrightarrow \Cal A
$
be the 
morphism of
differential graded left $U(\Cal A,{\Cal L})$-modules
given by 
$$
\varepsilon (a) = a,\quad
\varepsilon (a \alpha) =0,\quad
\varepsilon (\alpha a) = \alpha (a).
\tag6.1.6
$$
It is {\it not\/}
a morphism  of differential graded algebras
unless ${\Cal L}$ acts trivially on $\Cal A$, and
its kernel
is the differential graded left ideal in $U(\Cal A,{\Cal L})$ 
generated by ${\Cal L}$.
In particular, the composite 
${\varepsilon \iota_{\Cal A}}$
is the identity map of $\Cal A$ whence
$\iota_{\Cal A}$ is injective.
Henceforth we will identify $\Cal A$ with its image in $U(\Cal A,{\Cal L})$, 
and we will
not distinguish 
in notation between the elements of $\Cal A$ and their images in
$U(\Cal A,{\Cal L})$.
Furthermore, it is clear that,
given 
two differential graded 
Lie-Rinehart algebras $(\Cal A_1,\Cal L_1)$ and
$(\Cal A_2,\Cal L_2)$,
a  morphism
$$
(\phi,\psi)
\colon
(\Cal A_1,\Cal L_1)
\longrightarrow
(\Cal A_2,\Cal L_2)
$$
of differential graded Lie-Rinehart algebras
induces a
 morphism
$$
U(\phi,\psi)
\colon
U(\Cal A_1,\Cal L_1)
\longrightarrow
U(\Cal A_2,\Cal L_2)
$$
of differential graded $R$-algebras.
\smallskip
Under the circumstances of (2.8),
the universal graded algebra
$\Cal U = U(\Cal A,\Cal L)$
(with zero differential)
may be obtained as follows:
The graded left $(A,L)$-module structure on
$\Cal A$
induces a
graded left $U(A,L)$-module structure on
$\Cal A$,
where 
$U(A,L)$ refers to the ordinary universal algebra
of $(A,L)$ mentioned above.
Let
$$
\Cal U = \Cal A \otimes_AU(A,L);
\tag6.2.1
$$
further, given
$a\in \Cal A$ and $u\in U(A,L)$,
identified in notation with
$a \otimes_A 1$ and $1 \otimes_A u$, respectively,
define the product
$au$ in the obvious way and let
$$
ua =au +u(a)
\tag6.2.2
$$
where $(u,a) \mapsto u(a)$
refers to the $U(A,L)$-action on $\Cal A$.
Since $L$
acts on $\Cal A$ by derivations,
this construction yields a
graded $R$-algebra structure 
on $\Cal U=\Cal A \otimes_AU(A,L)$
and,
together
with the obvious morphisms
$$
\iota_{\Cal L}\colon
\Cal L \to
\Cal U,
\quad
\iota_{\Cal A}
\colon
\Cal A
\to
\Cal U,
\tag6.2.3
$$
the graded $R$-algebra
$\Cal U$
is the universal graded algebra
$U(\Cal A,\Cal L)$
for
$(\Cal A,\Cal L)$.
Thus, 
since
as an $R$-algebra, $U(A,L)$
is generated by
the $a \in A$ and the $\alpha \in L$ subject to the relations
$$
\alpha \beta - \beta \alpha = [\alpha,\beta],
\quad
\alpha a - a \alpha = \alpha(a),
\tag6.2.4
$$
we see that,
as a graded $R$-algebra,
$U(\Cal A,\Cal L)$
is generated by the $a \in \Cal A$ and
the $\alpha \in L$ subject to the relations
$$
\alpha \beta - \beta \alpha = [\alpha,\beta],
\quad
\alpha a - a \alpha = \alpha(a),
\tag6.2.5
$$
for $a \in \Cal A$ and
$\alpha,\beta \in L$.
For clarity we point out that the non-trivial fact to be verified here
is that the algebra abstractly defined
by 
the generators $a \in \Cal A$ and $\alpha \in L$ and the relations
(6.2.5)
indeed admits the concrete description given by
(6.2.1) and (6.2.2).

\smallskip
Under the circumstances of (2.9),
given a left
$(A,L)$-module $M$ which, as an $A$-module, is
finitely generated and projective,
as a graded $U(A,L)$-module and
as a graded $\Cal A$-module,
the universal algebra
$\Cal U(\Cal A,\Cal L)$
may be written
$\roman{Alt}_A(M,U(A,L))$.
More precisely,
given 
$\phi \in \roman{Alt}_A(M,A)$ and $w \in U(A,L)$,
define
$\phi_w \in \roman{Alt}_A(M,U(A,L))$
by
$$
\phi_w(\xi_1,\dots,\xi_m)
=
(\phi(\xi_1,\dots,\xi_m)) w,
$$
where
$\xi_1,\dots,\xi_m \in M$.
It is manifest that
the canonical morphism
$$
\roman{Alt}_A(M,A) \otimes_AU(A,L)
@>>>
\roman{Alt}_A(M,U(A,L))
\tag6.2.6
$$
of graded left
$\roman{Alt}_A(M,A)$-modules and
right $U(A,L)$-modules
given by
the assignment to
$\phi \otimes_A w \in \roman{Alt}_A(M,A) \otimes_AU(A,L)$ of
$\phi_w \in \roman{Alt}_A(M,U(A,L))$
is an isomorphism.
We note that
the multiplication
in $\roman{Alt}_A(M,U(A,L))$
is {\it not\/}
given a shuffle map, though,
and additional terms of the kind spelled out in (4.2.4)
come into play;
the shuffle map would not even be well defined
since $U(A,L)$ is not an $A$-algebra.

\smallskip
Let $(A,L',L'')$ be a twilled Lie-Rinehart algebra
having $L'$ and $L''$ finitely generated and projective as $A$-modules.
By Theorem 3.2, the Lie-Rinehart differential $d''$
turns the graded crossed product Lie-Rinehart algebra
$$
(\Cal A'',\Cal L')
 = (\roman{Alt}^*_A(L'',A),\roman{Alt}^*_A(L'',L'))
\tag6.3.1
$$
into a differential graded Lie-Rinehart algebra.
Consider its universal differential graded
$R$-algebra
$U(\Cal A'',\Cal L';d'')$.
Its underlying graded $R$-algebra structure 
has been given in
(6.2.1) and (6.2.2) above and,
in view of (6.2.6)
the underlying graded $A$-module
may be written
in the form
$\roman{Alt}_A(L'',U(A,L'))$.
We now seek an explicit description of the differential $d''$.
\smallskip
Extend the left
$(A,L'')$-module structure 
{\rm (1.4.4)}
on $L'$
(which is part of the structure of a twilled
Lie-Rinehart algebra)
to a pairing
$$
\cdot\colon
L''
\otimes_R
U(A,L')
@>>>
U(A,L')
\tag6.3.2
$$
by means of the
recursive formula
$$
\aligned
\xi \cdot (x_1 \dots x_m)
&=
(\xi \cdot x_1) x_2 \dots x_m
\\
&\quad
-
(x_1 \cdot \xi) \cdot(x_2\dots x_m)
+ 
x_1 (\xi \cdot (x_2 \dots x_m))
\endaligned
\tag6.3.3
$$
where
$x_1,\dots, x_m \in L'$ and $\xi \in L''$,
where
$x_1 \dots x_m$ and
$x_2 \dots x_m$ 
refer to the corresponding elements
of $U(A,L')$
written out as monomials,
and
where
\linebreak
$(\xi \cdot x_1) x_2 \dots x_m$
and
$x_1 (\xi \cdot (x_2 \dots x_m))$
are the corresponding elements
of $U(A,L')$,
$\xi \cdot (x_2 \dots x_m)$
being supposed already defined.
We note that,
there is initially 
no need for the expression on the right-hand side of (6.3.3)
to be well defined
since the element $x_1\dots x_m$ of $U(A,L')$ depends on the order of 
the factors $x_1,\dots, x_m$.
When $(A,L',L'')$ is only a pre-twilled Lie-Rinehart algebra,
the right-hand side of (6.3.3)
will in general not be well defined.

\proclaim{Lemma 6.3}
For $x_1,\dots, x_m\in L'$ and $\xi \in L''$,
the differential $d''$ of the universal
differential graded algebra $U(\Cal A'',\Cal L';d'')$
satisfies the formula
$$
(d''(x_1\dots x_m))(\xi) = \xi \cdot(x_1\dots x_m)
$$
and is determined by it.
\endproclaim

Before proving the Lemma, we observe that,
given $x \in L'$ and $w \in U(A,L')$,
under the isomorphism (6.2.6),
$
(d''x) \otimes_Aw \in \roman{Alt}_A(L'',A) \otimes_AU(A,L')
$
manifestly goes to
$\phi_{x,w}
\in
\roman{Alt}_A(L'',U(A,L'))
$
which,
for $\xi \in L''$, is
defined by
$\phi_{x,w} (\xi) = (\xi \cdot x) w$ 

\demo{Proof}
Let $x_1,\dots, x_m\in L'$
and write $x = x_1$ and $w = x_2 \dots x_m$.
With reference to the description (6.2.1) and (6.2.2)
of 
the differential graded algebra $U(\Cal A'',\Cal L')$,
we have 
$$
d''(x_1\dots x_m) = d''(x w) =
(d''x) w +x (d''w).
$$
Let $\xi \in L''$.
In view of the observation made just after the statement of the Lemma,
$$
((d''x) w)(\xi) = (\xi \cdot x) w \in U(A,L').
$$
We now assert that
$$
(x (d''w))(\xi) =
-(x \cdot \xi) \cdot w + x (dw''(\xi)).
$$
In order to see this, we write $d''w = \sum_j \alpha_j w_j$,
for suitable 
$\alpha_j \in
\roman{Alt}_A(L'',A)$
and
$w_j \in U(A,L')$.
Now the product
$x d'' w(\xi)$ of $x$ and $d'' w(\xi)$
in $U(\Cal A'',\Cal L')$
may be written
$x d'' w(\xi)= \sum_j x \alpha_j (\xi) w_j$ 
whence
$$
\align
x d'' w(\xi)&= 
\sum_j x \alpha_j (\xi) w_j 
=
\sum_j \alpha_j (\xi) x w_j + \sum_j x (\alpha_j(\xi)) w_j,
\\
(d'' w)(x \cdot \xi)&= 
\sum_j \alpha_j (x \cdot \xi)w_j;
\endalign
$$
further, $(x (d'' w)) = \sum_j x\alpha_j w_j =
\sum_j (\alpha_j x  + x \cdot \alpha_j) w_j $
whence
$$
\align
(x (d'' w)) (\xi)&= 
\sum_j (\alpha_j (\xi) x  + (x \cdot \alpha_j)(\xi)) w_j 
\\
&=
\sum_j \alpha_j (\xi) x w_j + \sum_j x (\alpha_j(\xi)) w_j 
-\sum_j  \alpha_j(x \cdot\xi) w_j 
\\
&=
x ((d'' w)(\xi))
-
(d'' w)(x \cdot \xi).
\endalign
$$
Hence $(x (d''w))(\xi) =-(x \cdot \xi) \cdot w + x (dw''(\xi))$
as asserted.
Consequently
$$
((d''x) w)(\xi)=
((d''x) w)(\xi) +(x (d''w))(\xi)
=
(\xi \cdot x) w 
-(x \cdot \xi) \cdot w + x (dw''(\xi)).
$$
By induction on the length of monomials
(in the $x_j$'s) we may assume that
$
dw''(\xi)
=
\xi \cdot w
$
whence the assertion. \qed
\enddemo

Since 
for an arbitrary Lie-Rinehart algebra $(A,L)$ and an $A$-module
$M$, $L$ and $M$ being projective as $A$-modules,
the structure of a differential in
$\roman{Alt}_A(L,M)$ is equivalent to a left
$(A,L)$-module structure on $M$,
the Lemma entails at once the following.

\proclaim{Theorem 6.4}
For a twilled Lie-Rinehart algebra
$(A,L',L'')$ 
having $L'$ and $L''$ finitely generated and projective as $A$-modules,
the pairing
{\rm (6.3.2)} is a
left $(A,L'')$-module structure  on $U(A,L')$ which extends the
$(A,L'')$-module structure 
{\rm (1.4.4)}
on $L'$ (which is part of the structure of a twilled Lie-Rinehart algebra),
and the Lie-Rinehart differential $d''$ on 
$\roman{Alt}_A(L'',U(A,L'))$
with respect to this left $(A,L'')$-module structure  on $U(A,L')$
turns
$$
U(\Cal A'',\Cal L';d'')
=
(\roman{Alt}_A(L'',U(A,L')),d'')
$$
into a differential graded $R$-algebra
in such a way that
$(U(\Cal A'',\Cal L'),d'')$
is the universal differential graded algebra
for
the differential graded crossed product Lie-Rinehart algebra
$(\Cal A'',\Cal L';d'')$. \qed
\endproclaim

In particular, we see that, when
$(A,L',L'')$
is a twilled Lie-Rinehart algebra,
the expression on the right-hand side of
(6.3.3) {\it is\/} well defined.
Here and below
we will write
$U(\Cal A'',\Cal L';d'')$
for the universal differential graded algebra;
when only the underlying universal graded algebra
is understood
(i.~e. when differentials are ignored),
we write
$U(\Cal A'',\Cal L')$.
We note again that,
for reasons explained before,
the multiplication in
$\roman{Alt}_A(L'',U(A,L'))$
is {\it not\/} a shuffle product.
It is clear that,
with the roles of $L'$ and $L''$ interchanged, 
the statement of the theorem is true as well.

\smallskip
The
left $(A,L'')$-module structure 
{\rm (6.3.2)}  on $U(A,L')$ 
may be explained in another way which we now explain briefly:
Let $L= L'' \bowtie L'$ be the twilled sum of $L'$ and $L''$.
From the Poincar\'e-Birkhoff-Witt theorem for
$U(A,L)$, cf. \cite\rinehone,
we deduce at once that
$U(A,L)$
may be written as the tensor product over $A$
of $U(A,L')$  and $U(A,L'')$.
As an $R$-algebra,
$U(A,L)$ is what might be called the {\it twilled product\/}
of $U(A,L'')$  and $U(A,L')$
but we do not explain this here since we will not need it.
Let $\varepsilon'' \colon U(A,L'') \to A$
be the corresponding morphism of left 
$U(A,L'')$-modules introduced
(in somewhat greater generality) in
(6.1.6) above.
Exploiting the fact that, as an $A$-module,
$U(A,L)$ is generated by monomials
$x_1x_2 \dots x_m \xi_1 \xi_2 \dots \xi_{\ell}$,
where
$x_1,x_2, \dots, x_m \in L'$ and
$\xi_1, \xi_2, \dots, \xi_{\ell}\in L''$,
we extend
$\varepsilon''$
to an $A$-module surjection
$$
\widetilde{\varepsilon}'' \colon
U(A,L)
@>>>
U(A,L')
\tag6.3.4
$$
by means of the assignments
$$
\widetilde{\varepsilon}''(x_1x_2 \dots x_m \xi_1 \xi_2 \dots \xi_{\ell}) =
\cases  x_1x_2 \dots x_m,  \quad & \ell = 0,
\\ 0, & \ell \geq 1.
\endcases
$$
When $U(A,L)$  is written as the tensor product
$U(A,L')\otimes_AU(A,L'')$,
$\widetilde{\varepsilon}''$
takes the form
$\roman {Id}\otimes_A\varepsilon''$.
The Poincar\'e-Birkhoff-Witt theorem 
for Lie-Rinehart algebras 
implies that $\widetilde{\varepsilon}''$ is well defined.
Moreover,
the kernel of
$\widetilde{\varepsilon}''$ is the left ideal of
$U(A,L)$
generated by
$L''$ whence the assignment
$$
v \cdot w=\widetilde{\varepsilon}'' (v  w),
\quad v \in U(A,L),\ w \in U(A,L'),
$$
$v  w$ being the product in $U(A,L)$,
endows $U(A,L')$ with a left
$(A,L)$-module structure
$$
\cdot\colon
U(A,L)
\otimes_R
U(A,L')
@>>>
U(A,L')
\tag6.3.5
$$
in such a way that,
$U(A,L)$ being endowed with its obvious
left
$U(A,L)$-module structure,
$\widetilde{\varepsilon}''$
is a surjective morphism
of $(A,L)$-modules.
The description of
$\widetilde{\varepsilon}''$
as the tensor product $\roman {Id}\otimes_A\varepsilon''$
is {\it not\/}
compatible with
the $(A,L)$-module structures, though.
For $w \in U(A,L')$ and $\xi \in L''$,
the action (6.3.5) plainly satisfies
$$
\xi \cdot w = \widetilde{\varepsilon}''[\xi,w],
$$
and the
$(A,L'')$-module structure
(6.3.2)
on $U(A,L')$ resulting from restriction
and
(6.3.2)
may as well be described as the composite
$$
L''
\otimes_R
U(A,L')
@>{\iota_{L''} \otimes \roman{Id}}>>
U(A,L'')
\otimes_R
U(A,L')
@>{[\cdot,\cdot]}>>
U(A,L) 
@>{\widetilde{\varepsilon}''}>>
U(A,L') 
\tag6.3.6
$$
where
$[\cdot,\cdot] \colon U(A,L'') \otimes_R U(A,L') \to U(A,L)$
refers to the restriction of the commutator bracket
on 
$U(A,L) \otimes_R U(A,L)$
to $U(A,L'') \otimes_R U(A,L')$
(viewed as a subspace of
$U(A,L) \otimes_R U(A,L)$).
It may then readily be seen directly
that (6.3.6)
yields a left $(A,L'')$-module structure on $U(A,L')$:
Since the commutator bracket in $U(A,L)$
satisfies the Jacobi identity, for $\xi ,\eta \in L''$ and
$w \in U(A,L')$, we plainly have
$$
[\xi,\eta] \cdot w =  \xi\cdot (\eta\cdot w) - \eta\cdot (\xi \cdot w ) .
$$
Furthermore,
given $a \in A$,  $\xi \in L''$, and
$w \in U(A,L')$,
$$
\align
(a \xi)\cdot w &=
\widetilde{\varepsilon}''(a \xi w  -w a \xi)
=
\widetilde{\varepsilon}'' (a \xi w )
=
a\widetilde{\varepsilon}'' (\xi w )
= a(\xi\cdot w)
\\
\xi\cdot (aw) &=
\widetilde{\varepsilon}'' (\xi aw - aw \xi)
=
\widetilde{\varepsilon}'' ((a\xi +\xi(a))w - aw \xi)
=
\xi(a) w + a (\xi\cdot w)
\endalign
$$
whence (6.3.6) yields
a left
$(A,L'')$-module structure 
on
$U(A,L')$,
and this structure
plainly extends the
left
$(A,L'')$-module structure 
{\rm (1.4.4)}
on $L'$.

\smallskip
We now proceed towards the differential homological algebra
interpretation of twilled Lie-Rinehart algebras
alluded to earlier.
To begin with, let $(A,L)$ be an (ordinary ungraded) Lie-Rinehart algebra.
Consider the graded
left $U(A,L)$-module
$U(A,L)\otimes_A\Lambda_AL$
where $A$ acts on the right of $U(A,L)$ by means of the canonical
map $\iota_A \colon A \longrightarrow U(A,L)$.
For $u \in U(A,L)$ and $\alpha_1,\dots,\alpha_n \in L$, let
$$
\aligned
d&(u\otimes_A (\alpha_1\wedge\ldots\wedge \alpha_n))
\\
&=
\quad\sum_{i=1}^n (-1)^{(i-1)}u\alpha_i\otimes_A
(\alpha_1\wedge\dots\widehat{\alpha_i}\ldots\wedge \alpha_n)
\\
&\phantom{=}+\quad
\sum_{j<k} (-1)^{(j+k)}u\otimes_A 
 (\lbrack \alpha_j,\alpha_k \rbrack \wedge
\alpha_1\wedge\dots\widehat{\alpha_j}\dots
\widehat{\alpha_k}\ldots \wedge\alpha_n),
\endaligned
\tag6.5
$$
Rinehart \cite\rinehone\ has shown
that this yields
an 
$U(A,L)$-linear differential
$$
d
\colon
U(A,L) \otimes _A \Lambda_AL 
\longrightarrow
U(A,L) \otimes _A \Lambda_AL,
\tag6.6
$$
that is, $dd$ is zero.
The non-trivial fact
to be verified here is that
the operator $d$ is well defined.
We will refer to
$$
K(A,L) = (U(A,L) \otimes _A \Lambda_AL,d)
\tag6.7
$$
as the
{\it Rinehart complex\/} for $(A,L)$.
Rinehart has also proved that,
when $L$
is projective as an $A$-module,
$K(A,L)$ is a projective resolution of $A$ in the category
of left $U(A,L)$-modules.
This resolution generalizes the 
Koszul resolution of the ground ring 
(or ground field) in ordinary Lie algebra cohomology.

\smallskip
Let now $\Cal A$ be a graded commutative $A$-algebra
which is endowed with a graded left $(A,L)$-module structure
on $\Cal A$
in such a way 
that the underlying $L$-action on $\Cal A$ is by graded derivations
and that the canonical map from $A$ to $\Cal A$ is a morphism
of $(A,L)$-modules.
Consider the graded crossed product Lie-Rinehart algebra
$(\Cal A,\Cal L)$
introduced in (2.8) above.
Let
$$
\Cal K = K(\Cal A,\Cal L) = \Cal A \otimes _A K(A,L).
\tag6.8
$$
Plainly, $\Cal K$ is of the form
$$
\Cal K=
(\Cal A\otimes_A U(A,L) \otimes _A \Lambda_AL,d)
=
(U(\Cal A,\Cal L) \otimes _A \Lambda_AL,d).
$$
With the obvious induced left 
$U(\Cal A,\Cal L)$-module structure,
$$
\Cal K_n @>d>>
\Cal K_{n-1}
@>d>>
\dots
@>d>>
\Cal K_1 @>d>>
U(\Cal A,\Cal L)
@>\varepsilon>>
\Cal A
\tag6.9
$$
is then an exact sequence of graded left 
$U(\Cal A,\Cal L)$-modules.
Thus, when $\Cal A$ is projective as an $A$-module,
(6.8) yields a projective resolution
of $\Cal A$ in the category of graded left
$U(\Cal A,\Cal L)$-modules.
Hence we may then use this resolution to compute,
for any graded left $(\Cal A,\Cal L)$-module $\Cal M$,
the {\it cohomology\/}
$$
\roman H^*(\Cal L,\Cal M) 
= \roman{Ext}^*_{U(\Cal A,\Cal L)}(\Cal A,\Cal M)
\tag6.10.1
$$
{\it of\/} $\Cal L$ {\it with coefficients in\/} 
$\Cal M$ 
and, 
for any graded right $(\Cal A,\Cal L)$-module $\Cal N$,
the {\it homology\/}
$$
\roman H_*(\Cal L,\Cal N) 
= \roman{Tor}_*^{U(\Cal A,\Cal L)}(\Cal N,\Cal A)
\tag6.10.2
$$
{\it of\/} $\Cal L$ {\it with coefficients in\/} 
$\Cal N$. 
Since the resolution (6.8) has the form
$\Cal A \otimes_AK(A,L)$,
we see at once that
the canonical isomorphisms
$$
\roman{Hom}_{U(A,L)}(K(A,L),\Cal M)
@>>>
\roman{Hom}_{U(\Cal A,\Cal L)}(K(\Cal A,\Cal L),\Cal M)
$$
and
$$
\Cal N \otimes_{U( A, L)} K( A, L)
@>>>
\Cal N \otimes_{U(\Cal A,\Cal L)} K(\Cal A,\Cal L)
$$
induces isomorphisms
$$
\roman H^*(L,\Cal M)
@>>>
\roman H^*(\Cal L,\Cal M)
\quad 
\text{and}
\quad
\roman H_*(L,\Cal N)
@>>>
\roman H_*(\Cal L,\Cal N) 
\tag6.11
$$
of $R$-modules.
Thus the homology and cohomology of
$\Cal L$ boil down to
the homology and cohomology of $L$ with graded coefficients.
\smallskip
Let now
$(A,L',L'')$ 
be a twilled Lie-Rinehart algebra
having $L'$ and $L''$ finitely generated and projective as $A$-modules.
Consider
the differential graded crossed product Lie-Rinehart algebra
$$
(\Cal A'',\Cal L';d'')
 = (\roman{Alt}^*_A(L'',A),\roman{Alt}^*_A(L'',L');d'');
$$
cf. Theorem 3.2 above.
The corresponding Rinehart complex
$K(\Cal A'',\Cal L';d'')$,
cf. (6.8) above,
is plainly of the form
$$
(\roman{Alt}^*_A(L'',K(A,L')),d')
=
(\roman{Alt}^*_A(L'',U(A,L')\otimes_A\Lambda_A L'),d'),
\tag6.12
$$
that is, may be written
$$
\dots
@>>>
K_j(\Cal A'',\Cal L')
@>{d'}>>
K_{j-1}(\Cal A'',\Cal L')
@>{d'}>>
\dots
@>{d'}>>
K_1(\Cal A'',\Cal L')
@>{d'}>>
K_0(\Cal A'',\Cal L')
\tag6.13
$$
where for $j \geq 0$, 
$K_j(\Cal A'',\Cal L')=\roman{Alt}_A(L'', K_j(A, L'))$,
and the latter, in turn, is isomorphic to
$\roman{Alt}_A(L'', U(A,L')\otimes_A\Lambda^j_A L')$;
as additional piece of structure, each $K_j(\Cal A'',\Cal L')$
now also carries the differential $d''$
(with respect to the Lie-Rinehart structure on
$(A,L'')$ and the left
$(A,L'')$-module structure (6.3.2) on $U(A,L')$).
Plainly,
each $\roman{Alt}_A(L'', U(A,L')\otimes_A\Lambda^j_A L')$
is isomorphic to
$\roman{Alt}_A(L'', U(A,L'))\otimes_A\Lambda^j_A L'$ which,
in turn, 
in view of the structure of
$U(\Cal A'',\Cal L';d'')$
elucidated in (6.4),
is just a rewrite of
$U(\Cal A'',\Cal L';d'') \otimes _A\Lambda^j_A L'$.
The complex (6.13) 
(in the category of $R$-chain complexes)
is a proper projective resolution
(cf. e.~g. \cite\maclaboo\ (XII.11) p. 397 for the notion of a proper
projective resolution)
of
$\Cal A''=(\roman{Alt}_A(L'',A),d'')$
in the category of differential graded left
$U(\Cal A'',\Cal L';d'')$-modules.
Thus, from (6.13) we may compute,
for any {\it differential graded left\/} 
$(\Cal A'',\Cal L')$-{\it module\/} $\Cal M$,
the {\it cohomology\/}
$$
\roman H^*(\Cal L',\Cal M) 
= \roman{Ext}^*_{U(\Cal A'',\Cal L';d'')}(\Cal A'',\Cal M)
\tag6.14.1
$$
{\it of\/} $\Cal L'$ {\it with coefficients in\/} 
$\Cal M$ 
and, 
for any {\it differential graded right\/} 
$(\Cal A'',\Cal L')$-{\it module\/} $\Cal N$,
the {\it homology\/}
$$
\roman H_*(\Cal L',\Cal N) 
= \roman{Tor}_*^{U(\Cal A'',\Cal L';d'')}(\Cal N,\Cal A'')
\tag6.14.2
$$
{\it of\/} $\Cal L''$ {\it with coefficients in\/} 
$\Cal N$.

\proclaim{Theorem 6.15}
For a twilled Lie-Rinehart algebra
$(A,L',L'')$ 
having $L'$ and $L''$ finitely generated and projective as $A$-modules,
the graded $R$-modules
$\roman H^*(\Cal L',\Cal A'')$
and
$\roman H^*(L,A)$
(=$\roman{Ext}^*_{U(A,L)}(A,A)$)
are canonically isomorphic.
\endproclaim

It is clear that the statement of the theorem holds as well
with the roles of $L'$ and $L''$ interchanged.

\demo{Proof of the Theorem}
Consider
the
bicomplex 
$$
\roman{Hom}^*_{U(\Cal A'',\Cal L';d'')}(K(\Cal A'',\Cal L';d''),\Cal A'')
\tag6.15.1
$$
computing
$\roman{Ext}^*_{U(\Cal A'',\Cal L';d'')}(\Cal A'',\Cal A'')$.
The bigraded $A$-module
which underlies (6.15.1) has the form
$$
\roman{Hom}_{\Cal A''}(\roman{Alt}_A(L'', \Lambda_A L'),\Cal A'')
\cong
\roman{Hom}_{\Cal A''}(\Cal A''\otimes_A\Lambda_A L',\Cal A'')
$$
and this is clearly canonically isomorphic to
$$
\roman{Hom}_A(\Lambda_A L',\Cal A'')
\cong
\roman{Alt}_A(L', \Cal A'')
\cong
\roman{Alt}_A(L',\roman{Alt}_A(L'',A))
\cong
\roman{Alt}_A(L,A).
\tag6.15.2
$$
The operator
$d'$ on
$\roman{Hom}_{\Cal A''}(\roman{Alt}_A(L'', \Lambda_A L'),\Cal A'')$
plainly amounts to the
Lie-Rinehart differential
$d'$ on
$\roman{Alt}_A(L', \Cal A'')$
with reference to the Lie-Rinehart structure
on $(A,L')$ and the graded
left
$(A,L')$-module structure on $\Cal A''$;
this
Lie-Rinehart differential, in turn,
corresponds to the operator
on $\roman{Alt}_A(L,A)$
denoted by the same symbol.
By construction,
the 
operator
$d''$ on
$\roman{Hom}_{\Cal A''}(\roman{Alt}_A(L'', \Lambda_A L'),\Cal A'')$
is compatible with the operator $d'$.
Moreover, since
$K(\Cal A'',\Cal L')$ has the form
$\roman{Alt}_A(L'', K(A,L'))$,
the 
bigraded $A$-module underlying the
bicomplex
(6.15.1)
may as well be written
$$
\roman{Alt}_A(L'', \roman{Hom}_{U(A,L')}(K(A,L'),A))
\cong
\roman{Alt}_A(L'', \roman{Alt}_A(L',A))
\tag6.15.3
$$
and, using this description,
we see that the 
operator
$d''$ on 
(6.15.3),
that is, the resulting operator on
$\roman{Hom}_{\Cal A''}(\roman{Alt}_A(L'', \Lambda_A L'),\Cal A'')$,
amounts to the
Lie-Rinehart differential
$d''$ on
$\roman{Alt}_A(L'', \Cal A')$
with reference to the Lie-Rinehart structure
on $(A,L'')$ and the graded
left
$(A,L'')$-module structure on $\Cal A' =\roman{Alt}_A(L',A)$;
thus
the operator
$d''$ on 
(6.15.3)
corresponds to the operator
on $\roman{Alt}_A(L,A)$
denoted by the same symbol.
Consequently
the total differential on
$\roman{Hom}_{\Cal A''}(\roman{Alt}_A(L'', \Lambda_A L'),\Cal A'')$
amounts 
to the total differential on
$\roman{Alt}_A(L,A)$
arising from the bicomplex (1.15.4). \qed
\enddemo

The theorem provides an interpretation of the bicomplex
(1.15.4):
The object in the middle of
(6.15.1)
is just the bigraded $A$-module underlying
(1.15.4), with the roles of
$L'$ and $L''$ interchanged.
Thus, the differential $d'$ 
(in (1.15.4))
alone computes the graded
$\roman{Ext}^*_{U(\Cal A'',\Cal L')}(\Cal A'',\Cal A'')$,
the differential $d''$ on every object in sight being ignored.
However, the compatibility between $d'$ and $d''$
entails that
(6.13)
is a resolution in the differential graded category,
and the bicomplex
computing 
$\roman{Ext}^*_{U(\Cal A'',\Cal L';d'')}(\Cal A'',\Cal A'')$
boils down to (1.15.4).
This provides, in particular, a new interpretation of the
Dolbeault complex.
Rinehart has shown that the ordinary de Rham cohomology groups
may be written as Ext-groups over the algebra of differential
operators.
Theorem 6.15
includes the corresponding result for the Dolbeault cohomology
groups,
which now appear as differential graded Ext-groups.

\medskip\noindent{\bf 7. Duality and generators of dBV algebras}
\smallskip\noindent
Let $(A,L)$ be a Lie-Rinehart algebra. 
In \cite\bv\ we have shown that an exact generator
of a Gerstenhaber algebra of the kind
$\Lambda_AL$ yields precisely the
differential in
the standard complex computing
the {\it homology\/} of the Lie-Rinehart algebra $(A,L)$
with values in $A$, endowed with a
right $(A,L)$-module structure corresponding to the generator.
This relies on a notion of homological duality.
Our present aim is to generalize this 
notion and the
relationship between exact generators and differentials
in the standard complex to the 
differential graded
setting.
This will 
give conceptual explanations of some of the results in earlier sections
and will
elucidate the at first somewhat mysterious concept of generator
of a differential bigraded Gerstenhaber bracket.
\smallskip
Let $(A,L)$ be a Lie-Rinehart algebra.
Recall from Section 2 above that a (graded) left $L$-module structure
$L\otimes_R M \to M$
on a graded $A$-module $M$, written
$(\alpha,x) \mapsto \alpha(x)$,
is called a {\it left\/}
$(A,L)$-module structure provided
$$
\align
\alpha (ax) &= \alpha(a) x + a \alpha(x),
\tag7.1.1
\\
(a\alpha)(x) &= a(\alpha (x)),
\tag7.1.2
\endalign
$$
where $a \in A,\ x \in M, \ \alpha \in L$.
More generally,
cf. (1.10),
such an assignment
$L\otimes_R M \to M$,
not necessarily a left $L$-module structure
but still satisfying (7.1.1) and (7.1.2),
is referred to as an
$(A,L)$-{\it connection\/}, cf. \cite\poiscoho,
or, somewhat more precisely,
{\it left\/} $(A,L)$-{\it connection\/};
in this language,
a (graded) left $(A,L)$-module structure
is a (graded left)
{\it flat\/} $(A,L)$-{\it connection\/}.
Likewise, let $N$ be a graded $A$-module,
and let there be given an assignment
$N\otimes_R L \to N$,
written
$(x,\alpha) \mapsto x \circ \alpha$
or, somewhat simpler,
$(x,\alpha) \mapsto x \alpha$
(when there is no risk of confusion);
it is
called a {\it (graded) right\/}
$(A,L)$-module structure provided
it is a (graded) right 
$L$-module structure
and, moreover, satisfies
$$
\align
(ax)\alpha &= a(x \alpha) - (\alpha(a))x,
\tag7.2.1
\\
x(a\alpha) &= a(x \alpha) - (\alpha(a))x,
\tag7.2.2
\endalign
$$
where $a \in A,\ x \in N, \ \alpha \in L$;
an assignment 
$N\otimes_R L \to N$
of this kind
is referred to as a {\it (graded) right\/}
$(A,L)$-{\it connection\/} provided 
it satisfies only (7.2.1) and (7.2.2)
without
necessarily being a (graded) right $L$-module structure.
A (graded) right $(A,L)$-module structure
is also said to be a (graded)
{\it flat right\/} $(A,L)$-{\it connection\/}.
Graded left- and right $(A,L)$-modules
correspond to
graded left- and right $U(A,L)$-modules,
and vice versa.
More generally,
graded left- and right $(A,L)$-connections may be shown to correspond
bijectively to
graded left- and right $U(A,E)$-module structures,
for suitable $(R,A)$-Lie algebras $E$
mapping surjectively onto $L$.
\smallskip
Return to the situation at the beginning of (5.3).
Thus $\Cal A$ is a graded commutative $A$-algebra
together with an $L$-action
$L \otimes_R\Cal A \to \Cal A$ by derivations
such that the canonical map from $A$ to $\Cal A$ is
a morphism of left 
$(A,L)$- and hence
$U(A,L)$-modules.
A (left) $(A,L)$-connection
$$
L\otimes_R \Cal M @>>> \Cal M
\tag7.3
$$
on an induced (graded) $\Cal A$-module
of the kind
$\Cal M=\Cal A \otimes _A M$ 
where $M$ is an $A$-module
will be said to be {\it compatible\/}
(with the $L$-action on $\Cal A$)
provided
$$
\alpha(a x) = \alpha (a) x + a \alpha (x),
\quad
\alpha \in L,
\ a \in \Cal A,
\ x \in M
\tag7.4
$$
where
$(\alpha, a) \mapsto \alpha (a)$
refers to the graded left $(A,L)$-module structure
on $\Cal A$ and, accordingly, we define
a {\it compatible\/}
(left) $(A,L)$-module structure
on an induced (graded) $\Cal A$-module.
Furthermore, a right $(A,L)$-connection
$\Cal M\otimes_R L\to \Cal M$
on an induced (graded) $\Cal A$-module
$\Cal M=\Cal A \otimes _A M$, 
$M$ being an $A$-module,
will be said to be {\it compatible\/}
(with the $L$-action on $\Cal A$)
provided
$$
(ax)\alpha = a(x \alpha) - (\alpha(a))x,
\quad
\alpha \in L,
\ a \in \Cal A,
\ x \in M
\tag7.5
$$
where
$(\alpha, a) \mapsto \alpha (a)$
still refers to the graded left $(A,L)$-module structure
on $\Cal A$ and
we can  accordingly talk about
a {\it compatible\/}
right $(A,L)$-module structure
on an induced (graded) $\Cal A$-module.
It is clear that,
given  a left $A$-module $M$,
any compatible left or right $(A,L)$-connection
or
left or right $(A,L)$-module structure
on an induced module
of the kind
$\Cal A \otimes _A M$,
is determined by its restriction to $M$.
\smallskip
Theorem 1 of \cite\bv\  now generalizes in the following way.

\proclaim {Theorem 7.6}
There is a bijective correspondence between
right $(A,L)$-connections on 
$A$ and hence compatible right $(A,L)$-connections on
$\Cal A$
and $R$-linear operators generating the bigraded
Gerstenhaber bracket on
$\Lambda_{\Cal A}\Cal L$.
Under this correspondence,
flat 
right $(A,L)$-connections on 
$A$ and hence compatible ones on
$\Cal A$,
that is,
compatible right $(A,L)$-module structures
on $\Cal A$,
correspond to exact generators.
More precisely:
Given an $R$-linear operator $D$ generating the bigraded Gerstenhaber 
bracket on 
$\Lambda_{\Cal A}\Cal L$
the formula
$$
a \circ \alpha = a(D\alpha) - \alpha(a),\quad
a \in \Cal A,\ \alpha \in L,
\tag7.6.1
$$
defines
a right
$(A,L)$-connection on $\Cal A$,
and this
right
$(A,L)$-connection on $\Cal A$
is determined by
$$
a \circ \alpha = a(D\alpha) - \alpha(a),\quad
a \in A,\ \alpha \in L.
\tag7.6.2
$$
Conversely, given a
compatible right $(A,L)$-connection 
$(a,\alpha) \mapsto a \circ \alpha$
on $\Cal A$ 
\linebreak
($a \in \Cal A, \alpha \in L)$,
the operator
$D$ on $\Lambda_{\Cal A}\Cal L$ defined by means of the formula 
$$
\aligned
D (a\alpha_1 \wedge \ldots \wedge \alpha_n)
&=
\quad\sum_{i=1}^n (-1)^{(i-1)}(a\circ \alpha_i)
\alpha_1 \wedge \ldots\widehat{\alpha_i}\ldots \wedge \alpha_n 
\\
&\quad
+
\sum_{j<k} (-1)^{(j+k)}
a \lbrack \alpha_j,\alpha_k \rbrack \wedge
\alpha_1 \ldots\widehat{\alpha_j}\ldots
\widehat{\alpha_k}\ldots \wedge \alpha_n,
\endaligned
\tag7.6.3
$$
where
$a \in \Cal A$
and
$\alpha_1, \ldots, \alpha_n \in L$,
yields an $R$-linear operator $D$ generating the bigraded Gerstenhaber 
bracket on 
$\Lambda_{\Cal A}\Cal L$.
\endproclaim

\demo{Proof}
The argument given
for the proof of Theorem 1
in \cite\bv\ 
carries readily over.
In fact,
compatible right $(A,L)$-connections on
$\Cal A$ correspond
bijectively to 
right $(A,L)$-connections on 
$A$ and,
these, in turn,
correspond bijectively
to 
generators for the ordinary Gerstenhaber algebra
$\Lambda_AL$,
by virtue of 
Theorem 1
in \cite\bv\.
Furthermore,
generators for the ordinary Gerstenhaber algebra
$\Lambda_AL$
correspond bijectively
to generators
for the bigraded Gerstenhaber algebra
$\Lambda_{\Cal A}\Cal L$. \qed
\enddemo

We note that, when $\Cal A$ is just $A$,
the statement of the present Theorem 7.6
boils down verbatim to
Theorem 1
in \cite\bv. 
\smallskip
Given an exact generator 
$\partial$
for the bigraded Gerstenhaber
algebra 
$\Lambda_{\Cal A}\Cal L$,
we will acccordingly write
$\Cal A_{\partial}$
for $\Cal A$
together with the graded right
$(A,L)$-module structure
given by
(7.6.1).
Theorem 2 of \cite\bv\ 
now extends as follows
where $K(A,L)$ refers to the {\it Rinehart complex\/}  for
$(A,L)$ 
(reproduced as (6.7) above).

\proclaim{Theorem 7.7}
Given an exact generator 
$\partial$
for the bigraded Gerstenhaber
algebra 
$\Lambda_{\Cal A}\Cal L$,
the chain
complex underlying the bigraded Batalin-Vilkovisky algebra
$(\Lambda_{\Cal A}\Cal L,\partial)$
coincides with the chain complex
$({\Cal A}_{\partial}\otimes_{U(A,L)}K(A,L),d)$.
In particular,
when $L$ is projective as an $A$-module,
the bigraded Batalin-Vilkovisky algebra
$(\Lambda_{\Cal A}\Cal L,\partial)$
computes
$$
\roman H_*(L,\Cal A_\partial) 
\left(= \roman{Tor}_*^{U(A,L)}(\Cal A_\partial,A)\right),
$$
the homology of $L$ with coefficients in $\Cal A_\partial$. 
\endproclaim

\demo{Proof}
The argument given
for the proof of Theorem 2
in \cite\bv\ 
carries  readily over.
Details are left to the reader. \qed
\enddemo

When $\Cal A$ is just $A$,
the statement of the present Theorem 7.7
come down verbatim to
Theorem 2
in \cite\bv. 

\smallskip\noindent
{\smc Remark 7.8.}
The two above theorems
{\it reveal\/} the significance
of a generator of the Gerstenhaber bracket
of a Gerstenhaber algebra of the kind $\Lambda_{\Cal A} \Cal L$
(and in particular of the kind $\Lambda_AL$):
Indeed,
the defining property (5.1)
of an exact generator
{\it precisely incorporates a description
of the Lie-Rinehart
differential in the corresponding
complex\/}
$({\Cal A}_{\partial}\otimes_{U(A,L)}K(A,L),d)$
{\it in terms of the Lie bracket in\/} $L$
and the {\it corresponding right\/} $(A,L)$-module structure 
on $\Cal A$.
We will show below that the extension of this observation
to the differential graded setting provides a conceptual
explanation of the isomorphism
(5.4.6.1).
For a generator $\Delta$ which is not necessarily exact,
a corresponding remark can still be made:
The 
corresponding right
$(A,L)$-connection
on
${\Cal A}$
still induces an operator
on ${\Cal A}\otimes_{U(A,L)}K(A,L)$,
and the defining property (5.1)
merely yields a description
thereof
in terms of the Lie-Rinehart
structure and the right
$(A,L)$-connection on $\Cal A$.
\smallskip

The above considerations entail that,
under appropriate circumstances,
generators for bigraded Gerstenhaber brackets exist.
In order to explain this,
we suppose that, as an $A$-module, $L$ is finitely generated and projective
of finite constant rank $n$
and, for intelligibility,
recall the following,
spelled out in \cite\bv\ as Theorem 3.

\proclaim{Proposition 7.9}
There is a bijective
correspondence between
$(A,L)$-connections on $\Lambda_A^nL$
and
right
$(A,L)$-connections
on $A$.
Under this correspondence,
left $(A,L)$-module structures on $\Lambda_A^nL$
(i.~e. flat connections)
correspond to
right
$(A,L)$-module structures
on $A$.
More precisely: Given an
$(A,L)$-connection $\nabla$ on $\Lambda_A^nL$,
the negative of the (generalized) Lie-derivative
on
$A \cong \roman {Hom}_A(\Lambda_A^n L, M)
$
with reference to the
connection $\nabla$ on $M=\Lambda_A^n L$, that is, the formula
$$
(\phi \alpha) x = 
\phi(\alpha x) -\nabla_{\alpha} (\phi(x)),
\tag7.9.1
$$
where $x \in \Lambda_A^n L,\ \alpha \in L,
\ \phi \in \roman {Hom}_A(\Lambda_A^n L, \Lambda_A^n L) \cong A$,
yields 
a right $(A,L)$-connection on $A$.
Conversely,
given a
a right $(A,L)$-connection on $A$
(written
$(a,\alpha) \mapsto a \alpha$),
on $\Lambda_A^nL \cong \roman {Hom}_A(C_L, A)$,
the assignment
$$
L \otimes _R \Lambda_A^nL \to \Lambda_A^nL,
\quad
(\alpha, \psi) \mapsto \nabla_\alpha \psi
$$
where
$$
(\nabla_\alpha \psi) x
= \psi (x\alpha)
- (\psi x) \alpha,
\quad
x \in C_L,\ \alpha \in L,
\ \psi \in \roman {Hom}_A(C_L, A),
\tag7.9.2
$$
yields an $(A,L)$-connection $\nabla$.
\endproclaim

Combining Theorem 7.6 with Proposition 7.9, we obtain:

\proclaim{Corollary 7.10}
There is a bijective correspondence between
$(A,L)$-connections
on
$\Lambda_A^n L$
and linear operators $D$ generating the Gerstenhaber 
bracket on 
$\Lambda_{\Cal A} \Cal L$.
Under this correspondence, flat connections correspond to 
exact generators, that is, to differentials.
The relationship is made explicit
by means of
{\rm (7.6.1), (7.6.3), (7.9.1)} and
{\rm (7.9.2)}. \qed
\endproclaim

The relationship 
in terms of
{\rm (7.6.1), (7.6.3), (7.9.1)} and
{\rm (7.9.2)}
comes down precisely to that spelled out in (5.3.7),
which involves (5.3.2) and (5.3.6).
This observation provides another proof of Theorem 5.3.7.

\proclaim{Corollary 7.11}
As an $A$-module, $L$ being finitely generated and projective
of finite constant rank $n$,
the
bracket on the bigraded crossed product Gerstenhaber algebra
$\Lambda_{\Cal A} \Cal L$
always has a generator,
and these generators are classified by
$(A,L)$-connections
on
$\Lambda_A^n L$.
\endproclaim

Thus, under the circumstances of Theorem 7.7,
since $\Cal L$ is just the graded crossed product extension
of $\Cal A$ and $L$,
the (graded) right $(A,L)$-module structure
on $\Cal A$ extends to a
graded right $(\Cal A,\Cal L)$-module structure
on $\Cal A$,
and
the bigraded Batalin-Vilkovisky algebra
$(\Lambda_{\Cal A}\Cal L,\partial)$
computes as well
the {\it homology\/}
$$
\roman H_*(\Cal L,\Cal A_\partial) 
= \roman{Tor}_*^{U(\Cal A,\Cal L)}(\Cal A_\partial,\Cal A)
$$
{\it of\/} $\Cal L$ {\it with coefficients in\/} 
$\Cal A_\partial$
where
$U(\Cal A,\Cal L)$ 
is the
universal algebra
for the graded Lie-Rinehart algebra
$(\Cal A,\Cal L)$
explained in the previous Section.
\smallskip
With these preparations out of the way,
let $(A,L',L'')$ 
be a twilled Lie-Rinehart algebra
having $L'$ and $L''$ finitely generated and projective as $A$-modules,
and consider
the corresponding differential graded crossed product Lie-Rinehart algebra
$$
(\Cal A'',\Cal L';d'')
 = (\roman{Alt}^*_A(L'',A),\roman{Alt}^*_A(L'',L');d'').
$$
Suppose that, as an $A$-module,
$L'$ (being finitely generated and projective)
is of finite constant rank $n$ (say).
The statements of Theorem 7.7,
Theorem 7.8,
Remark 7.9,
Propostion 7.10,
and
Corollary 7.11,
extend to this situation
where the role of
$L$ in (7.7) -- (7.11) is played
by $L'$, and the statements
are in fact
compatible with the additional $(A,L'')$-module structures.
Rather than spelling out the details, we confine ourselves
to describing the consequences thereof
for the differential bigraded Batalin-Vilkovisky algebras
studied in earlier sections.
\smallskip
Extending the notion
of dualizing module introduced in Section 2 of \cite\duality,
let 
$$
C_{\Cal L'} = \roman{Alt}^*_A(L'',C_{L'})=
\roman{Alt}^*_A(L'',\roman{Hom}_A(\Lambda^n_A L',A)).
$$
The reasoning in Section 2 of \cite\duality,
adapted to the present differential bigraded setting, 
shows that
$
C_{\Cal L'}
$
is canonically isomorphic to
$$
\roman H^n(\Cal L',U(\Cal A'',\Cal L';d''))
$$
and hence inherits a differential graded right
$U(\Cal A'',\Cal L';d'')$-module structure
from
the obvious right
$U(\Cal A'',\Cal L';d'')$-module structure
on
$U(\Cal A'',\Cal L';d'')$
which remains free when
the construction
of
$\roman H^n(\Cal L',U(\Cal A'',\Cal L';d''))$
is carried out.
The theory of homological duality
developed in \cite\duality\ 
now carries over verbatim and yields
natural isomorphisms
$$
\roman H^k(\Cal L',\Cal M)
\cong
\roman H_{n-k}(\Cal L',C_{\Cal L'} \otimes_{\Cal A''} \Cal M)
\tag7.12.1
$$
for all non-negative integers $k$ and all
left $(\Cal A'',\Cal L')$-modules $\Cal M$ and, furthermore,
natural isomorphisms
$$
\roman H_{k}(\Cal L',\Cal N)
\cong
\roman H^{n-k}(\Cal L,\roman{Hom}_{\Cal A''}(C_{\Cal L'},\Cal N))
\tag7.12.2
$$
for all non-negative integers $k$ and all
right $(\Cal A'',\Cal L')$-modules $\Cal N$.
We therefore refer to
$C_{\Cal L'}$ as the {\it dualizing module\/}
of
$\Cal L'$.
\smallskip
Under the circumstances of (5.4.6),
let $\partial'$ be an exact generator
turning 
\linebreak
$(\roman{Alt}_A(L'',\Lambda_A L'),\partial',\parti'')$
into a differential bigraded Batalin-Vilkovisky algebra.
By Theorem 7.6, this generator $\partial'$ endows
$\Cal A'' = \roman{Alt}_A(L'',A)$
with a right
$U(\Cal A'',\Cal L')$-module structure
and we denote the resulting
right
$U(\Cal A'',\Cal L')$-module
by
$\Cal A''_{\partial'}$.
Since 
$\partial'$
turns
$(\roman{Alt}_A(L'',\Lambda_A L'),\partial',\parti'')$
into a differential bigraded Batalin-Vilkovisky algebra
(not just into a bigraded Batalin-Vilkovisky algebra),
the 
right
$U(\Cal A'',\Cal L')$-module
$\Cal A''_{\partial'}$
is a differential graded
right $U(\Cal A'',\Cal L';d'')$-module.
Inspection shows 
the following:

\proclaim{Proposition 7.13}
The chain complex
$\Cal A''_{\partial'} \otimes _{U(\Cal A'',\Cal L';d'')}K(\Cal A'',\Cal L')$
calculating
$$
\roman H_*(\Cal L',\Cal A''_{\partial'}) 
= \roman{Tor}_*^{U(\Cal A'',\Cal L';d'')}(\Cal A''_{\partial'},\Cal A'')
$$
boils down
to the 
chain complex which underlies the
differential bigraded Batalin-Vilkovisky algebra
$(\roman{Alt}_A(L'',\Lambda_A L'),\partial',\parti'')$
(coming into play in {\rm (5.4.6)}).
Thus the exact generator
$\partial'$
amounts 
to the differential graded Lie-Rinehart
differential in the corresponding
standard complex
$\Cal A''_{\partial'} \otimes _{U(\Cal A'',\Cal L';d'')}K(\Cal A'',\Cal L')$,
with reference to the differential graded right
$U(\Cal A'',\Cal L';d'')$-module structure on
$\Cal A''_{\partial'}$. \qed
\endproclaim

Likewise,
by Corollary 7.10,
the generator $\partial'$ endows
$\Lambda^n_{\Cal A''}\Cal L' =\roman{Alt}^*_A(L'',\Lambda^n_A L')$
with a left
$(\Cal A'',\Cal L')$-module structure,
and we denote the resulting
left $(\Cal A'',\Cal L')$-module
by
$\Lambda^n_{\Cal A''}\Cal L'_{\partial'}$.
Since 
$\partial'$
turns
$(\roman{Alt}_A(L'',\Lambda_A L'),\partial',\parti'')$
into a differential bigraded Batalin-Vilkovisky algebra,
the 
left $(\Cal A'',\Cal L')$-module
$\Lambda^n_{\Cal A''}\Cal L'_{\partial'}$
is  a differential graded
left $(\Cal A'', \Cal L'; d'')$-module, i.~e. a differential graded left
$U(\Cal A'',\Cal L';d'')$-module.
Again inspection shows the following.

\proclaim{Proposition 7.14}
The chain complex
$
\roman{Hom}_{U(\Cal A'',\Cal L';d'')}(K(\Cal A'',\Cal L'),
\Lambda^n_{\Cal A''}\Cal L'_{\partial'})
$
computing
$$
\roman H^*(\Cal L',\Lambda^n_{\Cal A''}\Cal L'_{\partial'}) 
=\roman{Ext}^*_{U(\Cal A'',\Cal L';d'')}(\Cal A'',
\Lambda^n_{\Cal A''}\Cal L'_{\partial'})
$$
comes down to
$(\roman{Alt}_A(L'',\roman{Alt}_A(L',\Lambda_A^nL')),\parti',\parti'')$.
Moreover, 
$\roman{Hom}_{\Cal A''}(C_{\Cal L'},\Cal A''_{\partial'})$
is canonically isomorphic to
$\Lambda^n_{\Cal A''}\Cal L'_{\partial'}$,
and
the isomorphism
$$
(\roman{Alt}_A(L'',\Lambda^*_A L'),\parti',\parti'')
@>>>
(\roman{Alt}_A(L'',\roman{Alt}_A^{n-*}(L',\Lambda_A^nL')),\partial',\parti'')
\tag7.15
$$
of chain complexes 
spelled out as {\rm (5.4.6.1)} 
induces the corresponding
duality isomorphism
$$
\roman H_*(\Cal L',\Cal A''_{\partial'}) 
@>>>
\roman H^{n-*}(\Cal L', 
\roman{Hom}_{\Cal A''}(C_{\Cal L'},\Cal A''_{\partial'})) 
\cong
\roman H^{n-*}(\Cal L',\Lambda^n_{\Cal A''}\Cal L'_{\partial'}) 
\tag7.16
$$
given as {\rm (7.12.2)} above,
where the roles of
$\Cal L,\Cal A,\Cal N$
in {\rm (7.12.2)} are played by,
respectively,
$\Cal L', \Cal A'', \Cal A''_{\partial'}$. \qed
\endproclaim

Thus, in view of the remarks about the Tian-Todorov Lemma
(5.4.8) made in Section 5 above,
this Lemma comes down to differential graded homological duality.
\smallskip
We conclude with the following observation:
When
the twilled Lie-Rinehart algebra $(A,L',L'')$ 
has $L'$ and $L''$ abelian,
with trivial actions of $L'$ and $L''$ on $A$
and on $L''$ and $L'$ (respectively),
the duality isomorphism
(7.16) for this special case
is just the isomorphism (7.15),
the operators $d',d'',\partial'$
being ignored.
Thus,
for a general twilled Lie-Rinehart algebra $(A,L',L'')$,
($L'$, $L''$ still finitely generated and projective as
$A$-modules and $L'$ of constant rank),
the isomorphism 
of bigraded $A$-modules underlying
(7.15) is obtained when the (non-trivial)
true twilled Lie-Rinehart
structure is ignored.
The true twilled
Lie-Rinehart
structure being considered as a \lq\lq perturbation\rq\rq \ 
of the trivial
twilled Lie-Rinehart
structure,
for the duality isomorphism (7.16),
this perturbation
amounts to insertion of the operators 
$d',d'', \partial '$
which, in turn, may be viewed as perturbations of the
trivial operators.

\medskip\noindent {\bf 8. Globalization}\smallskip\noindent
Let $M$ be a smooth manifold, 
let $A$ be the ring   $C^{\infty} M$ of smooth functions on $M$,
and let
$\zeta'$ and $\zeta''$ be 
Lie algebroids over $M$,
that is,
$\zeta'$ and $\zeta''$ are smooth real vector bundles
together with 
$(\Bobb R,A)$-Lie algebra structures on
the spaces
of sections
$L' = \Gamma (\zeta')$ and $L'' =\Gamma (\zeta'')$.
Given a twilled Lie-Rinehart algebra structure turning
$(A,L',L'')$ into a twilled Lie-Rinehart algebra,
we will say that the pair
$(\zeta',\zeta'')$ is a {\it twilled Lie algebroid\/}.
In \cite\mackfift\  and in \cite\mokritwo\ 
these objects are referred to as matched pairs of Lie algebroids.
Likewise,
we can consider the ring $A^{\Bobb C} =C^{\infty} (M,\Bobb C)$ of smooth 
complex functions on $M$
and two complex vector bundles
$\zeta'$ and $\zeta''$;
let
$L' = \Gamma (\zeta')$ and $L'' =\Gamma (\zeta'')$
be their spaces of sections.
Given a twilled Lie-Rinehart algebra structure turning
$(A^{\Bobb C},L',L'')$ into a twilled Lie-Rinehart algebra,
we will say that the pair
$(\zeta',\zeta'')$ is a {\it complex twilled Lie algebroid\/}.
An example of a complex twilled Lie algebroid
arises from a complex structure on $M$.
Another example arises from 
Cauchy-Riemann structures.
\smallskip
Any Lie groupoid $G \text{ \lower 7pt
\vbox to 1.15 pc
 {
  \hbox {$@>{\,\,\, \,\,\,}>>$}
  \vskip-1.2pc
  \hbox {$@>{\,\,\,\,\,\,}>>$}
                }
}
P$ gives rise to a Lie algebroid $AG$.
What is the corresponding 
object for a
twilled Lie algebroid?
To provide an answer to this question, we recall that, by
Theorem 8.3 of \cite\mackxu,
for 
any Poisson groupoid
$G \text{ \lower 7pt
\vbox to 1.15 pc
 {
  \hbox {$@>{\,\,\, \,\,\,}>>$}
  \vskip-1.2pc
  \hbox {$@>{\,\,\,\,\,\,}>>$}
                }
}
P$, 
the pair $(AG,A^*G)$ 
consisting of the Lie algebroid $AG$ and its dual $A^*G$
inherits a 
Lie bialgebroid structure.
Let 
$(\zeta',\zeta'')$
be a twilled Lie algebroid;
in view of
Corollary 4.9,
$(\zeta' \ltime (\zeta'')^*,\zeta''\ltime (\zeta')^*)$
then inherits a Lie bialgebroid structure.
We define a {\it corresponding Lie groupoid\/} to be
a Poisson groupoid
$G \text{ \lower 7pt
\vbox to 1.15 pc
 {
  \hbox {$@>{\,\,\, \,\,\,}>>$}
  \vskip-1.2pc
  \hbox {$@>{\,\,\,\,\,\,}>>$}
                }
}
P$ 
such that
the pair $(AG,A^*G)$ 
is isomorphic to
$(\zeta' \ltime (\zeta'')^*,\zeta''\ltime (\zeta')^*)$
as a
Lie bialgebroid.
Such a Poisson groupoid globalizes
the notion of
twilled Lie algebroid or of
matched pair of Lie algebroids.
What remains to be done is first to single out explicitly those
Poisson groupoids
$G \text{ \lower 7pt
\vbox to 1.15 pc
 {
  \hbox {$@>{\,\,\, \,\,\,}>>$}
  \vskip-1.2pc
  \hbox {$@>{\,\,\,\,\,\,}>>$}
                }
}
P$ 
such that 
the pair $(AG,A^*G)$ 
is of the kind
$(\zeta' \ltime (\zeta'')^*,\zeta''\ltime (\zeta')^*)$,
and thereafter
to give an intrinsic description
of the structure which thus emerges
in terms of groupoids alone.
We hope to return to this at another occasion.
This kind of groupoid might also lead to a concept of groupoid
which integrates a general complex Lie algebroid.
It will certainly integrate those
complex Lie algebroids
$\eta$ which come together with their complex conjugate
$\overline \eta$
in such a way that
$\eta \oplus\overline \eta$
carries the Lie algebroid structure which corresponds
to a twilled sum,
for example those arising from
a complex structure on a smooth manifold or from
a Cauchy-Riemann structure.
See for example \cite\canhawei\ (15.4)
for a discussion of complex Lie algebroids
and how 
a Cauchy-Riemann structure
gives rise to a complex Lie algebroid.

\bigskip
\widestnumber\key{999}
\centerline{References}
\smallskip\noindent

\ref \no \barakont
\by S. Barannikov and M. Kontsevich
\paper Frobenius manifolds and formality of Lie algebras of polyvector fields
\paperinfo alg-geom/9710032
\jour Internat. Res. Notices
\vol 4
\yr 1998
\pages 201--215
\endref

\ref \no \canhawei
\by A. Cannas de Silva, K. Hartshorn, A. Weinstein
\paper Lectures on Geometric Models for Noncommutative Algebras
\paperinfo U of California at Berkeley, June 15, 1998
\endref

\ref \no \cheveile
\by C. Chevalley and S. Eilenberg
\paper Cohomology theory of Lie groups and Lie algebras
\jour  Trans. Amer. Math. Soc.
\vol 63
\yr 1948
\pages 85--124
\endref

\ref \no \evluwein
\by S. Evens, J.-H. Lu, and A. Weinstein
\paper Transverse measures, the modular class, and a cohomology pairing
for Lie algebroids
\jour Quart. J. of Math. (to appear)
\endref

\ref \no \gersthtw
\by M. Gerstenhaber
\paper The cohomology structure of an associative ring
\jour Ann. of Math.
\vol 78
\yr 1963
\pages  267-288
\endref

\ref \no \geschthr
\by M. Gerstenhaber and S. D. Schack
\paper Algebras, bialgebras, quantum groups and algebraic
deformations
\paperinfo In: Deformation theory and quantum groups with
applications to mathematical physics, M. Gerstenhaber and J. Stasheff, eds.
\jour Cont. Math.
\vol 134
\pages 51--92
\publ AMS
\yr 1992
\publaddr Providence 
\endref

\ref \no \getzltwo
\by E. Getzler
\paper Batalin-Vilkovisky algebras and two-dimensional topological field
theories
\jour Comm. in Math. Phys.
\vol 195
\yr 1994
\pages 265--285
\endref

\ref \no \golmilfo
\by W. M. Goldman and J. Millson
\paper The homotopy invariance of the Kuranishi space
\jour Illinois J. of Mathematics
\vol 34
\yr 1990
\pages 337--367
\endref

\ref \no \grothtwo
\by A. Grothendieck
\paper On the de Rham cohomology of algebraic varieties
\jour Pub. Math. Sci. I. H. E. S.
\vol 29
\yr 1966
\pages 351--359
\endref

\ref \no \poiscoho
\by J. Huebschmann
\paper Poisson cohomology and quantization
\jour J. f\"ur die Reine und Angew. Math.
\vol 408
\yr 1990
\pages 57--113
\endref

\ref \no \duality
\by J. Huebschmann
\paper 
Duality for Lie-Rinehart algebras and the modular class
\paperinfo dg-ga/9702008
\jour Journal f\"ur die Reine und Angew. Math. (to appear)
\endref

\ref \no \bv
\by J. Huebschmann
\paper Lie-Rinehart algebras, Gerstenhaber algebras, and Batalin-
Vilkovisky algebras
\jour Annales de l'Institut Fourier
\vol 48
\yr 1998
\pages 425--440
\endref

\ref \no \extensta
\by J. Huebschmann
\paper 
Extensions of Lie-Rinehart algebras and the Chern-Weil construction
\paperinfo in: Festschrift in honor of J. Stasheff's 60th birthday
\jour Cont. Math. (to appear)
\publ Amer. Math. Soc.
\publaddr Providence R. I.
\endref

\ref \no \crosspro
\by J. Huebschmann
\paper Crossed products and twilled Lie-Rinehart algebras
\paperinfo in preparation
\endref

\ref \no \liribi
\by J. Huebschmann
\paper The modular class and master equation for Lie-Rinehart bialgebras
\paperinfo in preparation
\endref

\ref \no \huebstas
\by J. Huebschmann and J. D. Stasheff
\paper Formal solution of the master equation via HPT and
deformation theory
\paperinfo in preparation
\endref

\ref \no \kosmathr
\by Y. Kosmann-Schwarzbach 
\paper Exact Gerstenhaber algebras and Lie bialgebroids
\jour  Acta Applicandae Mathematicae
\vol 41
\yr 1995
\pages 153--165
\endref

\ref \no \kosmafou
\by Y. Kosmann-Schwarzbach 
\paper From Poisson algebras to Gerstenhaber algebras
\jour Annales de l'Institut Fourier
\vol 46
\yr 1996
\pages 1243--1274
\endref

\ref \no \kosmafiv
\by Y. Kosmann-Schwarzbach 
\paper The Lie bialgebroid of a Poisson-Nijenhius manifold
\jour  Letters in Math. Physics
\vol 38
\yr 1996
\pages 421--428
\endref

\ref \no \kosmagtw
\by Y. Kosmann-Schwarzbach and F. Magri 
\paper Poisson-Lie groups and complete integrability. I.
Drinfeld bigebras, dual extensions and their
canonical representations
\jour  Annales Inst. H. Poincar\'e S\'erie A (Physique th\'eorique)
\vol 49
\yr 1988
\pages 433--460
\endref

\ref \no \koszulon
\by J. L. Koszul
\paper Crochet de Schouten-Nijenhuis et cohomologie
\jour Ast\'erisque,
\vol hors-s\'erie,
\yr 1985
\pages 251--271
\paperinfo in E. Cartan et les Math\'ematiciens d'aujourd'hui, 
Lyon, 25--29 Juin, 1984
\endref

\ref \no \liazutwo
\by B. H. Lian and G. J. Zuckerman
\paper New perspectives on the BRST-algebraic structure
of string theory
\jour Comm. in Math. Phys.
\vol 154
\yr 1993
\pages  613--646
\endref

\ref \no \liwextwo
\by Z.-J. Liu, A. Weinstein and P. Xu
\paper  Manin triples for Lie bialgebroids
\jour J. Diff. Geom.
\vol 45
\yr 1997
\pages 547--574
\endref 

\ref \no \liweixu
\by Z.-J. Liu, A. Weinstein and P. Xu
\paper Dirac structures and Poisson homogeneous spaces
\paperinfo preprint
\endref 

\ref \no \luweinst
\by J.-H. Lu and A. Weinstein
\paper Poisson Lie groups, dressing transformations, and Bruhat decompositions
\jour J. of Diff. Geom.
\vol 31
\yr 1990
\pages 501--526
\endref

\ref \no \mackfift
\by K. Mackenzie
\paper Double Lie algebroids and the double of a Lie bialgebroid
\paperinfo preprint 1998; math.DG/9808081
\endref

\ref \no \mackxu
\by K. C. Mackenzie and P.  Xu
\paper Lie bialgebroids and Poisson groupoids
\jour Duke Math. J.
\vol 73
\yr 1994
\pages 415--452
\endref

\ref \no \maclaboo
\by S. Mac Lane
\book Homology
\bookinfo Die Grundlehren der mathematischen Wissenschaften
 No. 114
\publ Springer
\publaddr Berlin $\cdot$ G\"ottingen $\cdot$ Heidelberg
\yr 1963
\endref

\ref \no \majidtwo
\by S. Majid
\paper Matched pairs of Lie groups associated to solutions
of the Yang-Baxter equation
\jour Pac. J. of Math.
\vol 141
\yr 1990
\pages 311--332
\endref

\ref \no \maninfiv
\by Yu. I. Manin
\paper Three constructions of Frobenius manifolds
\paperinfo math.QA/9801006
\endref

\ref \no \mokritwo
\by T. Mokri
\paper Matched pairs of Lie algebroids
\jour Glasgow Math. J.
\vol 39
\yr 1997
\pages 167--181
\endref

\ref \no \palaione
\by R. S. Palais
\paper The cohomology of Lie rings
\jour  Proc. Symp. Pure Math.
\vol III
\yr 1961
\pages 130--137
\paperinfo Amer. Math. Soc., Providence, R. I.
\endref

\ref \no \rinehone
\by G. Rinehart
\paper Differential forms for general commutative algebras
\jour  Trans. Amer. Math. Soc.
\vol 108
\yr 1963
\pages 195--222
\endref

\ref \no \schecone
\by V. Schechtman
\paper Remarks on formal deformations and Batalin-Vilkovisky algebras
\paperinfo math.AG/9802006
\endref

\ref \no \tianone
\by G. Tian
\paper A note on Kaehler manifolds with $c_1=0$
\paperinfo preprint
\endref

\ref \no \todortwo
\by A. N. Todorov
\paper The Weil-Petersson geometry
of the moduli space of $\fra {su}(n \geq 3)$ (Calabi-Yau)
manifolds, I.
\jour Comm. Math. Phys.
\vol 126
\yr 1989
\pages 325--346
\endref

\ref \no  \weinsfte
\by A. Weinstein
\paper The modular automorphism group of a Poisson manifold
\paperinfo in: Special volume in honor of A. Lichnerowicz
\jour J. of Geometry and Physics
\vol 23
\yr 1997
\pages 379--394
\endref

\ref \no \wittetwe
\by E. Witten
\paper Mirror manifolds and topological field theory
\paperinfo in: Essays on mirror manifolds, 
S. T. Yau, ed.
\publ International Press Co.
\publaddr Hong Kong
\yr 1992
\pages 230--310
\endref

\ref \no \xuone
\by P. Xu
\paper 
Gerstenhaber algebras and BV-algebras
in Poisson geometry
\paperinfo preprint, 1997
\endref
\enddocument